%% file: pdf.tex
\newcommand\SEKImasterusepackages{%
  \usepackage[T1]{fontenc}%
  \usepackage[normalem]{ulem}%
  \usepackage{makeidx}%
  \makeindex%
}
\newcommand\SEKIusersusepackages
{
\usepackage{named}
\bibliographystyle{languageindependentcpsnamedwithapalikesorting}
\renewcommand\refname{Bibliography}
}

\newcommand\SEKIREVIEWSTATUS{1}

\newcommand\SEKISUBEDITOR
{{\sc Erica Melis} 
\\FR Informatik, Universit\"at des Saarlandes, D--66123 Saarbr\"ucken, Germany
\\E-mail: {\tt melis@activemath.org}
\\WWW: \url{http://www.activemath.org/~melis}
}

\newcommand\SEKIREVIEWER
{\siegnamenoindex, \CMU, \Dept\ of Philosophy
  \getittotheright{Baker Hall 161, 5000 Forbes Avenue Pittsburgh, PA 15213}
\\E-mail: {\tt sieg@cmu.edu}
\\WWW: \url{https://www.cmu.edu/dietrich/philosophy/people/faculty/sieg.html}
}

\newcommand\SEKIREFEREEPROCESS
{}
\newcommand\SEKIDOCUMENTCLASS{0}
\newcommand\SEKIYEAR{2017}
\newcommand\SEKINUMBEROFYEAR{01}
\newcommand\SEKITYPE{1}
\newcommand\SEKIPUBLISHEDTYPE{0}
\title{A Most Interesting
\mbox{Draft for \hilbertplain} and \bernaysplain' ``\grundlagendermathematiknoindex''
\\\discretionary{\mbox{that never found its way into any publi-}}
                {\mbox{cation, and 2 CV\,of \hskip.1em\hasenjaegernameplain
                }}
                {\mbox{\begin{tabular}{@{}c@{}}%
                         that never found its way into any publication,
                       \\and \hskip.1em two \hskip.15em CV
                       of \hskip.2em\hasenjaegernameplain
                       \end{tabular}}}
}
\author{\wirthnamenoindex\\\Institute\\\emailcp}
\date{\small
  First Published: \
  \Mar\,4, 2018\\
  Thoroughly\,\rev\,\&\,largely\,\extd\,(title,
  \sectrefss
    {section Our Typescript}
    {section Trying to Find Hints on the Time of Writing}
    {section hasenjaeger and bernays},
           CV, \refname, \etc): \ 
  \Jan\,20, 2020\\
  Thoroughly\,\rev\,\&\,largely\,\extd\,(
  all\,sections, now\,all\,texts\,with\,English\,translation): \
  \Jun\,9, 2021\\
  Minor updates to \sectref{section Joint Work} and the Bibliography: \
  \Jun\,25, 2021}
\input seki-deckblatt-5
\input namedhelper

\input headerhot
\input headernames

\input headerforformulas
\input headersugarterms

\input header12pt

\input headerproof

\input quotation

\input hbdictionary

\newcommand\germanfont
{\rm}
\newcommand\germanfontfootnote
{\rm\footnotesize}
\renewcommand\namefont{\sc}
\newcommand\thispaper{this paper}
\newcommand\germanpointertogoedelsincompletness
{\germanfont Tats\ae chlich hat e\es\ sich
im Fall der \germanbeweistheorie\ herau\esi gestellt,
da\sz\ f\ue r die gew\ue nschten Nachweise 
der \germanwiderspruchsfreiheit\ formalisierter Theorien
die \germanfinit en Methoden
nicht zul\ae nglich sind.}
\newcommand\englishpointertogoedelsincompletness
{\mbox{For the case} of \englishbeweistheorie,               \hskip.1em
it has indeed turned out that the \englishfinit\ methods
are not sufficient for the \englishgewuenscht\ proofs of
\englishwiderspruchsfreiheit\ of formalized theories.}
\newcommand\germanpointertohasenjaegersdiss
{Dabei ist in\esi besondere der Umstand mit\-bestimmend,
da\sz\ neuerding\es\ die \germanbeweistheoretisch en Unter\-suchungen
in einen engeren Kontakt getreten sind
mit den allgemeinen Theorien der \mbox{\germanabstrakt en}
Algebra und Topologie, \hskip.3em
\actualpagebreakpagenumber{14}
so da\sz\ die Au\esi sicht sich er\oe ffnet,
da\sz\ die \germanbeweistheoretisch en Methoden
zu einem wirkung\esi vollen Hilf\esi mittel in diesen Gebieten sich entwickeln.}
\newcommand\englishpointertohasenjaegersdiss[1]{%
A particular reason for this non-limitation is 
the fact that recently the \englishbeweistheoretisch\ investigations
have come into closer contact
with the general theories of abstract algebra and topology,#1\       \hskip.3em
\actualpagebreakpagenumber{14}
and thus the prospect is opening up
that the \englishbeweistheoretisch\ methods
may develop into a powerful tool in these fields.}
\newcommand\selfquotedsentenceattheendofsectiononeofourtypescript
{In view of this \englishsachlage, \hskip.1em
the appropriate \englishverfahrenprocedure\ seems to be  that
we will \mbox{keep track} of the
\englishgesichtspunkt\ of \englisheinfinit\ \englishbetrachtungbehandlung\
as a guideline to determine our method of choice, \hskip.2em
but do not strictly limit ourselves to this method.}
\newcommand\textongermanoriginalanditsenglishtranslationrectoverso{%
The German original is found on the pages with even page numbers and
our English translation on the respectively following pages with odd
page numbers.\par
To this end, the remainder of the current page is left blank here.}

\catcode`\@=11
\renewcommand\tableofcontents{%
    \section*{\contentsname
        \@mkboth{%
           \MakeUppercase\contentsname}{\MakeUppercase\contentsname}}%
    \vskip .1ex
    \@starttoc{toc}%
    }
\def\l@section#1#2{%
  \ifnum \c@tocdepth >\z@
    \addpenalty\@secpenalty
    \addvspace{1.2em \@plus\p@}
    \setlength\@tempdima{2.5em}
    \begingroup
      \parskip -1pt            
      \parindent \z@ \rightskip \@pnumwidth
      \parfillskip -\@pnumwidth
      \leavevmode
      \advance\leftskip\@tempdima
      \hskip -\leftskip
      \bf#1\nobreak\hfil \nobreak\hb@xt@\@pnumwidth{\hss #2}\par
    \endgroup
    \fi}
\def\@dottedtocline#1#2#3#4#5{%
  \ifnum #1>\c@tocdepth \else
    \vskip \z@ \@plus.2\p@
    {\leftskip #2\relax \rightskip \@tocrmarg \parfillskip -\rightskip
     \parindent #2\relax\@afterindenttrue
     \interlinepenalty\@M
     \leavevmode
     \@tempdima #3\relax
     \advance\leftskip \@tempdima \null\nobreak\hskip -\leftskip
     {#4}\nobreak
     \leaders\hbox{$\m@th
        \mkern \@dotsep mu\hbox{.}\mkern \@dotsep
        mu$}\hfill
     \nobreak
     \hb@xt@\@pnumwidth{\hfil\normalfont \normalcolor #5}%
     \notop\halftop  
     \par
     }%
  \fi}
\def\l@subsection{\@dottedtocline{2}{2.5em}{3em}}
\def\l@subsubsection{\@dottedtocline{3}{3.8em}{3.2em}}
\catcode`\@=12
\begin{document}
\makecover
\maketitle
\begin{abstract}
  In\,1934, in \bernays' preface to the {\em first edition}\/ 
  of the first volume of
  \hilbert\ and \bernays' monograph ``\grundlagendermathematik\closequotecomma
  a nearly completed draft of the the finally two-volume monograph is
  mentioned, which had to be revoked because of the completely
  changed situation in the area of \englishbeweistheorie\ after
  \herbrand\ and \goedel's revolutionary results. 
  Nothing \atall\ seems to be known about this draft and its whereabouts. \par
  A third of a century later,
  \bernays' preface to the {\em second edition}\/ (1968)
  of the first volume of
  \hilbert\ and \bernays' ``\grundlagendermathematik'' mentions
  joint work\linebreak of \hasenjaeger\ and \bernays\ on the second edition.
  \hskip.2em \bernays\
  states there that\linebreak 
  ``it became obvious that the integration of the many new results in the area 
  of \englishbeweistheorie\ would have required
  a \mbox{complete} re\-organization of the book\closequotecommaextraspace
  \ie\ that the inclusion of the 
  intermediately found new results in the area of \mbox\englishbeweistheorie\
  \mbox{turned out} to be unobtainable by a revision,
  but would have required
  {\em a complete reorganization of the \mbox{entire textbook}}.
  \hskip.3em
  \mbox{We will document that}
  \mbox{---~even after}
  the need for a complete reorganization had become obvious~---
  this joint work went \nolinebreak on
  \mbox{\em to a considerable extent}. \hskip.3em
  Moreover, \hskip.1em
  we will document {\em when \hasenjaegerplain\ stayed in \Zuerich}\/
  to \nolinebreak assist \bernays\
  in \nolinebreak the completion of the second edition.\par
  In May 2017, \hskip.1em
  we identified an incorrectly filed text 
  in \bernays' scientific \nolinebreak legacy at the archive of the \ETHZshort\
  as a candidate for the beginning of \nolinebreak
  the \nolinebreak revoked \nolinebreak draft for the first edition or
  of a revoked draft for the second edition.
  In a {\em\mbox{partial presentation} and
  careful investigation} of this text we gather only some minor evidence
  that this \nolinebreak text is the beginning of
  the nearly completed draft of the
  first edition, 
  but ample evidence that
  this text is part of the work of \hasenjaeger\ and \bernays\
  on the second \nolinebreak edition.
  We provide some evidence that this work
  has covered
  {\em a complete reorganization of the entire first volume,
       including a completely new version of its last chapter on the \math\iota}.
  \end{abstract}

\vfill\pagebreak

\tableofcontents\vfill\cleardoublepage

\section{Introduction}
\yestop\subsection{\hilbertbernaysplain\ in general}
By the 1930s, ground-breaking work had been achieved by German scientists,
especially in philosophy, psychology, 
physics, 
chemistry,
and mathematics. 
With the Nazis' seizure of power in 1933,
the historical tradition of German research
was discontinued in most areas,
and, as a further consequence,
many achievements of German science 
in the first half of the \nth{20}\,century
have still not been
sufficiently recognized.
This is the case especially for those developments
that had not been completed before
the Nazis covered Germany under twelve years 
of intellectual darkness. \hskip.3em
\begin{sloppypar}
\hilbertname\ \hilbertlifetime\
is one of the most outstanding representatives of
mathematics, mathematical physics, and logic-oriented foundational sciences
in general \cite{reid-hilbert}. \hskip.3em
From the end of the \nth{19}\,century to the
erosion of the University of \Goettingen\
by the Nazis, \hskip.2em
\hilbertindex\hilbert\ formed and reshaped many areas of 
applied and pure mathematics. \hskip.3em
Most well-known and highly acknowledged are his
\mbox{``\englishgrundlagendergeometrie''\,\,\cite{grundlagen-der-geometrie}}.

After initial work at the very beginning of the \nth{20}\,century, \hskip.2em
\hilbertindex\hilbert\ re-intensified his 
research into the logical foundations of
mathematics in\,\,1917, \hskip.2em
together with his new assis\-tant
\bernaysname\ \bernayslifetime. \hskip.4em
Supported by their \PhD\ student
\ackermannname\ \ackermannlifetime, \hskip.25em
\bernaysindex\bernays\  and \hilbertindex\hilbert\
developed the field of {\em proof theory}\/
\mbox{(or {\em metamathematics}\/)},
where formalized mathematical proofs 
become themselves the objects
of mathematical operations and investigations ---
just as numbers are the objects of number theory. \hskip.3em
The goal of \hilbertindex\hilbert's endeavors in this field
was to prove the consistency of the customary methods in mathematics
once and for all,
without the loss of essential theorems as in the competing
\englishintuitionistisch\ movements of
\kroneckerindex\kronecker, \brouwerindex\brouwer, \weylindex\weyl, 
and \heytingindex\heyting. \hskip.3em
The proof of the consistency of mathematics 
was to be achieved by sub-division into the following three tasks:
\begin{itemize}\noitem\item
{\em Arithmetization}\/ of mathematics.\noitem\item
{\em Logical formalization}\/ of arithmetic.\noitem\item
{\em Consistency proof}\/ in the form of a {\em proof of impossibility}\/:
\hskip.2em It cannot occur in arithmetic
that there are
formal derivations of 
a formula \nlbmath A and also of its negation 
\nlbmaths{\mediumheadroom\overline A}.
\noitem\end{itemize}
The problematic step in this \programme\
(nowadays called {\em\hilbertsprogram}\/) \hskip.1em
is the consistency proof.\end{sloppypar}

\hilbertsprogram\ was nourished by the hope that mathematics
---~as the foundation of natural sciences, and especially of modern physics~---
could thus provide the proof of its own groundedness.
This was a paramount task of the time,
not least because of the foundational crisis
in mathematics
(which had been evoked among others by \russellsparadox\
 at the beginning of the \nth{20}\,century) \hskip.2em
and the vivid philosophic discussions of the formal sciences
stimulated {\it inter alia}\/ by
\wittgensteinindex\wittgenstein\
and
the Vienna Circle.\footnote{
\Cfnlb\ \eg\ 
\cite{wittgenstein-wiener-ausgabe-1},
\cite{wittgenstein-wiener-kreis-mcguinness}.%
}\pagebreak

It should be recognized that \hilbertindex\hilbert's primary goal was neither 
a reduction of mathematical reasoning and writing to formal
logic as in the seminal work of \citet{PM}, \hskip.2em
nor a formalization of larger parts of mathematics
as in the publications of the famous French 
\citet{bourbaki} \hskip.1em group of mathematicians. \hskip.5em
His ambition was to secure ---~once and for all~---
the foundation of mathematics with consistency proofs,
in \nolinebreak which an intuitively consistent, ``\englishfinit'' part
of mathematics was 
to be
used for showing 
that no contradiction 
could be formally derived in larger and larger parts of 
non-constructive and axiomatic mathematics.

\hilbertsprogram\ fascinated an elite of young outstanding mathematicians,
among them
\neumannname\ \neumannlifetime,   \hskip.2em
\goedelname\ \goedellifetime,     \hskip.2em
\herbrandname\ \herbrandlifetime, \hskip.2em
and
\gentzenname\ \gentzenlifetime,   \hskip.2em
whose contributions
essentially shaped the fields of 
modern mathematical logic
and proof theory.

We know today that
\mbox{\hilbertindex\hilbert's} quest
to establish a foundation for the whole 
scientific edifice could not be successful to the proposed extent: \hskip.2em
\goedel's incompleteness theorems dashed the broader hopes of \hilbertsprogram.
Without the emphasis that \hilbertindex\hilbert\
has put on the foundational issues, however,
our negative and positive knowledge
on the possibility of a logical grounding
of mathematics (and thus of all exact sciences) would hardly have been achieved
at his time.
\yestop\halftop\halftop\subsection{Drafts for \hilbertbernays}
The central and most involved presentation of
\hilbertsprogram\ and \hilbertindex\hilbert's proof theory
is found in the two-volume monograph ``\grundlagendermathematik''
of \makeacitetoftwo
{grundlagen-first-edition-volume-one}
{grundlagen-first-edition-volume-two},
and its second revised edition 
\makeaciteoftwo
{grundlagen-second-edition-volume-one}
{grundlagen-second-edition-volume-two}.

We should not forget the historical context of the original writing of
these texts in the late 1920s and the 1930s: \hskip.5em
First,
\herbrandsfundamentaltheorem\
and \goedelsincompletenesstheorems\ hit the new field of
\englishbeweistheorie\ like a hurricane. \hskip.3em
Moreover,
after the Nazi takeover of Germany in January\,1933, \hskip.2em
\bernaysindex\bernays\ was expelled from his academic position in \Goettingen\ in \Apr\,1933,
\hskip.1em and had to \nolinebreak leave the \nolinebreak country. \hskip.3em
As a consequence,
both volumes show strong signs of reorganization and rewriting,
sometimes even the signs of hurry to meet the publication deadlines.

Because of this editorial and historical context,
the drafts for both editions of the two \hilbertbernays\
volumes are of special interest,
in particular the drafts that did not find their way into any of the editions.

Until recently,
however,
we did not know anything substantial about these revoked drafts
for \hilbertbernays. \hskip.3em
In \May~2017, however,
we identified an incorrectly filed text in \bernaysindex\bernays' scientific legacy
in the archive of the \ETHZshort\ as 
a candidate for such a revoked draft. \hskip.2em
We will describe and investigate this text in \thispaper,
and discuss at what time it was probably written, and by whom, \etc

There are two infamous revoked drafts for \hilbertbernays\
of which \bernaysindex\bernays\ stated that they
at least had existed
---~one for the first and one for the second edition:\vfill\pagebreak
\subsubsection{The Mentioning of the Draft for the First Edition}%
\halftop\halftop\noindent           
\bernaysindex\bernays' ``{\germanfont Vorwort zur er\esi ten Au\fli age}''
(``Preface to the First Edition'') \hskip.1em
of \cite[\p VII\,\f]{grundlagen-first-edition-volume-one}, \hskip.2em
begins as follows:\begin{quote}
``{\germanfont Eine Darstellung der \germanbeweistheorie,
welche au\es\ dem \hilbertindex\hilbert schen Ansatz zur Behandlung der mathematisch-logischen
Grundlagen\-probleme erwachsen ist,
wurde schon seit l\ae ngerem von \hilbertindex\hilbert\ angek\ue ndigt.\par
Die Au\esi f\ue hrung diese\es\ Vorhaben\es\ hat eine wesentliche
Verz\oe gerung dadurch erfahren,
da\sz\ in einem Stadium,
in dem die Darstellung schon ihrem Abschlu\sz\ nahe war,
durch da\es\ Erscheinen der Arbeiten von \herbrand\ und \goedel\
eine ver\ae nderte Situation im Gebiet der \germanbeweistheorie\ entstand,
welche die Ber\ue cksichtigung neuer Einsichten \actualpagebreak
zur Au\fgi abe machte.
Dabei ist der Umfang de\es\ Buche\es\ angewachsen,
so da\sz\ eine Teilung
in zwei B\ae nde angezeigt erschien.}''\end{quote}
In the translation of \hskip.1em
\cite[\p VII.b,\,engl.]{grundlagen-german-english-edition-volume-one-one}
\hskip.05em
(comments omitted):\begin{quote}
``Some time ago,
\hilbertindex\hilbert\ announced a presentation of the \englishbeweistheorie\
that developed from the \hilbertindex\hilbert ian
approach to the problems in the foundations of mathematics and logic.
\par
The execution of this enterprise received considerable delay because 
the whole field of \englishbeweistheorie\ was 
changed by the publication of the works of \herbrand\ and \goedel\
when our work was already close to completion; \hskip.1em
and this change put the \englishhinzuziehendebetrachtung\ of new insights 
\actualpagebreak onto the agenda. \hskip.2em
As a consequence of \nolinebreak this, 
the size of the book grew to the extent
that a separation into two volumes seemed appropriate.\end{quote}
Nothing 
about this work
``already close to completion''
and its whereabouts seems to be known.
\yestop
\subsubsection{The Mentioning of the Given-Up Work on the Second Edition}\label
{section The Mentioning of the Given-Up Work on the Second Edition}%
\halftop\halftop\noindent
\bernaysindex\bernays' ``{\germanfont Vorwort zur zweiten Au\fli age}''
(``Preface to the Second Edition'') \hskip.1em
of \hskip.1em\cite[\p V]{grundlagen-second-edition-volume-one} \hskip.2em
begins as follows:\begin{quote}\sloppy
``{\germanfont Schon vor etlichen Jahren haben der verstorbene
\scholzname\ \hskip.1em und Herr \mbox\fkschmidtshortname\ \hskip.1em
mir vorgeschlagen,
eine zweite Au\fli age der ``\grundlagendermathematik'' vorzunehmen,
und Herr \hasenjaegerindex\hasenjaegershortname\ \hskip.1em
war auch zu meiner Unter\-st\ue tzung bei dieser
Arbeit auf einige Zeit nach \Zuerich\ gekommen. \hskip.2em
E\es\ \nolinebreak zeigte sich jedoch bereit\es\ damal\es,
da\sz\ eine Einarbeitung der vielen im Gebiet der \germanbeweistheorie\
hinzugekommenen Ergebnisse eine v\oe llige Umgestaltung 
de\es\ Buche\es\ erfordert h\ae tte. \hskip.2em
Erst recht kann bei der jetzt vorliegenden zweiten Au\fli age,
zu der wiederum Herr \fkschmidtshortname\ den Ansto\sz\ gab,
nicht davon die Rede sein,
den Inhalt dessen,
wa\es\ seither in der \germanbeweistheorie\ erreicht worden ist,
zur Darstellung zu bringen.%
}''\end{quote}
\vfill\pagebreak
\par\noindent
In the translation of \hskip.1em
\cite[\p V\hspace*{-.1em},\,engl.]{grundlagen-german-english-edition-volume-one-one}:\begin{quote}
  ``A number of years ago,
the late \scholzname\ and \MrUS\shorth\fkschmidtshortname\ 
suggested the undertaking of a second edition of the
``Foundations of Mathematics''; \ 
and moreover, 
to assist me in this work,
\MrUS\shorth\hasenjaegerindex\hasenjaegershortname\ came to \Zuerich\ for some time. \hskip.3em
Already back then, 
it became obvious that the integration of the many new results in the area 
of \englishbeweistheorie\ would have required
a complete reorganization of the book. \hskip.3em
Furthermore,
the 
\englishjetztvorliegendpresent\
second edition 
(the impetus 
for which
came again from \MrUS\nolinebreak\shorth\nolinebreak\fkschmidtshortname) 
\hskip.1em
can by no means 
present the 
substance 
of the achievements in \englishbeweistheorie\ since the
appearance of the first edition.''\end{quote}
In \sectref{section Joint Work}, \hskip.2em
we will document when \hasenjaegername\ came to \Zuerich\ for this purpose,
and provide some evidence that the work of \hasenjaegerindex\hasenjaeger\ and \bernaysindex\bernays\ has covered
\mbox{\em a complete reorganization of the entire first volume}, \hskip.2em
which means that the main document
of this collaboration and its whereabouts
still remain completely unknown, \hskip.1em
and also that \bernaysindex\bernays\ may give
us a wrong impression by using the subjunctive
``erfordert h\ae tte'' (``would have required'')
in the above quotation.\vfill\pagebreak
\section{Our Typescript}\label
{section Our Typescript}
The unpublished typescript we discuss in \thispaper\ will briefly be called
``our typescript\closequotefullstopnospace
\subsection{Form, Location, and Incorrect Filing of Our Typescript}\label
{section bernays-no-typewriter}%
Its {\em outer form}\/ is as follows:
\hskip.3em
It is an untitled typescript,
with corrections by \bernaysindex\bernays' hand. \hskip.3em
Our typescript has 34~pages,
with page numbers 2--34 on the respective page headers. \hskip.3em
Its spelling is German--Austrian,
and it includes the German letter ``\sz'' 
not found on typewriters with Swiss layout. 
\hskip.3em
In this context it may be relevant that
---~according \bernaysludwigname\ \cite{bernays-no-typewriter}~---
his close uncle and legator \bernaysname\
never had a typewriter and is not known to have ever used one.

The {\em location}\/ of our typescript is the archive
(\ETHZofficialarchivereference) of the
\ETHZshort\ (\englishETHZlongwithoutshort) (Switzerland). \hskip.3em

The folder in which our typescript was found there by \wirthname\
on \May\,12, 2017, \hskip.1em
was one of two folders
in the legacy of \bernaysname\ under the label
``Hs\,973:\,41\closequotecommasmallextraspace
which is listed in the inventory \cite{bernays-1986} \hskip.1em
on page~7 as
\notop\begin{quote}``41.~~Texte und Korrekturen zur Neuau\fli age de\es\
  ``Grundlagenbuche\es'' Bd.\,II von \nolinebreak \hilbertnameplain\ und
  \bernaysindex\bernaysnameplain. 1970
  \getittotheright{2~Mappen}''\notop\end{quote}
\noindent In English:
\notop\begin{quote}
  ``41.~~Texts and corrections for the new edition of the ``Foundations Book''
  \Vol\,II by \hilbertnameplain\ und \bernaysindex\bernaysnameplain. 1970
  \getittotheright{2~folders}''\notop\halftop\end{quote}
This is a wrong place for our typescript because its contents
refer exclusively to \Vol\,I, \hskip.2em
published in 1934~(\nth 1\,\edn)\ and in 1968~(\nth 2\,\edn)\
\ ---~~neither to \Vol\,II, \hskip.2em nor to the year\,1970\@. \hskip.3em
Thus, regarding its subject,
our typescript was {\em incorrectly filed}, and probably still is.
\subsection{Two Additional Copies of \hskip.1em Our Typescript, With Footnotes!}\label
{section carbon copies}\label
{section additional footnotes}%
Moreover,
in January~2018,
in the legacy of \hasenjaegername\commanospace\footnote{%
  By the end of the year~2018,
  \beatebeckername\ gave her part of the
  scientific legacy of her father \hasenjaegername\ in 9~boxes
  to the following archive:
  \hasenjaegerlegacy.%
}
his daughter \beatebeckername\
succeeded 
in finding
two carbon copies of our typescript.

Both of these carbon copies come without
the corrections {\em by \bernaysindex\bernays' hand}\/ found in our typescript, \hskip.1em
but one of them has most of these corrections added
(mostly with a typewriter), \hskip.1em
and the other includes the typewriting of
even more of the corrections by \bernaysindex\bernays' hand found in our typescript.

One of these carbon copies
comes with two extra typewritten pages containing all footnotes
for our typescript,
again corrected by \bernaysindex\bernays' hand. \hskip.3em
Our typescript, however, does not contain any footnotes at all,
but only raised closing parentheses
(without any numbers, letters, asterisks or other signs)
to indicate their respective positions
in the text.

Regarding the font and the actual typing of the letters,
we did not find any differences between our typescript
and the extra footnote pages found by \beatebeckername.\pagebreak
\subsection{Presentational Form of Our Excerpts
}
\halftop\halftop\noindent
The deletions and additions by \bernaysindex\bernays' handwritten remarks
are not documented in the following excerpts from our typescript; \hskip.2em
only the version of our typescript after
the application of \bernaysindex\bernays' hand-written corrections
is presented here. \par
Therefore,
the rare deletions (\xout{this text is deleted}) \hskip.1em
and additions (\edcomment{this text is added}) \hskip.1em
indica\-ted in these excerpts
are ours, not \bernaysindex\bernays\fullstopnospace'

Moreover ---~in these excerpts~---
the footnotes of the form \edcomment\ldots\
are our additions, whereas the other footnotes
are those from the two extra pages found with one of the two carbon copies
of our typescript in the \hasenjaegerindex\hasenjaeger\ legacy
(\cfnlb\ \sectref{section carbon copies}). \par
All footnote numbers in these excerpts are introduced by us. \hskip.3em
This form of presentation seems appropriate
in \nolinebreak particular because our typescript has no footnotes
(but only raised closing parentheses
to indicate their positions, \cfnlb\ \sectref{section carbon copies}). \par
The positions of the original page breaks in our typescript are
indicated by the sign~``\actualpagebreaknospace'' in our excerpts,
with the number of the new page in a lower index, \ie\
``\actualpagebreaknospacepagenumber 5'' for the page break from page~4 to page~5.
\yestop\yestop\yestop\yestop
\subsection{Discussion and Excerpts
            of \hskip.2em\litsectref 1 of \hskip.1em Our Typescript}\label
{section Discussion and Excerpts of 1}
\halftop\halftop\noindent
The first section of our typescript has the headline
\mbox{``{\germanfont Einleitung}.''} (``Introduction.'') \hskip.1em
and the additional subsection headline \hskip.1em
``{\germanfont\litsectref 1.~Einf\ue hrung in die Fragestellung}''
(``\litsectref 1.~Introduction to the \englishFrageStellungdefinition'').%
\yestop\yestop\subsubsection
{Parallel introductory part, several paragraphs missing in our typescript}\label
{section Parallel introductory part}%
\halftop\halftop\noindent
The text starts literally with the first paragraph of \hskip.15em
\litsectref 1 of
\makeaciteoftwo{grundlagen-first-edition-volume-one}{grundlagen-second-edition-volume-one} (\p\,1), \hskip.3em
including the enumerated list of three items.

The next paragraph of \makeaciteoftwo{grundlagen-first-edition-volume-one}{grundlagen-second-edition-volume-one}, \hskip.2em
however, \hskip.1em
is not present in our typescript, \hskip.1em
namely the paragraph on the
``{\germanfont\germanverschaerft en methodischen Anforderungen}''
(``\englishverschaerft\ \englishmethodischeanforderung{}'')
resulting in
``{\germanfont eine neue Art der Au\esi einander\-setzung mit dem Problem
  de\es\ Unendlichen}''
(``to deal with the problem of the infinite in a new way'').

Then, however, the text continues almost literally with the
penultimate paragraph and the first sentence of the last paragraph of \hskip.1em
\p\,1
of \makeaciteoftwo{grundlagen-first-edition-volume-one}{grundlagen-second-edition-volume-one}. \hskip.5em
The part on the subject of the
``\hspace*{-.2em}{\em exi\esi tentiale Form}''
(``\hspace*{-.2em}{\em existential form}'')
(which starts in the middle of the last paragraph of \hskip.1em \p\,1 and
runs up to the end of the \nth 1\,paragraph starting on \p\,2) \hskip.1em
\mbox{is again missing} in our typescript
--~just as the introduction to this subject was omitted before,
namely in the paragraph on ``the problem of the infinite\closequotecomma
right after the enumerated list of three items. \hskip.3em
Instead of this part, the paragraph ends in our typescript with a
digression into \englisherfuellbarkeit\ and \englishwiderspruchsfreiheit.
\pagebreak

\begin{sloppypar}Then the last paragraph on \p\,2 of our typescript
follows almost literally the text from the
\nth 2\,paragraph starting on \p\,2 of \hskip.1em
\makeaciteoftwo{grundlagen-first-edition-volume-one}{grundlagen-second-edition-volume-one} \hskip.1em
up to the end of the \nth 2\,paragraph starting on \p\,3 of
\hskip.1em
\makeaciteoftwo{grundlagen-first-edition-volume-one}{grundlagen-second-edition-volume-one}, \hskip.2em
where only the last part of the last sentence%
\notop\begin{quote}
  ``{\germanfont al\es\ g\ue ltig vorau\esi gesetzt,
  und wir kommen so zu der Frage,
  welcher Art diese \germangeltung\ \nolinebreak ist.}''\par
  ``presupposed to be valid.
  And so we come to the \englishfragequestion\ what nature of
  \englishgeltung\ this is.''  
\notop\end{quote}
is missing in our typescript, \hskip.1em
where the sentence ends with\notop\begin{quote}
``{\germanfont zugrundegelegt.}''\par
``taken as a basis.''\notop\end{quote}
instead, \hskip.2em
at the {\em very end of page~4}\/ of our typescript.\end{sloppypar}

\halftop\yestop\subsubsection{Over 16~pages entirely missing in the typescript}
\halftop\halftop\noindent
The remainder of \hskip.1em \litsectref 1 of \hskip.1em
\makeaciteoftwo{grundlagen-first-edition-volume-one}{grundlagen-second-edition-volume-one} \hskip.1em
is entirely omitted in our typescript, \hskip.2em
\ie\ everything from the penultimate paragraph on \p\,3 \hskip.1em
to the very end of\,~\p\,19 \hskip.1em
is entirely missing. \hskip.3em
This remainder covers --~ among others subjects~-- the following:
the logical symbolism, 
axiomatizations of geometry, 
\englisherfuellbarkeit\ and \englishallgemeingueltigkeit\ of \formulae, 
the \englishzenoscheparadoxie, the existence of an infinite manifold,
the method of arithmetization, and the task of a proof of
\englishwiderspruchsfreiheit\ as \englisheinunmoeglichkeitsbeweis.

\halftop\yestop\subsubsection{A completely unknown introduction; \hskip.1em
from our typescript, \hskip.1em\PP 5{14}}\label
{section A completely new introduction}%
\halftop\halftop\noindent
In \litsectref 1 of our typescript, however,
right after the very end of page~4
(just mentioned at the very end of
 \sectref{section Parallel introductory part}), \hskip.25em
there is a completely unknown, most interesting, and well-written
introduction to
the foundations of mathematics, \hilbertindex\hilbert's proof theory,
and the \englishfinit\ \englishstandpunkt,
which seems to be entirely un\-published up to now.

We now present this introduction in the following way: \hskip.2em
\textongermanoriginalanditsenglishtranslationrectoverso



\vfill\pagebreak

{\germanfont

\yestop\actualpagebreakindentpagenumber 5
Wir kommen somit zu dem zweiten der anfang\es\ genannten Themata der
Grundlagen\-untersuchungen.
In der Begr\ue ndung der \germananalysis\ ist e\es\ ja
im 19.\,Jahrhundert zuerst
durch die Unter\-suchungen von \bolzanoindex\bolzano\ und \cauchyindex\cauchy\
und hernach deren Weiter\-f\ue hrung und Vollendung durch
\dedekindindex\dedekind, \cantorindex\cantor\ und \weierstrassindex\weierstrass\
gelungen,\footnote{%
  \merayindex\meray, Enc.\,I.1, 3 (1904).%
}
\hskip.1em
die Methoden der Infinitesimal\-rechnung,
die ja in ihren Anf\ae ngen einer vollen Deutlich\-keit entbehrten und
mehr nur instinktiv gehandhabt wurden,
im Sinne einer st\ae rkeren Ankn\ue pfung an die klassischen Methoden
der griechischen Mathematiker \eudoxosindex\mbox\eudoxos\ und
\archimedesindex\archimedes\
zu \xout{einer} pr\ae zise \germanmitteilbar en und lehrbaren zu gestalten.

\yestop\indent
Indem diese Deutlichkeit erreicht wurde,
traten zugleich die zugrunde\-liegenden methodischen \germanvoraussetzung en
mehr hervor,
und man ging auch dazu \ue ber,
diese \germanvoraussetzung en \ue ber die Zielsetzung
der Infini\-tesi\-mal\-rechnung hinau\es\ sy\esi tematisch zu verwerten,
wie e\es\ ja vor allem in der \cantorindex\cantor'schen Mengen\-lehre geschah.
Die hier stattfindende starke \Ue berschreitung de\es\ mathematisch
Gewohnten weckte vielerseit\es\ Kritik,
die dann noch durch die Entdeckung der mengentheoretischen Paradoxien
best\ae rkt wurde.\par

\yestop\indent
Wenngleich e\es\ sich nun auch bei n\ae herem Zusehen erwie\es,
da\sz\ e\es\ zur Verh\ue tung der Paradoxien gen\ue gte,
gewisse extreme \germanbegriffsbildungen\ zu vermeiden,
die tats\ae chlich f\ue r den \germanaufbau\ der Mengen\-lehre 
und erst recht f\ue r die Methoden der \germananalysis\
gar nicht erforderlich sind,
so ist doch seitdem die Di\esi kussion \ue ber die Grundlagen der Mathematik
nicht zur Ruhe gekommen,
und man hat sich auch jener Paradoxien al\es\ Argument bedient,
um viel weiter gehende Einschr\ae nkungen de\es\ mathematischen
\germanverfahren\es\ zu motivieren,
al\es\ sie zur Behebung der \germanwidersprueche\ \actualpagebreakpagenumber 6
direkt erfordert werden.\par

\yestop\indent
F\ue r eine gr\ue ndliche Stellungnahme zu dieser
Grundlagen\-di\esi kussion
erscheint eine ein\-gehen\-de \germanbetrachtung\ der
{\em logischen Struktur der mathematischen Theorien}\/
al\es\ geboten.\par

\yestop\indent
In der Tat bemerkt man,
da\sz\ e\es\ sich bei den zur Di\esi kussion stehenden \germanverfahren\
der Mathematik um Methoden de\es\ Folgern\es\ und der
\germanbegriffsbildung\ handelt,
da\sz\ also hier eine Art der Erweiterung der gew\oe hnlichen Logik
zur \germangeltung\ kommt.
Zugleich zeigt sich eine enge Verflochtenheit de\es\ Mathematischen 
mit dem Logischen:
einer\-seit\es\ tritt die Mengen\-lehre ihrem Gegenstand nach,
durch die \germanbeziehung\ von Mengen und Pr\ae dikaten
(\dasheisst\ durch da\es\ \germanverhaeltnis\ von Umfang und Inhalt der Begriffe)
in eng\esi te Ber\ue hrung mit der Logik;
andererseit\es\ wird man in der sy\esi tematischen Untersuchung
der logischen Bildung\esi formen und \germanschlussweisenoldspelling\
mit Notwendigkeit auf mathematische \germanbetrachtungen\ gef\ue hrt.
So ist ja bereit\es\ die traditionelle Lehre von den kategorischen
Schl\ue ssen eine typisch mathematische Untersuchung,
wa\es\ nur durch ihre hi\esi torische Einordnung in die Philosophie leicht
verdeckt wird.
Mit dieser mathematischen Seite de\es\ Logischen h\ae ngt e\es\ auch zusammen,
da\sz\ die logischen Schl\ue sse
--~wie sie in\esi besondere bei der reichhaltigeren Anwendung der Logik
in den mathematischen Theorien zur Verwendung kommen~--,
in einer mathematischen Weise fixierbar und au\es\ einer Reihe von
wenigen Elementar\-prozessen zusammensetzbar sind.\par

}

\vfill\pagebreak

\yestop\actualpagebreakindentpagenumber 5
Let us now turn to the second topic of the foundational investigations
listed at the beginning,
regarding the \englishbegruendungtheorie\ of \englishanalysis. \hskip.3em
As is well known, in the \nth{19}\,century,
first through the investigations of \bolzanoindex\bolzano\ and
\cauchyindex\cauchy\ and
then through their continuation and completion by \dedekindindex\dedekind,
\cantorindex\cantor\ and \weierstrassindex\weierstrass,\footnote{%
  \meray, Enc.\,I.1, 3 (1904). \par
  \edcomment{%
    It is difficult to say what this abbreviated footnote text by \bernaysindex\bernays' hand means,
    but it probably refers 
    to section ``6.~Point de vue de \merayshortname\closequotecomma\PP{147}{149}
    in \cite{pringsheim-1904}. \hskip.2em
    This would make sense because our typescript omits the
    first name to be mentioned in the list, namely \merayname\ \meraylifetime,
    \hskip.2em
    and \cite{pringsheim-1904} is one of the first texts that mentions the primacy of
    this neglected mathematician. \hskip.3em
    Note that the German original \cite{pringsheim-1898} does not contain this section on
    \meray\ at all.\par
    The footnote is not executed in any of the two carbon copies
    (\cf\ \sectref{section carbon copies})
    of our typescript, \hskip.1em
    but \bernaysindex\bernays\ probably wanted to add ``\meray,'' before ``\dedekind'' and have
    the expanded proper citation in a footnote.%
  }%
} \hskip.1em
the methods of infinitesimal calculus,
which lacked clarity 
and were applied more or less
instinctively in their beginnings,
were given a precisely \englishmitteilbar\ and teachable form,
in the sense of a closer orientation toward the classical methods
of the Greek mathematicians \eudoxosindex\eudoxos\ and
\archimedesindex\archimedes.

\yestop\indent
As this clarity was achieved,
the underlying \englishvoraussetzungen\ of \englishmethodischthesemethods\
became more obvious \aswell.
Moreover, 
these \englishvoraussetzungen\ were then systematically applied
beyond the field of infinitesimal calculus,                            \hskip.1em
in particular in \cantorindex\cantor's set theory.                     \hskip.1em
The strong \englishueberschreitung\ of the common usage of these
\englishvoraussetzungen\ in mathematics aroused \mbox{criticism}
which came from many sides and was then further encouraged
by the discovery of the \englishmengentheoretisch\ paradoxes.

\yestop\indent
On a closer inspection,                                                \hskip.1em
however,                                                               \hskip.1em
it turned out that the paradoxes can be averted
by the exclusion of certain extreme \englishbegriffsbildungen,         \hskip.2em
which are not actually required for the constructions of set theory,   \hskip.2em
and even less for the methods of \englishanalysis.                     \hskip.3em
And yet,                                                               \hskip.1em
the discussion on the foundations of mathematics
has not settled down since then.                                       \hskip.4em
Moreover, \linebreak
those paradoxes were employed as motivational arguments
for restrictions on the mathematical \englishverfahrenvorgehen\
that go far beyond of what is 
immediately required for eliminating the \englishwidersprueche.
\actualparpagebreakpagenumber 6

\yestop\indent
A careful positioning in this discussion on the foundations of mathematics
requires a more detailed consideration of the
{\em logical structure of the mathematical theories}.

\yestop\indent
Let us notice that  
the mathematical \englishverfahrenvorgehenplural\
under discussion, in fact, \hskip.1em
consist of methods of inference and \englishbegriffsbildung. 
This means that some extension of the\linebreak usual logic 
\englishzurgeltungkommenthirdsingular\
in these \englishverfahrenvorgehenplural. \hskip.2em
At the same time,
we notice that the mathematical and the logical are tightly intertwined here:
\hskip.2em
On the one hand,
considering its objects, \hskip.1em
set theory comes very close to logic
through the \englishbeziehung\ between sets and predicates,
\ie\ through the \englishverhaeltnis\
between \englishumfangslogikohnelogik\ and content of \englishbegriffenotion.
\hskip.2em
On the other hand,
any systematic investigation of the modes of logical formation and inference
will necessarily raise mathematical \englishbetrachtungenproblem. \hskip.2em
As a matter of fact,
the traditional teaching of the categorical \englishschluesse\ 
is already a typical mathematical investigation; \hskip.2em
this fact is only superficially concealed
by the historical classification of this teaching into philosophy.
It is also connected with this mathematical side of the logical that the logical
\englishschluesse\
\mbox{---~in particular} as they occur in comprehensive logic application
   in mathematical theories~--- \hskip.1em
can be captured in a mathematical way 
as being composed of a small number of elementary processes.

\vfill\pagebreak

{\germanfont

\yestop\indent
Dieser \germansachverhalt\ wurde zur vollen Deutlich\-keit gebracht
\actualpagebreakpagenumber 7
durch die Entwicklung der \germansystem e der symbolischen Logik,
wie sie,
vorbereitet durch den \booleindex\boole'schen Logik\-kalkul,\linebreak
\edcomment{{\em sic}\/!}
um die Jahrhundertwende in\esi besondere von \peirceindex\peirce,
\fregeindex\frege, \schroederindex\schroeder, \peanoindex\peano,
\whiteheadindex\whitehead\ u.\ \russellindex\russell\ geschaffen wurden.
Bei der Konstruktion dieser \germansystem e ging man teil\es\ darauf au\es,
eine handliche Symbolik zu gewinnen,
die zugleich eine genauere Kontrolle der \germanschlussfolgerungoldspelling en
erm\oe glichte,
teil\es\ bezweckte man eine Einordnung der Mathematik in die Logik.\par

\begin{sloppypar}
\yestop\indent
E\es\ war der Gedanke \hilbertindex\hilbert\es,
die logische Symbolik dazu zu verwerten,
die mathematischen Bewei\esi methoden
zum Gegenstand einer mathematischen Untersuchung,
einer \mbox{\glqq\germanbeweistheorie\grqqcommanospace} \
zu machen.\footnote{%
  Der er\esi te Ansatz in dieser Richtung war der \hilbertindex\hilbert'sche
  \Heidelberg er Vortrag 
  1904 \cite{hilbert-grundlagen-logik},
  der freilich noch ganz im Fragmentarischen blieb.
  (In diesem wurde auch schon der Gedanke eine\es\ gemeinsamen
  \germanaufbau s von Mathematik und Logik zur \germangeltung\ gebracht.)
  \hskip.3em Eine er\esi te Weiterf\ue hrung dieser Gedanken findet sich,
  noch vor \hilbertindex\hilbert\es\ sp\ae teren Unter\-suchungen,
  in dem Werk von \koenigfathername:
  \koenigfathertitle\ (\Leipzig,\,1914) \cite{koenig-1914}.%
} \hskip.3em
Der wesentliche \germangesichtspunkt\ dabei ist,
die Methode der formalen \germanaxiomatik\ auch auf da\es\ logische
\germanschliessen\ selbst,
wie e\es\ in den Theorien der Arithmetik und Mengen\-lehre au\esi ge\ue bt wird,
anzuwenden und somit an die Stelle der Prozesse der logische\edcomment n
\germanbegriffsbildung\ und Folgerung
formal angesetzte Operationen treten zu lassen. \hskip.3em
Hierdurch gewinnen wir den Vorteil, \hskip.1em
da\sz\ wir bei strittigen Begriffen und \germanschlussweisenoldspelling\
nicht die \germaninhaltlich e Bedeutung in Betracht zu ziehen brauchen,
sondern nur den formalen Effekt, \hskip.1em
der durch ihre Anwendung in den deduktiven Prozessen bewirkt wird. \hskip.3em
Dieser Effekt\linebreak
l\ae\sz t sich vom \germanstandpunkt\ einer ganz elementaren
\germanbetrachtung\ verfolgen. \hskip.1em
Wir haben so\linebreak die M\oe glichkeit, \hskip.2em
Methoden, \hskip.1em
die vom \germaninhaltlich en \germanstandpunkt\ problematisch
erscheinen,\linebreak \hskip.1em
al\es\ bewei\esi technische \germanverfahren\
zu akzeptieren und zu rechtfertigen.\par
\end{sloppypar}

\yestop\indent
In diesem Sinne hat
\hilbertindex\hilbert\ die Au\fgi abe gestellt und in Angriff genommen,
da\es\ \germansystem\ der \germananalysis\ und Mengen\-lehre al\es\
\actualpagebreakpagenumber 8
ein widerspruch\esi freie\es\ Gedanken\-geb\ae ude zu erweisen.
Diese Au\fgi abe gliedert sich in zwei Teile.

\yestop\indent
E\es\ handelt sich \germanzunaechst\ darum,
die Bewei\esi methoden der \germananalysis\
und Mengen\-lehre einer formalen \germanaxiomatik\
zu unterwerfen oder, \hskip.1em
wie wir e\es\ kurz nennen wollen, \hskip.2em
zu {\em for\-ma\-li\-sieren}. \hskip.3em
Hierf\ue r konnte sich \hilbertindex\hilbert\ 
auf die bereit\es\ au\esi gebildeten zuvor genannten\linebreak
\germansystem e der \germanlogistik\ st\ue tzen,
in denen eine solche Formalisierung
bereit\es\ geleistet war.\linebreak
Der \germangesichtspunkt\ der streng formalen Deduktion wurde zuerst bei \frege\
scharf herau\esi gestellt und f\ue r Teile der Mathematik zur Durchf\ue hrung
gebracht. \hskip.3em
Die Methode zur hand\-lichen Au\esi gestaltung 
einer Formalisierung wurde durch \peanoindex\peano\ entwickelt. \hskip.3em
Eine Verbindung von beiden fand in den {\em\PM}\/ 
durch \whiteheadindex\whitehead\ und \russellindex\russell\ statt.\par

\yestop\indent
Die hier\edcomment{mit} \germanvorliegend e Formalisierung
ist freilich f\ue r die Zwecke der
\germanbeweistheorie\ insofern nicht vorteilhaft, \hskip.2em
al\es\ sie keine Gliederung in elementarere und h\oe here Bereiche der
\germanbegriffsbildung\ und
de\es\ \germanschliessen\es\ erm\oe glicht. \hskip.4em
Da\es\ r\ue hrt davon her, \hskip.1em
da\sz\ in den \PM\ sowie bei \fregeindex\frege\
die Gewinnung der \germanzahlentheorie\
au\es\ der allge\-meinen Mengen\-lehre
al\es\ eine\es\ der Haup\tzi iele genommen ist. \hskip.3em
So k\oe nnen hier die Methoden einer elementareren Behandlung
der \germanzahlentheorie\ nicht in Erscheinung treten.\par}

\vfill\pagebreak

\yestop\indent
This \englishsachverhalt\ was brought to full clarity 
\actualpagebreakpagenumber 7
by the development of the systems of symbolic logic,
which,
prepared by the \booleindex\myBoolean\ logic calculus,
were created around the turn of the century in particular by
\peirceindex\peirce, \fregeindex\frege, \schroederindex\schroeder,
\peanoindex\peano, \whiteheadindex\whitehead\ \& \russellindex\russell.
\hskip.3em
The construction of these systems partly aimed at obtaining a handy symbolism,
which also facilitated
a more precise control over the \englishschlussfolgerungen, \hskip.1em
and partly aimed at a classification of mathematics into logic.

\yestop\indent
It was \hilbertindex\hilbert's idea to create a ``proof theory\closequotecomma
in which the mathematical proof methods,
formalized in a logical symbolism,
become the objects of mathematical investigation.\footnote{%
  The first approach in this direction was \hilbertindex\hilbert's
  \Heidelberg\ talk of 1904 \cite{hilbert-grundlagen-logik},
  which, of course, had to remain outright fragmentary.
  (In this talk also the idea of a joint construction of 
   mathematics and logic is already \englishzurgeltungbringenppp.) \hskip.3em
  A first continuation of these ideas can be found,
  even before \hilbertindex\hilbert's later investigations,
  in the work of \koenigfathername: \hskip.1em
  \koenigfathertitle\
  \edcomment{\koenigfathertitleenglish}
  (\Leipzig,\,1914)
  \cite{koenig-1914}.%
} \hskip.2em
The significant \englishgesichtspunkt\ is to apply
the method of formal \englishaxiomatik\ also to the logical reasoning itself
as it is applied in the theories of arithmetic and set theory, \hskip.1em
and thus to let formally specified operations take the places of 
the processes of logical \englishbegriffsbildung\ and
\englishfolgerungconclusion. \hskip.2em
This has the advantage that we need not consider the
\englishinhaltlich\ meaning of contentious
\englishbegriffenotion\ and \englishschlussweisen, \hskip.1em
but only the formal effect
of their application in deductive processes. \hskip.2em
This effect can be traced
by means 
of most elementary \englishbetrachtungenueberlegungen. \hskip.2em
We thus have the possibility of accepting and justifying methods
as proof-technical \englishverfahrenprocedures, \hskip.1em
no matter whether they
are problematic from the \englishinhaltlich\ \englishstandpunkt.

\yestop\indent
In this sense, \hskip.1em
\hilbertindex\hilbert\ has set and tackled the task of proving the system of
\englishanalysis\ and set theory to be
\actualpagebreakpagenumber 8
\englisheinwiderspruchsfrei\
construct of ideas. \hskip.2em
This task subdivides into two parts.

\yestop\indent
\englishZunaechstfirst\ we have to represent
the proof methods of \englishanalysis\ and set theory
in a formal \englishaxiomatik\ or, \hskip.1em
as we want to call it briefly,     \hskip.1em
to {\em formalize}\/ them.         \hskip.3em
To this end, 
\hilbertindex\hilbert\ could resort to the mentioned full-fledged systems of \englishlogistik,
in which such a formalization had already been accomplished. \hskip.3em
\frege\fregeindex\ was the first who clearly displayed the \englishgesichtspunkt\
of strict formal deduction and \englisheszurdurchfuehrungbringen\ for parts of
mathematics. \hskip.2em
It was \peanoindex\peano\ who first developed handy presentations for  
such formalizations. \hskip.3em
In the {\em\PM}, \hskip.2em
\whiteheadindex\whitehead\ and \russellindex\russell\
combined these accomplishments of
\fregeindex\frege\ and \peanoindex\peano.

\yestop\indent
The formalization \englishvorlegenpppgiven\ by them is not advantageous
for the purposes of \englishbeweistheorie, however,
as it lacks subdividability into elementary and more advanced
\englishbereiche\ of \englishbegriffsbildung\ and \englishschliessen. \hskip.2em
This lack of subdividability has its historic origin in one of the main goals 
of the \PM\ \aswellas\ \frege's work: \hskip.1em
the extraction of \englishzahlentheorie\ from general set theory. \hskip.3em
Consequently, \hskip.1em
methods for a more elementary treatment of \englishzahlentheorie\ cannot
appear in their formalizations.

\vfill\pagebreak

{\germanfont
\yestop\indent
Von einer Formalisierung f\ue r die Zwecke der \germanbeweistheorie\
ist zu w\ue nschen,
da\sz\ sie eine analoge axiomatische Gliederung der logisch-mathematischen
\germanbildung en und Prozesse liefert,
wie sie in der \ue blichen \germanaxiomatik\ durch die Sonderung der
Axiomen\-gruppen bewirkt wird.
Unter diesem \germangesichtspunkt\ erscheint ein schichtweiser \germanaufbau\
de\es\ deduktiven \germanformalismus\ 
al\es\ angemessen.\actualparpagebreakpagenumber 9\par

\yestop\indent
So wurde man veranla\sz t,
von neuem die Formalisierung der mathematischen Di\esi ziplinen vor\-zu\-nehmen,
und man ist dabei zur pr\ae zisen Beschreibung von naturgem\ae\sz\ abgegrenzten
Teilbereichen der logisch-mathematischen Deduktion gelangt,
welche sich al\es\ solche zum er\esi ten Mal 
an Hand de\es\ Formalisierung\esi prozesse\es\ darstellten.\par

\yestop\indent
Auf Grund der vollzogenen Formalisierung gewinnt nun die Au\fgi abe
de\es\ Nachweise\es\ der \germanwiderspruchsfreiheit\
f\ue r die arithmetischen Di\esi ziplinen eine bestimmtere mathematische Form.
E\es\ handelt sich jetzt darum,
zu erkennen,
da\sz\ die festgelegten Prozesse der \mbox{Au\esi sagen}\-bildung 
und der \germanschlussfolgerungoldspelling\
nicht zur Herleitung solcher S\ae tze f\ue hren k\oe nnen, \hskip.1em
die einander im Sinne der gew\oe hnlichen \germaninhaltlich en Interpretation
widersprechen. \hskip.2em
Da\es\ Widersprechen von S\ae tzen stellt sich
mittel\es\ der Formalisierung der Negation
durch eine einfache Beziehung der entsprechenden Satzformeln dar.\par

\yestop\indent
\Ue berdie\es\ ergibt sich bei einer zweck\-m\ae\sz igen Formalisierung der
\germanaussagenlogik\ noch eine Ver\-einfachung in der Weise,
da\sz\ au\esi gehend von irgend einem \germanwiderspruch\
jede beliebige \germanaussage\ de\es\ formalisierten Bereich\es\
herleitbar wird.
Auf Grund davon gen\ue gt e\es,
um die Unm\oe glichkeit eine\es\ \germanwiderspruch\es\ in dem genannten Sinne
erkennen zu lassen,
da\sz\ man da\es\ Gegenteil eine\es\ bestimmten
einzelnen elementar g\ue ltigen Satze\es\ al\es\ nicht herleitbar erweist.
Die Au\esi f\ue hrung eine\es\ solchen Nachweise\es\
--~anschlie\sz end an den schrittweisen Au\fsi tieg der Teilbereiche~--
bildet den zweiten Teil de\es\ \hilbertindex\hilbert'schen Programm\es.
Die Durchf\ue hrung ist freilich \hilbertindex\hilbert\ nicht gelungen,
und e\es\ ist auch heute noch nicht abzusehen,
ob
--~oder vielmehr in welchem Sinne~--
sie gelingen kann.\actualparpagebreakpagenumber{10}\par

\begin{sloppypar}%
\yestop\indent
E\es\ bestehen n\ae mlich hier nicht nur gro\sz e technische, \hskip.1em
sondern auch grund\-s\ae tzliche Schwierig\-keiten.
Diese erheben sich vor allem mit Bezug auf die Frage, \hskip.1em
welche Mittel
f\ue r den gew\ue nschten Nachwei\es\ zugelassen werden sollen. \hskip.3em
In der Tat ging ja da\es\ Bed\ue rfni\es\ f\ue r einen solchen Nachwei\es\
von einer Kritik der \ue blichen Bewei\esi methoden au\es. \hskip.3em
Soll dieser Kritik Rechnung getragen werden, \hskip.1em
so darf der Nachwei\es\ der \germanwiderspruchsfreiheit\
nicht seiner\-seit\es\
auf einer Verwendung der kritisierten Methoden beruhen.%
\par\end{sloppypar}

\yestop\indent
Durch diese Erw\ae gung erhalten wir aber \germanzunaechst\ nur eine
Abgrenzung im negativen Sinne, \hskip.2em
und e\es\ bleibt noch die Au\fgi abe, \hskip.1em
genauer zu bestimmen, \hskip.1em
auf welche Arten der \Ue ber\-legung
die \germanbeweistheorie\ sich st\ue tzen soll. \hskip.3em
F\ue r die Wahl diese\es\ \germanstandpunkt\es\ wird un\es\
ein Anhalt durch da\es\ Erforderni\es\ gegeben, \hskip.15em
da\sz\ zumindest ja die Behauptung der \germanwiderspruchsfreiheit\ 
f\ue r die formalisierten Theorien sich pr\ae zise
fassen lassen mu\sz.\par}

\vfill\pagebreak

\yestop\indent
It is to be desired from a formalization for the purposes of
\englishbeweistheorie\ that it admits
an axiomatic subdivision of the logico-mathematical
formations and processes that is analogous
to the effect of the separation of axioms
into different groups, \hskip.1em
as it is often found in \englishueblich\ \englishaxiomatik. \hskip.25em
From this \englishgesichtspunkt, \hskip.1em
it appears to be adequate to construct
the deductive formalism in layers.\actualparpagebreakpagenumber 9\par

\yestop\indent
Therefore, \hskip.15em
the formalization of the mathematical disciplines
was put on the agenda again.                          \hskip.3em
And,                                                  \hskip.2em
in the course of the resulting formalization process, \hskip.2em
naturally separated sub-\linebreak domains\,of logico-mathematical\,deduction
were brought into separate\,being for\,the\,first\,time,\linebreak 
together with their precise descriptions.

\yestop\indent
On the basis of the given formalization, \hskip.2em
the task of proving the \englishwiderspruchsfreiheit\
of the arithmetical disciplines now
gains a more \englishbestimmt\ mathematical form. \hskip.3em
It now means to realize that the \englisheindeutigfestgelegt\ processes of
\englishschlussfolgerung\ and of \englishaussage\ \englishbildungformation\
cannot lead to the derivation of \englishsaetzetheoreme\
that contradict each other according to
the common \englishinhaltlich\ \englishinterpretation. \hskip.3em
By means of the formalization of negation, \hskip.15em
a contradiction of \englishsaetzetheoreme\ can\linebreak be represented
simply as the occurrence
of two corresponding formalized \englishsaetzetheoreme.\footnote{\edcomment{%
  This sentence may mean two different things; \hskip.25em
  in the order of preference: \ \
  1.~~``The realization of \englishsaetzetheoreme\
  that contradict each other can be reduced to tracing the
  derivation of two \englishsaetzetheoreme\ of the forms
  \math A and \overline A.'' \hskip.4em
  Soundness of this problem reduction requires
  that the derivation process satisfies {\it ex falso quod\-libet}, \hskip.2em
  and completeness of the reduction requires a certain form of completeness
  of the derivation \nolinebreak process. The latter, for instance,
  is not given for the derivation process in
  \cite[\litsectref{2.2}]{Wirth_Stolzenburg_Series_Specificity_AMAI_2015}.%
  \linebreak
  2.~~``Contradiction 
  may be defined in this context as
  two derivable \englishsaetzetheoreme\ of the forms
  \math A and \overline A.''%
}}

\yestop\indent
Moreover,                                                      \hskip.2em
there is a further simplification by the fact that,            \hskip.1em
for a formalization of \englishaussagenlogik\
appropriate for our purposes,                                  \hskip.1em
an arbitrary contradiction will admit the derivation of each and every 
\englishaussage\ of the formalized \englishbereich.            \hskip.2em
To realize the impossibility of a contradiction
in the mentioned sense,                                        \hskip.2em
it thus suffices to pick any \englishbestimmt,
elementary \englishgueltig\ \englishsatztheorem\ 
and show that the opposite of this \englishsatztheorem\
cannot be derived.
After the stepwise construction of the sub\-domains,           \hskip.2em
it is the execution of such a proof which forms the second part of
\hilbertindex\hilbertsprogram.                                 \hskip.3em
Admittedly,                                                    \hskip.2em
this execution was not achieved by \hilbertindex\hilbert,
\hskip.25em
and even today it is still not foreseeable whether
\mbox{---~or rather in which sense~---}
\mbox{it may be achieved.}\actualparpagebreakpagenumber{10}\par

\yestop\indent
In fact,   
there are not only major technical difficulties here,
but also fundamental ones.                                        \hskip.3em
These arise mainly with regard to the question,
which means are to be admitted for the \englishgewuenscht\ proof. \hskip.5em
It was,                                                           \hskip.1em  
after all,                                                        \hskip.1em
a critique of the \englishueblich\ proof methods 
that prompted the demand for such a proof.                        \hskip.5em
If this critique is to be taken into proper account,              \hskip.1em
then the proof of \englishwiderspruchsfreiheit\
must not in turn rely on the application of the criticized methods.

\yestop\indent
By this \englisherwaegung, however,
we get only a first boundary in the negative sense; \hskip.25em
and the task still remains to determine more precisely, \hskip.2em
on which
\englishartensorts\ of \englishueberlegung\
\englishbeweistheorie\ is to be based. \hskip.3em
A hint on the choice of this \englishstandpunkt\
is given by the requirement\linebreak that
at least
the \englishbehauptung\ of
\englishwiderspruchsfreiheit\
must be precisely expressible
for
the formalized theories.

\vfill\pagebreak

{\germanfont

\yestop\indent
Diesem Erforderni\es\ wird bereit\es\ gen\ue gt durch eine Art der elementaren
mathemati-\linebreak schen \germanbetrachtungsweise,                 \hskip.1em
welche \hilbertindex\hilbert\ al\es\ den
{\em\germanfinit en \germanstandpunkt}\/                             \hskip.1em
bezeichnet hat.\linebreak 
E\es\ \nolinebreak ist diejenige Art
anschaulicher mathematischer \Ue berlegung,                          \hskip.25em
wie sie in der elementa-\linebreak ren
Kombinatorik angewandt wird.                                         \hskip.3em
Auch die elementare \germanzahlentheorie\ und Buchstaben-Algebra
l\ae\sz t sich auf diese Art behandeln.                              \hskip.5em
Da\es\ Kennzeichende f\ue r die \germanfinit e Betrachtung besteht in
folgenden Momenten:
\begin{enumerate}
\yestop\item
Al\es\ \germangegenstaende\ werden nur {\em endliche}\/ Gebilde genommen,
an denen auch nur {\em diskrete}\/
Gestalt\esi merkmale unterschieden werden.
\yestop\item
Die Formen de\es\ allgemeinen und de\es\ \germanexistentialenurteils\
kommen nur auf eine ein\-ge\-schr\ae nk\-te Art zur Anwendung,
im Sinne \actualpagebreakpagenumber{11}
der Vermeidung der Vorstellung von unendlichen Gesamtheiten;
n\ae mlich da\es\ allgemeine Urteil wird nur in hypothetischem Sinne gebraucht,
al\es\ eine \germanaussage\ \ue ber jedweden vorliegenden Einzelfall,
und da\es\ \germanexistentialeurteil\ al\es\ ein
(zweck\-m\ae\sz ig zu vermerkender) Teil einer n\ae her bestimmten Feststellung,
in der entweder ein bestimmt strukturier\-te\es\ Gebilde vorgewiesen 
oder ein allgemeine\es\ \germanverfahren\ au\fgi ezeigt wird,
nach dem man zu einem (gewisse Bedingungen erf\ue llenden)
\germangegenstand\footnote{\edcomment{%
   This must be the plural in general: ``\germangegenstaende n''.%
}}
einen anderen \germangegenstand\
mit verlangten Eigenschaften gewinnen kann.
\yestop\item
Alle Annahmen,
die man einf\ue hrt,
beziehen sich auf endliche Konfigurationen.
\yestop\end{enumerate}

\noindent Mit der genaueren Beschreibung und Er\oe rterung der \germanfinit en
\germanbetrachtungsweise\ werden wir un\es\ noch de\es\ N\ae heren zu befassen
haben.\par

\yestop\yestop\indent
E\es\ w\ae re zweife\lli o\es\ sehr befriedigend,
wenn wir un\es\ in der \germanbeweistheoretisch en Untersuchung
v\oe llig an diesen Rahmen elementarer Betrachtung halten k\oe nnten.
Die M\oe glichkeit hierf\ue r scheint \germanzunaechst\ insofern gegeben zu sein,
al\es\ ja mit Bezug auf eine formalisierte Theorie
die Behauptung ihrer \germanwiderspruchsfreiheit\
sich nach dem vorhin Bemerkten
in \germanfinit er Form dahin au\esi sprechen l\ae\sz t,
da\sz\ ein jeder formalisierte Bewei\es\ eine \germanendformel\ hat,
die verschieden ist von der Negation einer bestimmten\edcomment{,}
geeignet gew\ae hlten Satzformel.\par}

\vfill\pagebreak

\yestop\indent
This requirement is already satisfied by a form of
elementary mathematical \englishbetrachtungsweise,
which \hilbertindex\hilbert\
called the {\em\englishfinit\ \englishstandpunkt}. \hskip.5em
It is the form of those \englishanschaulich\linebreak mathematical
\englishueberlegungsingularasplural\
which are applied in elementary combinatorics. \hskip.35em
Also elementary \englishzahlentheorie\ and algebra
can be treated in this form. \hskip.4em
The characteristics of \englishfinit\ considerations
are given by the following moments:
\begin{enumerate}
\yestop\item
The objects must be {\em finite}\/ entities, \hskip.25em
distinguished only  
by their {\em discrete}\/ \hskip.05em attributes.
\yestop\item\englishAllgemeineundexistentialeurteile\
may occur \onlyif\ \hskip.2em\actualpagebreakpagenumber{11}
they do not refer to the \englishvorstellung\
of any infinite \englishgesamtheit: \hskip.3em
The \englishallgemeinesurteil\ may be used
only in the hypothetical sense,                                     \hskip.2em
as \englisheineaussage\ on \englishjedwede\ 
individual object to be \englishvorlegenpppgiven.                   \hskip.3em
Moreover,                                                           \hskip.25em
the \englishexistentialesurteil\ may occur only as                  \hskip.1em
(an appropriately noted)                                            \hskip.2em
part of \englisheinnaeherbestimmter\ \englishfeststellungstatement, \hskip.2em
which additionally presents
either an entity\linebreak of \englishbestimmt\ structure,          \hskip.2em
or else a general \englishverfahrenprocedure\ for generating an object
of the\linebreak asserted properties from any object\edcomment{s}
to be \englishvorlegenpppgiven\ (satisfying certain conditions).\footnote{%
  \edcomment{%
    The objects which the general \englishverfahrenprocedure\
    (to be presented here) \hskip.1em 
    must accept as input arguments are those individual objects
    which may be given to the \englishallgemeineurteile\
    on which the \englishexistentialesurteil\ depends, \hskip.2em \ie\
    in whose scope the \englishexistentialesurteil\ occurs. \hskip.4em
    See \makeaciteoftwo{wirthcardinal}{wirth-simplified-epsilon} \hskip.1em
    for several formal frameworks where this statement
    makes sense even in the absence of quantifiers.\par
    The conditions that these individual objects may be assumed to satisfy
    (according to
     ``gewisse Bedingungen erf\ue llenden''/``satisfying certain conditions'')
    \hskip.1em are those which arise from the
    logical context of the \englishexistentialesurteil,
    \eg\ for \englisheinexistentialesurteil\ in
    the conclusion of an implication the individual objects may be assumed
    to satisfy the condition of the implication \etc\ \hskip.4em
    The required form of reasoning with this kind of logical context, \hskip.1em
    called {\em hierarchical contextual reasoning}, \hskip.1em
    was nicely formalized in \cite{sergediss} \hskip.1em for the first time.%
}}\yestop\item
\englishAnnahmen\ may be introduced only if they refer to finite
configurations.\yestop\end{enumerate}

\noindent We will have to take a closer look on 
a more detailed description and discussion
of the \englishfinit\ \englishbetrachtungsweise\ later.

\yestop\yestop\indent
There can be no doubt that it would be very fulfilling if we were able to
confine our \englishbeweistheoretisch\ investigations to this framework
of elementary consideration. \hskip.2em
\englishZunaechstatfirst, 
this seems to be possible insofar as \hskip.15em
---~according to what we have previously noticed~---\linebreak
the \englishbehauptung\ of the
\englishwiderspruchsfreiheit\ of a formalized theory can be
expressed in \englishfinit\ form by the \englishbehauptung\
that every formalized proof
has \englisheineendformel\ that is different from the negation of
\englisheinbestimmterdefinite, properly chosen formalized sentence.

\vfill\pagebreak

{\germanfont

\yestop\indent
Die \germanformulierbarkeit\ eine\es\ Problem\es\
im Rahmen gewisser Au\esi druck\esi mittel
bietet aber noch keine Gew\ae hr daf\ue r,
da\sz\ seine L\oe sung sich mit diesen Mitteln bewerkstelligen l\ae\sz t.
\germanpointertogoedelsincompletness\ \hskip.3em
\actualpagebreakpagenumber{12}
Man hat sich so gen\oe tigt gesehen,
f\ue r die \germanbeweistheorie\ den ur\-spr\ue ng\-lichen \germanfinit en
\germanstandpunkt\ zu einem \glqq konstruktiven\grqq\ \germanstandpunkt\
zu erweitern,
der sich etwa so kennzeichen l\ae\sz t,
da\sz\ von den drei soeben genannten Forderungen nur die er\esi ten beiden
au\fri echterhalten werden,
die dritte aber fallen gelassen wird.
So kommt man dazu,
f\ue r die \germanbeweistheorie\ eine ungef\ae hr solche methodische Haltung
einzunehmen,
wie sie der \brouwerindex\brouwer'sche \germanintuitionismus\
f\ue r die Mathematik
\ue berhaupt al\es\ einzig zul\ae ssig ansieht.\footnote{\label{note griss}%
  Eine pr\ae zise Vergleichung ist darum hier nicht zu verlangen,
  weil die \germanintuitionistisch e\ Haltung nicht durch Gebrauch\esi regeln,
  sondern durch eine philosophische Einstellung charakterisiert ist. \hskip.1em
  Das wird auch von \heytingindex\heyting\ ungeachtet der von ihm
  durchgef\ue hrten Formalisierung der \germanintuitionistisch en Logik
  und Mathematik hervorgehoben. \hskip.3em
  F\ue r eine Konfrontierung der Methoden der konstruktiven
  \germanbeweistheorie\ mit denen de\es\ \germanintuitionismus\ sind auch die
  neueren Unter\-suchungen von \grissshortname\ \hskip.1em
  \ue ber die negation\esi freie intuitionistische Mathematik von Belang,
  bei welchen da\es\ Operieren mit irrealen (unerf\ue llten) \germanannahme n
  grunds\ae tzlich vermieden wird
  (\Proc\ Kon.\ Ned.\ Adad.\ v.\ Wetensch.\
  49 (1946) \cite{griss-1946},                                    \hskip.5em
  53 (1950) \cite{griss-1950}                                     \hskip.2em
  und                                                             \hskip.1em
  54 (1951) \makeaciteoffour
  {griss-1951a}{griss-1951b}{griss-1951c}{griss-1951d}\edcomment).%
}\par

\yestop\indent
Doch selbst bei dieser Erweiterung de\es\ finiten \germanstandpunkt\es\
hat e\es\ nicht allenthalben sein Bewenden.
Man sieht sich vielmehr bei der Behandlung gewisser Fragen dazu gedr\ae ngt,
auch die zweite der obigen Forderungen fallen zu lassen.
Die\es\ ist zum Beispiel bei der Behandlung der \germanfragedervollstaendigkeit\
de\es\ \germansystem\es\ der Regeln f\ue r die gew\oe hnliche
\germanpraedikatenlogik\ der Fall,
die zuerst von \goedelname\ im positiven Sinne gel\oe st worden ist.
Hier erfordert bereit\es\ die \germanformulierung\ de\es\ Ergebnisse\es,
wenig\esi ten\es\ wenn man sie in pr\ae gnanter und einfacher Form haben will,
die Einf\ue hrung einer nichtkonstruktiven \germanbegriffsbildung.
Eine \germanfinit e Fassung de\es\ Ergebnisse\es\
und auch ein Bewei\es\ im \germanfinit en Rahmen l\ae\sz t sich erzwingen,
ist aber mit technischen Komplikationen belastet.
\par

\yestop\indent
Diese\es\ Beispiel de\es\ \goedel'schen \germanvollstaendigkeit\esi satze\es\
ist zugleich daf\ue r charak\-teri\esi tisch,
da\sz\ die \germanbeweistheoretisch e \germanfragestellung\ in ihrer
nat\ue rlichen Au\esi gestaltung nicht bei dem Problem verbleibt,
welche\es\ au\es\ der Kritik der \ue blichen \germanverfahren\ der
klassischen Mathematik erwachsen ist. \hskip.3em
Auch \hilbertindex\hilbert\ hat ja von vornherein die Au\fgi abe der
\germanbeweistheorie\ sehr weit gefa\sz t.\footnote{%
  \majorheadroom
  \Vgl\ seine \Ae u\sz erungen in dem Vortrag:
  \titlehilbertnineteenhundredandseventeen\
  (Math.\,Ann.\,78, \PP{405}{415} (1918)) \cite{hilbert-1917}.%
}\actualparpagebreakpagenumber{13}\par

\yestop\indent
Unter den weitergehenden Fragestellungen sind nun auch etliche solche,
bei denen die Verbindlichkeit der Anforderung einer methodischen
Beschr\ae nkung al\es\ fraglich erscheint.
So sind in neuerer Zeit verschiedene\xout n erfolgreiche Unter\-suchungen,
die in weiterem Sinne zum Felde der \germanbeweistheorie\ geh\oe ren,
im Rahmen der \ue blichen mathematischen Methodik,
also ohne Beschr\ae nkung der \germanbegriffsbildungen\ und
Bewei\esi methoden, 
durchgef\ue hrt worden.
Andererseit\es\ haben verschiedene Autoren f\ue r jene speziellen
\germanbeweistheoretisch en Unter\-suchungen,
die e\es\ mit den Fragen der \germanwiderspruchsfreiheit\ zu tun haben,
einen solchen \germanstandpunkt\ gew\ae hlt,
bei welchem von den drei vorhin formulierten Forderungen
nur die er\esi te zugrunde gelegt wird.}
\vfill\pagebreak\par
The \englishformulierbarkeit\ of a problem
within the limits of certain means of expression,           \hskip.1em
however,                                                    \hskip.2em
does not guarantee that its solution
can be accomplished with these means \aswell.                 \hskip.3em
\englishpointertogoedelsincompletness\footnote{%
   \label{note goedel}%
  \edcomment{Can this refer to anything else but \cite{goedel}\,?}%
} \hskip.3em
\actualpagebreakpagenumber{12}
Therefore, \hskip.1em
it \nolinebreak became necessary to extend the original
\englishfinit\ \englishstandpunkt\
for \englishbeweistheorie\ to a \mbox{``constructive''}
one,                                                        \hskip.1em
which can be characterized roughly by maintaining only the
first two requirements of the three just mentioned ones,    \hskip.1em
but dropping the third.                                     \hskip.3em
Thus,                                                       \hskip.1em
we arrive at at position toward the methods in \englishbeweistheorie\   
which is roughly the one that 
\brouwerindex\brouwer's \englishintuitionismus\ considers to be
the only \englishzulaessig\ one in mathematics in general.\footnote{%
  \label{note griss english}%
  \majorheadroom
  A precise comparison is not to be demanded here,                \hskip.1em
  because the \englishintuitionistisch\ position
  is actually not characterized by rules of application,          \hskip.1em
  but by a philosophical attitude.                                \hskip.1em
  This is emphasized also by \heytingindex\heyting,               \hskip.1em
  notwithstanding his formalization of
  \englishintuitionistisch\ logic and mathematics.                \hskip.3em
  Moreover,                                                       \hskip.1em
  a most relevant comparison of the methods of constructive
  \englishbeweistheorie\ with those of \englishintuitionismus\
  is found in the more recent investigations of                   \hskip.1em
  \grissshortname\  \hskip.1em
  on negation-free \englishintuitionistisch\ mathematics,         \hskip.1em
  where the handling of unreal (\ie\,\englishunerfuellbar)  
  \englishannahmen\ is avoided on principle
  (\Proc\ Kon.\ Ned.\ Adad.\ v.\ Wetensch.\
  49 (1946) \cite{griss-1946},                                    \hskip.5em
  53 (1950) \cite{griss-1950}                                     \hskip.2em  
  and                                                             \hskip.1em
  54 (1951) \makeaciteoffour
  {griss-1951a}{griss-1951b}{griss-1951c}{griss-1951d}).%
}

\yestop\indent
Even this extension of the
\englishfinit\ \englishstandpunkt,                                \hskip.1em
however,                                                          \hskip.1em
does not suffice in all cases.                                    \hskip.3em
\Infact,                                                          \hskip.1em   
for treating certain \englishfragenproblems,                      \hskip.1em   
we feel urged
to drop the second of the above requirements \aswell.             \hskip.3em 
This is for example the case for the treatment of the
\englishfragedervollstaendigkeit\ 
of the rule system for ordinary \englishpraedikatenlogik,         \hskip.2em
first solved by \goedelname\
---~in the positive sense.                                        \hskip.4em
In this example,                                                  \hskip.1em
a concise and simple \englishformulierung\ of the result
already requires the introduction of a non-constructive
\englishbegriffsbildung.                                          \hskip.3em
\englishEinfinit\ version of the result can be enforced,          \hskip.1em
and also a proof by \englishfinit\ means,                         \hskip.15em 
but these versions are 
burdened
with technical complication.\footnote{\majorheadroom
\edcomment{This seems to be a reference to \herbrandsfundamentaltheorem\
  \cite{herbrand-PhD},
  which, together with the \loewenheimskolemtheorem\
  \cite{loewenheim-1915},
  yields \englisheinfinit\ version of \goedel's completeness theorem
  \cite{goedel-completeness},
  \cfnlb\ \eg\
  \makeaciteoftwo{herbrand-handbook}{SR--2009--01}, 
  \makeaciteoftwo{wirth-heijenoort}{wirth-heijenoort-SEKI}.}%
}\par

\yestop\indent
Thus,                                                          \hskip.1em
\goedel's \englishvollstaendigkeit\ theorem provides also
a characteristic example
that \englishbeweistheoretisch\ \englishfragestellungendefinition\ 
in their natural embodiment already exceed
the problem that resulted from the critique
of the \englishueblich\ \englishverfahrenprocedures\
in classical mathematics.                                      \hskip.3em
Accordingly,                                                   \hskip.1em
\hilbertindex\hilbert\ already gave a comprehensive outlook on the tasks of
\englishbeweistheorie\ from the very beginning.\footnote{\majorheadroom%
  \Cf\ his statements in his talk:
            \titlehilbertnineteenhundredandseventeen\
  \edcomment\titlehilbertnineteenhundredandseventeenenglish\
  (Math.\,Ann.\,78, \PP{405}{415} (1918))
  \cite{hilbert-1917}.%
}\actualparpagebreakpagenumber{13}

\yestop\indent
Moreover,                                                 \hskip.1em
for quite a few among
the more advanced \englishfragestellungendefinition,      \hskip.1em
the commitment to a restriction of the methods
may be questioned anyway.                                 \hskip.3em
\Infact,                                                  \hskip.1em
in more recent times,                                     \hskip.1em
in the field of \englishbeweistheorie\
in the broader sense,                                     \hskip.1em
various successful investigations
have been realized with the full methodology of
\englishueblich\ mathematics,                             \hskip.1em
restricting neither
\englishbegriffsbildungen\ nor proof methods.             \hskip.3em 
Furthermore,                                              \hskip.1em
for \englishbeweistheoretisch\ investigations
explicitly concerned with
the \englishfragenderwiderspruchsfreiheit,                \hskip.1em
however,                                                  \hskip.1em
various authors have chosen \englisheinstandpunkt\
that relies only on the first of
the three previously stated restrictions.

\vfill\pagebreak

{\germanfont

\yestop\indent
Eine endg\ue ltige Entscheidung der Methodenfrage
kann dieser \germansachlage\ gegen\ue ber jeden\-fall\es\ nur erwartet werden,
wenn man einen \germanueberblick\ dar\ue ber hat,
wa\es\ die verschiedenen Methoden zu leisten verm\oe gen.
E\es\ kann schwerlich behauptet werden,
da\sz\ gegenw\ae rtig eine Entscheidung jener Frage vorliegt,
die nicht blo\sz\ durch eine vorgefa\sz te philosophische Ansicht bestimmt ist.
Und bez\ue glich der verschiedenen sich bek\ae mpfenden und heute
\ue blicherma\sz en gegen\ue bergestellten philosophischen Lehr\-meinungen
besteht der Verdacht,
da\sz\ sie ungekl\ae rte \germanvoraussetzung en\ in sich schlie\sz en,
die ihrerseit\es\ vielleicht eher fragw\ue rdig sind
al\es\ die angefochtenen mathematischen Theorien.\par

\yestop\indent
Angesicht\es\ dieser \germansachlage\ erscheint e\es\ al\es\ da\es\
angemessene \germanverfahren,
da\sz\ wir einer\-seit\es\ die methodische Richtlinie,
die durch den \germangesichtspunkt\ der \germanfinit en Betrachtung gegeben wird,
im Auge behalten,
andererseit\es\ un\es\ \edcomment{aber dadurch}
\label{page finitism as guideline without method}%
nicht in der Methode fest\-legen.   \hskip.3em
\germanpointertohasenjaegersdiss\ \hskip.4em
Bei solchen Anwendungen fungiert die \germanbeweistheorie\ 
nicht in der Rolle der Bewei\esi kritik,
sondern im Rahmen der \ue blichen Methoden
de\es\ mathematischen \germanschliessen\es,
und e\es\ w\ue rde darum hier die Forderung der \germanfinit en
\germanbetrachtungsweise\ gar nicht am Platze sein.
\\[+1ex]\LINEnomath{{\bf---}}\par

\yestop\indent
Zur Vorbereitung unserer \germanbeweistheoretisch en Betrachtungen
ist e\es\ nun auf jeden Fall w\ue nschen\esi wert,
da\sz\ die spezifische Art der \germanfinit en \Ue berlegung
deutlich gemacht werde.
Zur Illustrierung eignet sich besonder\es\ da\es\ Gebiet 
der elementaren \germanzahlentheorie,
in welcher der \germanstandpunkt\ der direkten \germaninhaltlich en,
ohne axiomatische Annahmen sich vollziehenden \Ue berlegungen
am reinsten au\esi gebildet ist.\actualparpagebreakpagenumber{15}
}

\vfill\pagebreak
In any case
---~confronted with this overall \englishsachlage~--- \hskip.1em
a final answer to the \englishfragequestion\
which methods to choose can only be expected after
attaining a comprehensive insight into the full capability
of the different methods.                                         \hskip.3em
\mbox{It can hardly} be claimed
that an answer to that \englishfragequestion\
has been given until today,                                       \hskip.3em
unless determined by
preconceived philosophical beliefs.     \hskip.3em
And the various competing philosophical doctrines                 \hskip.1em
---\,\,nowadays mostly confronted with each other\footnote{\edcomment{%
  The English phrase ``confronted\,with\,each\,other'' is just as
  ambiguous as the German original
  ``gegen\-\ue ber\-ge\-stellt\closequotecomma
  meaning either to be brought into opposition
  or to be arranged face-to-face for a
  comparison.%
}}\,---                                                      \hskip.1em
are suspected to comprise uncleared \englishvoraussetzungen. \hskip.4em
Moreover,                                                    \hskip.1em
these \englishvoraussetzungen\ may well be more questionable
in turn than the mathematical theories under challenge.

\yestop\indent
\label{page english: finitism as guideline without method}%
\selfquotedsentenceattheendofsectiononeofourtypescript\ \hskip.5em
\englishpointertohasenjaegersdiss{\footnote{\label
{note hasenjaeger}\majorheadroom
  \edcomment{To what can this refer? \hskip.3em
            Is there any candidate prior to \cite{hasenjaeger-1950c},
            entitled ``{\germanfontfootnote\hasenjaegernineteenfiftyctitle}''
            (``\hasenjaegernineteenfiftyctitleenglish'')\,?}%
}} \hskip.5em
\mbox{In such applications,}                                         \hskip.1em 
\englishbeweistheorie\ does not act
in the \role\ of a proof critique,                                   \hskip.2em 
but within the framework of the
\englishueblich\ methods of mathematical \englishschliessen,         \hskip.1em 
and therefore the demand for \englisheinfinit\ \englishbetrachtungsweise\
would not be appropriate here \atall.\\\mbox{}\\\LINEnomath{{\bf---}}

\yestop\indent
For the preparation
of our \englishbeweistheoretisch\ \englishbetrachtungenueberlegungen, \hskip.1em
it is now desirable in any case
that the specific mode of \englishfinit\ \englishueberlegungsingularasplural\
is demonstrated.                                                      \hskip.3em
Particularly suitable for such an illustration is the field of
elementary \englishzahlentheorie,                                     \hskip.1em
where the \englishstandpunkt\ of direct \englishinhaltlich\
\englishueberlegung\ without axiomatic \englishannahmen\ 
is developed most purely.\actualparpagebreakpagenumber{15}


\vfill\pagebreak

With the previous paragraph ends our long excerpt from our typescript
and its \litsectref 1, presented here together with its English translation
on the right-hand side pages (\ie\ those with uneven numbers).
\vfill\cleardoublepage
\yestop
\subsection
{Discussion and Excerpts of \hskip.1em\litsectref 2 of \hskip.2em Our Typescript}
\begin{sloppypar}\yestop\yestop\noindent
\litsectref 2 of our typescript comes with the subsection headline
``{\litsectref 2.~\germanfont Die elementare \germanzahlentheorie}''
(``\litsectref 2.~Elementary \englishZahlenTheorie'') \hskip.1em and 
contains a version similar to
the first nine pages (\PP{20}{28}) \hskip.1em
of \hskip.1em\litsectref 2 of \hskip.05em
\makeaciteoftwo{grundlagen-first-edition-volume-one}
               {grundlagen-second-edition-volume-one}, \hskip.2em
which comes with a similar headline:
\par\halftop\indent
{``{\germanfont\litsectref 2.~~Die elementare \germanzahlentheorie. \ ---\ \
    Da\es\ \germanfinit e \germanschliessen\ und seine Grenzen.}'' \ \ }
\par\halftop\noindent In English:
\par\halftop\indent
{``\litsectref2.~~Elementary \englishZahlenTheorie. \ --- \ \ 
\englishFinit\ \mbox{Inference} and its Limits.'' \ \ }
\par\yestop\yestop\noindent
The first five paragraphs of \hskip.1em\litsectref 2 of \hskip.05em
\makeaciteoftwo{grundlagen-first-edition-volume-one}
               {grundlagen-second-edition-volume-one}, \hskip.2em
however, \hskip.1em 
are missing in our typescript. \hskip.4em
For the first of these paragraphs,
this is just a consequence of the
omission of several paragraphs 
at the beginning of
\nolinebreak\hskip.15em\nolinebreak \litsectref 1
(\cfnlb\ our \sectref{section Parallel introductory part}), \hskip.2em
because it is again on the subject of the
``\hspace*{-.2em}{\em existential form}\closequotefullstopextraspace
The further four of these  paragraphs missing in our typescript
contain a digression
into the subject of
geometry.\par
\end{sloppypar}\yestop\yestop\yestop
\subsubsection
{A sentence following the second edition instead of the first one}
\label{section A sentence following the second edition instead of the first one}
\yestop\noindent
A significant difference to
\cite{grundlagen-first-edition-volume-one}
is that
\ ---\ instead of the paragraph\notop\halftop\begin{quote}
``{\germanfont Diese Figuren bilden eine Art von \germanziffern;
 wir wollen hier da\es\ Wort 
 \glq{\em\germanziffer}\/\grq\ \hskip.1em
 schlechtweg zur Bezeichnung 
 {\em dieser}\/ Figuren gebrauchen.}''\par
 ``These figures constitute a kind of \englishziffer; 
  and we will simply 
  use the word
  `\nolinebreak\hspace*{-.1em}\nolinebreak{\em\englishziffer}\/' \hskip.1em 
  to \englishbezeichnen\ just {\em these}\/ figures.''\notop\halftop\end{quote}
of \cite[\p\,21]{grundlagen-first-edition-volume-one}''~~--- \   
we find in our typescript the following sentence of 
\cite[\p\,21]{grundlagen-second-edition-volume-one}:\notop\halftop\begin{quote}
``{\germanfont Wir wollen diese Figuren,
  mit einer leichten Abweichung vom gewohnten \germansprachgebrauch, \hskip.1em
  al\es\ \glq{\em\germanziffern}\/\grq\ \hskip.1em bezeichen.}''\par
``Deviating slightly 
from the \englishgewoehnlichersprachgebrauch, 
we will call these figures
`{\nolinebreak\hspace*{-.12em}\nolinebreak
  \em\englishziffern}\/'.''\yestop\end{quote}
It is unlikely that this text was just copied from the second edition,
however, because our typescript contains the addition
``{\germanfont(in Ermangelung eine\es\ besseren kurzen \germanausdruck\es)}''
(``(lacking a better short \englishausdruck)''), \hskip.2em
which is deleted by \bernaysindex\bernays' hand.\vfill\pagebreak
\subsubsection{A significant improvement compared to both editions}
\label{subsection A significant improvement compared to both editions}

\halftop\indent
Another significant difference in our typescript occurs
in the enumerated list of five elements
found on page~21\f\ in both editions of \Vol\,I
\makeaciteoftwo
{grundlagen-first-edition-volume-one}
{grundlagen-second-edition-volume-one}. \hskip.3em
This list describes the syntactic form of the 
``{\germanfont Zeichen zur Mi\tti eilung}''
(``symbols for communication''): \hskip.3em
In \lititemref 1 we \nolinebreak read\notop\halftop\begin{quote}
``{\germanfont kleine deutsche Buchstaben zur Bezeichnung f\ue r
               unbestimmt gelassene \germanziffern;}''\par
``small German letters for designating \englishziffern\
that are left undetermined;''\halftop\notop\end{quote}
in our typescript, \hskip.1em
instead of\halftop\notop\begin{quote}
``{\germanfont Kleine deutsche Buchstaben zur Bezeichnung f\ue r
   irgendeine nicht festgelegte \germanziffer;}''\par
``Small German letters for designating \englishirgendein, not determined \englishziffer;''\notop\halftop\end{quote}
found in both editions of ``\grundlagendermathematik\closequotefullstopnospace

\halftop
Together with other hints on the \englishfinit\ \englishstandpunkt\
found in our typescript,
this correction was crucial for a change in translation 
in the third English edition
\cite{grundlagen-german-english-edition-volume-one-one}, \hskip.2em
as compared to the first two English editions where we translated the
latter German \nolinebreak term
\mbox{---~after} communication with several members
of the advisory board of the \hilbertbernaysproject,
who agreed that the meaning is ambiguous~---
with the similarly ambiguous term
``\englishirgendein\ indeterminate
\englishziffer\closequotefullstopnospace\footnote{%
  As long as we were not able to give the German term a unique reading,
  we had to be careful and use a translation that does not cut off
  any possibly meaningful reading. \hskip.3em
  It \nolinebreak was our typescript here
  that provided the information for such  a
  unique reading. \hskip.3em
  Now we can translate that reading and
  do not have to bother our readers
  with an ambiguous English \nolinebreak term. \hskip.3em
  \par
  As long as there are several possible readings,
  a translator must not give 
  only his favorite interpretation of the original in the translation. \hskip.3em
  If the translator knows exactly what is meant, however,
  then that meaning should to be captured in the translation
  as clearly and unambiguously as possible; \hskip.3em
  otherwise the translator would propagate the weed of
  imperfect expression
  that always comes with the crop,
  in particular regarding the fruits of science.%
} \hskip.5em
Our new translation in the third edition is\notop\halftop\begin{quote}``{%
  small German letters for \englishbezeichneningform\ \englishirgendwelche,
  not determined \englishziffern;}''\notop\halftop\end{quote}
This new version clearly disambiguates these small German letters
from {\em formal}\/\footnote{%
  \majorheadroom
  Whereas German letters
  ---~all over \makeaciteoffour
  {grundlagen-first-edition-volume-one}
  {grundlagen-first-edition-volume-two}
  {grundlagen-second-edition-volume-one}
  {grundlagen-second-edition-volume-two}~---
  do not denote formal variables
  (as Latin letters do), \hskip.2em 
  in the metalogical, \englishfinit\ framework,
  depending on the natural-language context,
  German letters may be used both
  for {\em free atoms}\/ and
  for {\em free variables}\/ in the terminology of \cite{wirth-simplified-epsilon}
  (both of which are substantially different from the (bound and free)
   \englishindividuenvariablen\
   and the free \englishformelvariablen\ of \hilbertbernays).\par
  The technical term ``\hspace*{-.2em}{\em atom}'' is from set theories
  (with atoms or urelements) and comes implicitly with
  a universal quantification (such as in the equations on \p\,30 of
  \makeaciteofthree
   {grundlagen-first-edition-volume-one}
   {grundlagen-second-edition-volume-one}
   {grundlagen-german-english-edition-volume-one-one}), \hskip.2em
  whereas the term ``\hspace*{-.1em}{\em free variable}\/''
  is from free-variable semantic
  tableaus \makeaciteoftwo{Fitting90}{fitting} \hskip.1em
  and comes implicitly with an
  (\math\varepsilon-restricted) existential quantification
  (such as in the equations on \p\,28 of
   \makeaciteofthree
    {grundlagen-first-edition-volume-one}
    {grundlagen-second-edition-volume-one}
    {grundlagen-german-english-edition-volume-one-one}).%
}
variables \aswellas\ from
arbitrary objects in the sense of \citet{fine-arbitrary-objects}.

\halftop\indent
For a more detailed discussion see \litnoteref{21.6} on \p\,21.b in
the third English edition
\cite{grundlagen-german-english-edition-volume-one-one}.
\vfill\pagebreak
\halftop\subsection
{Discussion and Excerpts of \hskip.2em\litsectref 3 of \hskip.1em Our Typescript}

\begin{sloppypar}
\yestop\halftop\noindent
\litsectref 3 is the final section of our typescript. \hskip.5em
It comes with the subsection headline
\ ``{\germanfont\litsectref 3.~\germanueberschreitung\
   de\es\ \germanfinit en \germanstandpunkt e\es\
   im mathematischen \germanschliessen.}'' \hskip.3em
(``\litsectref 3.~\englishUeberschreitung\
   of the \englishfinit\ \englishstandpunkt\
   in mathematical \englishschliessen.''). \hskip.9em

\yestop
\subsubsection
{From our typescript: the entire text of \hskip.1em\litsectref 3,
 \hskip.1em \PP{32}{34}}%
\halftop\halftop\noindent
In its first two paragraphs of \hskip.1em\litsectref 3, \hskip.2em
our typescript announces to discuss the most interesting and critical
\englishverfahrenvorgehen\ to \englishfinitismus\ 
chosen in the penultimate paragraph of \hskip.1em\litsectref 1
of our typescript. \hskip.3em
In \nolinebreak the penultimate paragraphs
on \spagerefs{page finitism as guideline without method}
{page english: finitism as guideline without method}
here, \hskip.1em
we have quoted the sentence on this choice and translated it 
as follows:\notop\halftop\begin{quote}
``\selfquotedsentenceattheendofsectiononeofourtypescript''\halftop\end{quote}
The bad news is that \litsectref 3 is obviously truncated in our typescript:
It
breaks off abruptly
after the first discussion of an abstract example on
\englishvollstaendigeinduktion,
not covering all what was announced in its first two paragraphs.

As this truncation results
in a very short,
but most interesting section (\PP{32}{34}),         \hskip.2em 
we \nolinebreak present this section here in total, \hskip.2em
in the same way we presented the new introduction from \litsectref 1
in \sectref{section A completely new introduction}: \hskip.4em
\textongermanoriginalanditsenglishtranslationrectoverso

\end{sloppypar}


\vfill\pagebreak
  
{\germanfont

\yestop\indent
Unsere au\esi gef\ue hrte \germanbetrachtung\ der \germananfangsgruende\ der
\germanzahlentheorie\ diente dazu,
un\es\ da\es\ direkte \germaninhaltlich e,
in Gedanken\-experimenten an \germananschaulichvorgestellt en Objekten
sich vollziehende und von axiomatischen Annahmen freie \germanschliessen\
in seiner Anwendung und Handhabung vorzuf\ue hren.
Wir haben un\es\ dabei an die methodische Einstellung gehalten,
die wir anfang\es\ nach \hilbertindex\hilbert\ al\es\ den
``\germanfinit en \germanstandpunkt'' bezeichnet haben.
Wir wollen nun de\es\ n\ae heren betrachten,
wie man dazu veranla\sz t wird,
den \germanfinit en \germanstandpunkt\ zu \germanueberschreiten.
Dabei wollen wir ankn\ue pfen an die fr\ue her gegebene Kennzeichnung de\es\
\germanfinit en \germanstandpunkt e\es,
die ja mittel\es\ der drei charakteri\esi tischen Momente erfolgte:
\hskip.4em 1.~Beschr\ae nkung der
\germangegenstaende\ auf endliche di\esi krete Gebilde; \hskip.4em
2.~Beschr\ae nkung der Anwendung der logischen Formen de\es\ allgemeinen
und de\es\ exi\esi tentialen \germanurteil\es\ im Sinne der Vermeidung der
\germanvorstellung\ von fertigen unendlichen Gesamt\-heiten; \hskip.4em
3.~Beschr\ae nkung der Annahmen
auf solche \ue ber endliche Kon\-figura\-tio\-nen.\par

\yestop\indent
Diese Momente sind geordnet im Sinne einer zunehmenden Anforderung.
Wir werden nun bei der \germanbetrachtung\ der \germanueberschreitung\
de\es\ \germanfinit en \germanstandpunkt e\es\
naturgem\ae\sz\ in entgegen\-gesetzter Reihenfolge,
im Sinne einer schritt\-weisen Abstreifung der Anforderungen, 
verfahren. \actualparpagebreakpagenumber{33}\par

\yestop\indent
Ein Versto\sz\ gegen die dritte Forderung,
wonach alle Annahmen sich auf endliche Kon\-figura\-tio\-nen beziehen sollen,
liegt bereit\es\ \ue berall da vor,
wo man die Annahme der \germangueltigkeit\ eine\es\ allgemeinen Satze\es\
\ue ber \germanziffern\ einf\ue hrt.\par

\yestop\indent
Eine Veranlassung dazu ist in\esi besondere gegeben bei Anwendung der
\germanvollstaendigeninduktion\ zum Beweise von S\ae tzen,
welche eine Beziehung \nolinebreak\hskip.05em
\nlbmath{\apptotuple{\mathfrak A}{\pair{\mathfrak m}{\mathfrak n}}} \hskip.05em
f\ue r \germanbeliebig e \germanziffern\
\nolinebreak\hskip.05em\nlbmath{\mathfrak m},
\nolinebreak\hskip.05em\nlbmath{\mathfrak n} \hskip.1em behaupten. \hskip.3em
Soll die Induktion,
etwa nach \nolinebreak\hskip.05em\nlbmath{\mathfrak n}, \hskip.1em
im \germanfinit en Sinne erfolgen,
so mu\sz\ bei dem \germanschlussoldspelling\ von \nolinebreak\hskip.05em
\nlbmath{\apptotuple{\mathfrak A}{\pair{\mathfrak m}{\mathfrak n}}} \hskip.05em
auf \nolinebreak\hskip.1em
\nlbmath{\apptotuple{\mathfrak A}
{\pair{\mathfrak m}{\mathfrak n\hskip.04em\tight+1}}} \hskip.05em
die \germanziffer\ \nolinebreak\nlbmath{\mathfrak m} \hskip.05em
festgehalten werden.
In dieser Weise sind wir auch im vorigen Paragraphen bei den Beweisen der
Rechengesetze f\ue r Summe und Produkt verfahren.\par

\yestop\indent
H\ae ufig wird aber die \germanvollstaendigeinduktion\ so angewandt,
da\sz\ man \germanzunaechst\ zeigt,
da\sz\ f\ue r jede \germanziffer\ \nlbmath{\mathfrak m} \hskip.05em
die Beziehung \nolinebreak\hskip.05em
\nlbmath{\apptotuple{\mathfrak A}{\pair{\mathfrak m}1}} \hskip.05em
besteht,
und sodann beweist,
da\sz,
fall\es\ f\ue r die \germanziffer\ \nlbmath{\mathfrak n} \hskip.05em
bei jeder \germanbeliebig en \germanziffer\ \nlbmath{\mathfrak m} \hskip.5em
\apptotuple{\mathfrak A}{\pair{\mathfrak m}{\mathfrak n}} \hskip.05em
besteht,
dann auch bei jeder \germanziffer\ \nlbmath{\mathfrak m} \hskip.5em
\apptotuple{\mathfrak A}{\pair{\mathfrak m}{\mathfrak n\hskip.04em\tight+1}}
\hskip.05em besteht.
Man schlie\sz t darau\es\ nach der \germanvollstaendigeninduktion,
da\sz\ f\ue r jede \germanziffer\ \nlbmath{\mathfrak n} \hskip.05em gilt,
da\sz\ f\ue r jede \germanziffer\ \nlbmath{\mathfrak m} \hskip.5em
\apptotuple{\mathfrak A}
{\pair{\mathfrak m}{\mathfrak n\hskip.04em
\tight+1
}}
\hskip.05em besteht.\par

\yestop\indent
Hier hat man in der zweiten zu beweisenden \germanbehauptung\ einen
Allsatz al\es\ \germanpraemisse; \hskip.2em
e\es\ wird ja angenommen,
da\sz\
(f\ue r den fixierten Wert \nlbmath{\mathfrak n}) \hskip.05em
bei jeder \germanziffer\ \nlbmath{\mathfrak m} \hskip.05em
die Beziehung \nolinebreak\hskip.05em
\nlbmath{\apptotuple{\mathfrak A}{\pair{\mathfrak m}{\mathfrak n}}} \hskip.05em
bestehe.
Diese\es\ Vorau\esi gesetzte k\oe nnen wir un\es\ nicht
in der \germanvorstellung\ eigentlich ver\-gegenw\ae rtigen.\par}
\vfill\pagebreak

\yestop\indent
Our \englishbetrachtungbehandlungplusgenitive\ 
the basics of \englishzahlentheorie\ 
was meant to demonstrate 
the application and the \englishhandhabung\ of
direct \englishinhaltlich\ inference
in thought experiments on 
\englishanschaulichvorgestellt\ \englishobjekte, \hskip.1em
free of axiomatic assumptions.\footnote{\edcomment{%
  This sentence is very similar to the first sentence of the third
  paragraph on \p\,32 in \makeaciteofthree
    {grundlagen-first-edition-volume-one}
    {grundlagen-second-edition-volume-one}
    {grundlagen-german-english-edition-volume-one-one}.%
}} \hskip.3em
In this \englishbetrachtungbehandlung, \hskip.1em
we have observed the \englishmethodischeeinstellung\
we initially called the ``\englishfinit\ \englishstandpunkt''
according to \hilbertindex\hilbert. \hskip.3em
\mbox{We now} want to \englishbetrachteneroertern\ in more detail
what may be the cause for \englisheineueberschreitung\ of 
the \englishfinit\ \englishstandpunkt. \hskip.3em
We want to follow here the previously given characterization of the
\englishfinit\ \englishstandpunkt\
by means of three characteristic moments: \hskip.5em
1.~limitation of any object to be a finite discrete entity; \hskip.4em
2.~limitation of any application of the logical forms of 
\englishdasallgemeineunddasexistentialeurteil\ 
to avoid the \englishvorstellung\
of any infinite \englishfertigegesamtheit; \hskip.4em
\mbox{3.~limitation} of any \englishannahme\
to refer only to finite configurations.

\yestop\indent
These moments are ordered by increasing demand. \hskip.2em
Now, \hskip.1em
in our \englishbetrachtungeroerterung\ of the
\englishueberschreitung\ of the \englishfinit\ \englishstandpunkt\
---~in the sense of a stepwise discarding of demands~---
\mbox{we will} naturally proceed in the \englishumgekehrtorder\ order.
\actualparpagebreakpagenumber{33}

\yestop\indent
A violation of the third requirement, \hskip.1em
according to which
all \englishannahmen\ should refer to finite configurations, \hskip.1em
is already \englishvorliegendgiven\ wherever one introduces
the \englishannahme\ of the \englishgueltigkeit\
of \englisheinallgemeinuniversal\ \englishsatztheorem\ about \englishziffern.

\yestop\indent
In particular, \hskip.1em
a cause for such a violation may be given in an application of 
\englishvollstaendigeinduktion\
for the proof of \englisheinsatzproposition\ that \englishbehauptenthirdsingular\
\englisheinebeziehung\ \nolinebreak\hskip.05em
\nlbmath{\apptotuple{\mathfrak A}{\pair{\mathfrak m}{\mathfrak n}}} \hskip.05em
for \englishbeliebig\ \englishziffern\
\nolinebreak\hskip.05em\nlbmath{\mathfrak m},
\nolinebreak\hskip.05em\nlbmath{\mathfrak n}. \hskip.5em
If the induction, \hskip.1em
say on \nolinebreak\hskip.05em\nlbmath{\mathfrak n}, \hskip.15em
is to be carried out in the \englishfinit\ sense, \hskip.1em
then the \englishziffer\ \nolinebreak\nlbmath{\mathfrak m} \hskip.05em
must be held constant in the conclusion step from
\nolinebreak\hskip.05em
\nlbmath{\apptotuple{\mathfrak A}{\pair{\mathfrak m}{\mathfrak n}}} \hskip.05em
to \nolinebreak\hskip.1em
\nlbmath{\apptotuple{\mathfrak A}
{\pair{\mathfrak m}{\mathfrak n\hskip.04em\tight+1}}}. \hskip.5em
In this way we proceeded also in the previous \englishparagraphsection\
in the proofs of the laws of calculation for sum and product.

\yestop\indent
Frequently, \hskip.1em
however, \hskip.1em
\englishvollstaendigeinduktion\ is applied in such a way that one
\englishzunaechstfirst\ shows that the \englishbeziehung\
\nolinebreak\hskip.05em
\nlbmath{\apptotuple{\mathfrak A}{\pair{\mathfrak m}1}} \hskip.05em
holds for every \englishziffer\ \nlbmath{\mathfrak m}, \hskip.25em
and then proves for the \englishziffer\ \nlbmath{\mathfrak n} \hskip.05em
that \hskip.05em
\apptotuple{\mathfrak A}{\pair{\mathfrak m}{\mathfrak n\hskip.04em\tight+1}}
\hskip.05em holds for every \englishziffer\ \nlbmath{\mathfrak m}, \hskip.25em
under the assumption that \hskip.05em
\apptotuple{\mathfrak A}{\pair{\mathfrak m}{\mathfrak n}} \hskip.05em
holds for \englishjedebeliebige\
\englishziffer\ \nlbmath{\mathfrak m}. \hskip.5em
From this, \hskip.1em
we may conclude by \englishvollstaendigeinduktion\ 
for every \englishziffer\ \nlbmath{\mathfrak n}\,%
:
\hskip.5em\apptotuple{\mathfrak A}
{\pair{\mathfrak m}{\mathfrak n\hskip.04em
\tight+1
}}
\hskip.05em holds 
for every \englishziffer\ \nlbmath{\mathfrak m}.

\yestop\indent
In\,the second \englishbehauptung\,to\,be shown\,here, \hskip.1em
we\,have\,a\,universal \englishsatzproposition\
as \englishpraemisse; \hskip.2em
indeed,\linebreak
we assume that
(for the fixed\,value \nlbmath{\mathfrak n}) \hskip.1em
the \englishbeziehung\ \nolinebreak\hskip.05em
\nlbmath{\apptotuple{\mathfrak A}{\pair{\mathfrak m}{\mathfrak n}}} \hskip.05em
holds for every \englishziffer\ \nlbmath{\mathfrak m}. \hskip.5em
\mbox{In our} \englishvorstellung, \hskip.1em
\mbox{we are}
not really able to bring clearly to our mind
what we have to \englishvoraussetzen\
here.\footnote{\majorheadroom\edcomment{%
  Indeed, \hskip.1em
  we cannot visualize
  ---~element by element~---
  the set \hskip.2em
  \setwith
      {\apptotuple{\mathfrakfootnote A}
                  {\pair{\mathfrakfootnote z}{\mathfrakfootnote n}}}
      {\mathfrakfootnote z\tightin\N}   \hskip.3em
  of all potential assumptions, \hskip.2em
  but only effectively given, finite subsets of it.%
}}

\vfill\pagebreak
{\germanfont

\yestop\indent
Freilich l\ae\sz t sich in vielen F\ae llen die genannte Form der Anwendung
der \germanvollstaendigeninduktion\ vom \germanfinit en \germanstandpunkt\
motivieren. \hskip.2em
Da\es\ ist \zB\ dann der Fall,
wenn beim Bewei\es\ de\es\ Bestehen\es\ von \nolinebreak\hskip.05em\nlbmath{
\apptotuple{\mathfrak A}
           {\pair{\mathfrak m}{\mathfrak n\hskip.04em\tight+1}}} \hskip.05em
f\ue r ein bestimmte\es\ \nlbmath{\mathfrak m} \hskip.05em
die \germanvoraussetzung\ de\es\ Bestehen\es\ von \nolinebreak\hskip.05em
\nlbmath{\apptotuple{\mathfrak A}{\pair{\mathfrak z}{\mathfrak n}}} \hskip.05em
f\ue r \germanbeliebig e \nlbmath{\mathfrak z} \hskip.05em nur \hskip.05em
\actualpagebreakpagenumber{34}
in solcher Weise zur Anwendung kommt,
da\sz\ eine durch \math{\mathfrak m} und \nlbmath{\mathfrak n}
bestimmte endliche \germananzahl\ von Beziehungen
\\[+.9ex]\noindent(1)\LINEmaths{
  \apptotuple{\mathfrak A}{\pair{\mathfrak k_1}{\mathfrak n}}\comma\ldots\comma
  \apptotuple{\mathfrak A}{\pair{\mathfrak k_{\mathfrak r}}{\mathfrak n}}
  }{}\mbox{~~~}
\\[+.8ex]\noindent benutzt wird,
worin \math{\mathfrak k_1,\ldots,\mathfrak k_{\mathfrak r}} gewisse au\es\
\math{\mathfrak m} und \nlbmath{\mathfrak n}
zu ermittelnde \germanziffern\ sind.\par 

\yestop\indent
In diesem Falle kommt ja der Bewei\es\ de\es\
hypothetischen
Satze\es\ mit der
All\-pr\ae misse darauf hinau\es,
da\sz\ man zeigt,
da\sz\ sich die Feststellung der Beziehung \nolinebreak\hskip.05em\nlbmath{
\apptotuple{\mathfrak A}
           {\pair{\mathfrak m}{\mathfrak n\hskip.04em\tight+1}}} \hskip.05em
(f\ue r die fixierten \germanziffern\ \nlbmath{\mathfrak m},
\nlbmath{\mathfrak n}) \hskip.05em
zur\ue ck\-f\ue hren l\ae\sz t auf die Feststellung der Beziehungen~(1),
und damit ist ein Regre\sz\ gegeben,
der die Feststellung von \nolinebreak\hskip.05em\nlbmath{
\apptotuple{\mathfrak A}
           {\pair{\mathfrak m}{\mathfrak n\hskip.04em\tight+1}}} \hskip.05em
in einer begrenzten \germanzahl\ von Schritten auf die Feststellung 
von Beziehungen
\\[+.2ex]\noindent\LINEmaths{
  \apptotuple{\mathfrak A}{\pair{\mathfrak z_1}1}\comma\ldots\comma
  \apptotuple{\mathfrak A}{\pair{\mathfrak z_{\,\mathfrak\es}}1}}{}
\\[+.1ex]\noindent zur\ue ckf\ue hrt,
und f\ue r diese wird durch den anf\ae nglichen Bewei\es\ de\es\ Bestehen\es\
von \nolinebreak\hskip.05em\nlbmath{
\apptotuple{\mathfrak A}{\pair{\mathfrak m}1}} \hskip.05em
f\ue r \germanbeliebig e \nlbmath{\mathfrak m} \hskip.05em
da\es\ allgemeine \germanverfahren\ gegeben.\par

\yestop\indent
Diese Art der Rechtfertigung der erweiterten Form der
\germanvollstaendigeninduktion\ mit All\-s\ae tzen al\es\ \germanpraemisse n
ist jedoch nicht generell anwendbar.
In\esi besondere erw\ae chst eine Schwierigkeit au\es\ dem Umstand,
da\sz\ derartige erweiterte Induktionen in komplizierter Weise
ineinander\-geschachtelt sein k\oe nnen.
Wir wollen einen typischen Fall dieser Art n\ae her betrachten.
E\es\ handelt sich dabei um einen konstruktiv behandelbaren Teil
der \cantorindex\cantor schen Theorie der tran\esi finiten Ordinalzahlen.%
\actualparpagebreakpagenumber{35}}
\vfill\pagebreak

\yestop\indent
In many cases it is \infact\ possible
to justify the considered \englishanwendung\
of \englishvollstaendigeinduktion\,from\,the\,\englishfinit\,\englishstandpunkt
.\footnote{%
  \edcomment{For a formal treatment of this motivation see
  \PP{348}{351} of \makeaciteoftwo{grundlagen-second-edition-volume-one}
  {grundlagen-german-english-edition-volume-one-three}.}%
} \hskip.2em
For\,instance,
such\,a\hskip.22em justification\,is\,always\,possible\hskip.18em
if\linebreak
---~in the proof \hskip.1em that \nolinebreak\hskip.1em\nlbmath{
\apptotuple{\mathfrak A}
           {\pair{\mathfrak m}{\mathfrak n\hskip.04em\tight+1}}} \hskip.05em
holds for \englisheinbestimmterparticular\ \nlbmath{\mathfrak m}\,\,~---
\hskip.5em
the application the 
hypothesis that \hskip.08em
\nolinebreak\hskip.05em
\nlbmath{\apptotuple{\mathfrak A}{\pair{\mathfrak z}{\mathfrak n}}} \hskip.05em
holds for \englishbeliebig\ \nlbmath{\mathfrak z} \hskip.05em
can be restricted
\hskip.05em \actualpagebreakpagenumber{34}
to a finite \englishanzahldrei\ of \hskip.15em
\edcomment{instances of this universal \englishsatzproposition\ in the form of\,}
\hskip.15em\englishbeziehungen\
\\[+.9ex]\noindent(1)\LINEmaths{
  \apptotuple{\mathfrak A}{\pair{\mathfrak k_1}{\mathfrak n}}\comma\ldots\comma
  \apptotuple{\mathfrak A}{\pair{\mathfrak k_{\mathfrak r}}{\mathfrak n}}
  },\mbox{~~~}
\\[+.8ex]\noindent
where \hskip.1em\math{\mathfrak r} \hskip.05em\aswellas\ \hskip.1em
\math{\mathfrak k_1,\ldots,\mathfrak k_{\mathfrak r}} \hskip.15em
are certain \englishziffern\ that must be 
\englishbestimmenppp\ from 
\math{\mathfrak m} and \nlbmath{\mathfrak n}.

\yestop\indent
In this case the proof of the 
proposition\footnote{\majorheadroom\edcomment{%
  The occurrence of \hskip.15em``hypothetischen'' in the German original
  is to make clear that ``Satzes'' (``sentence'') does not refer to a theorem,
  but just to a proposition, namely the induction step here. \hskip.3em
  The German wording clashes with the standard notion of an
  {\em induction hypothesis}. \hskip.3em
  Indeed, \hskip.1em
  we cannot translate this occurrence as ``hypothetical sentence''
  here because we already translated ``\germanvoraussetzung''
  (``\englishvoraussetzung'') in the previous paragraph as
  ``hypothesis\closequotefullstopnospace%
}}
with the universal \englishpraemisse\ 
amounts to showing that the \englishfeststellungdurchnachweis\
of the \englishbeziehung\ \nolinebreak\hskip.05em\nlbmath{
\apptotuple{\mathfrak A}
           {\pair{\mathfrak m}{\mathfrak n\hskip.04em\tight+1}}} \hskip.05em
(for the fixed \englishziffern\ \nlbmath{\mathfrak m},
\nlbmath{\mathfrak n}) \hskip.05em
can be \englishzurueckfuehreniSlogRedppp\
to the \englishfeststellungdurchnachweis\ of the
\englishbeziehungen~(1)\@. \hskip.5em
And therefore, \hskip.1em
we are then given a regression procedure
that \englishzurueckfuehreniSlogRedthirdsingular\
the \englishfeststellungdurchnachweis\ of
\nolinebreak\hskip.15em\nlbmath{\apptotuple{\mathfrak A}
           {\pair{\mathfrak m}{\mathfrak n\hskip.04em\tight+1}}} \hskip.05em
to the \englishfeststellungdurchnachweis\ of \englishbeziehungen
\\[+.2ex]\noindent\LINEmaths{
  \apptotuple{\mathfrak A}{\pair{\mathfrak z_1}1}\comma\ldots\comma
  \apptotuple{\mathfrak A}{\pair{\mathfrak z_{\,\mathfrak\es}}1}}{}
\\[+.1ex]in a limited number of steps. \hskip.4em
Moreover, \hskip.1em
for these \englishbeziehungen, \hskip.1em
there is a general \englishverfahrenprocedure\ given by the initial
proof that \nolinebreak\hskip.1em\nlbmath{
\apptotuple{\mathfrak A}{\pair{\mathfrak m}1}} \hskip.05em
holds for \englishbeliebig\ \nlbmath{\mathfrak m}.

\yestop\indent
This way of justifying the extended form of \englishvollstaendigeinduktion\
with universal \englishsaetzepropositions\ as \englishpraemissen, \hskip.1em
however, \hskip.1em
is not applicable in general. \hskip.3em
In particular, \hskip.1em
one of the difficulties that may arise here is the circumstance
that such extended inductions can be \englisheingeschachtelt\
into each other in a complicated way. \hskip.3em
Let us \englishbetrachteneroertern\ a typical case of this kind
in more detail. \hskip.3em
We will deal here with a part of \hskip.05em \cantorindex\cantor's theory
of transfinite ordinal numbers that can be treated constructively.
\actualparpagebreakpagenumber{35}


\vfill\pagebreak

With the previous paragraph ends our excerpt presenting
\litsectref 3 of our typescript in total
---~here together with its English translation
on the right-hand side pages (\ie\ those with uneven numbers).

The original text suddenly breaks off here. \hskip.3em
It is not clear whether the remainder
is just
incorrectly filed or whether \bernaysindex\bernays\ further draft
(in \englishgabelsbergerkurzschrift?)\linebreak
---~which must have been written before the introduction to this
subsection fragment~---
was never put into typewriting.

\vfill\cleardoublepage
\section{Trying to Find Hints on the Time of Writing}\label
{section Trying to Find Hints on the Time of Writing}%
\subsection{\litsectref 1 of \hskip.1em Our Typescript
(Hints: Before\,1929, 1931, or 1951--1977)}\label
{1 of Our Typescript}\subsubsection
{Possible explanations for the missing of the two paragraphs
\\on the ``existential form'' in the introductory part}\label
{section Parallel introductory part II}%
\notop\halftop\noindent To explain why
the two paragraphs from the beginning of \hskip.15em\litsectref 1
of \hskip.15em\makeaciteoftwo{grundlagen-first-edition-volume-one}
{grundlagen-second-edition-volume-one} \hskip.15em
are missing in our typescript
(as described in \sectref{section Parallel introductory part}), \hskip.2em
we see the following {\em three options}\/:\begin{enumerate}\noitem\item
The missing of the two paragraphs may indicate
that our typescript was written before\,1929
for the following reason: \hskip.3em
Soon later the subject of these paragraphs,
the ``existential form\closequotecomma
has become an integral part of the standard presentation,
which was never dropped in publications; \hskip.3em
\cfnlb\ the initial sections
of \hskip.1em
\cite{bernays-phil-math-hilbert-beweistheorie},
\makeaciteoftwo{grundlagen-first-edition-volume-one}
{grundlagen-second-edition-volume-one}.\item
It \nolinebreak could have been the case
that \bernaysindex\bernays\ or the \hilbertindex\hilbert\ school in logic
temporarily dropped the idea of the
``
exi\esi tentiale Form''
after the shock of \goedelsincompletenesstheorems,
say in\,1931. 
\item
A third option is that the whole text was written by \hasenjaegername\
during his stay as \bernaysindex\bernays' assistant in \Zuerich, or later as a consequence
of this stay. \hskip.3em
What speaks for this option is the overall improved pedagogical quality of
our typescript as compared to the corresponding sections in
\makeaciteoftwo{grundlagen-first-edition-volume-one}
{grundlagen-second-edition-volume-one}. \hskip.7em
Then these paragraphs may be not actually missing, \hskip.1em
but just deferred to a later section where they are more appropriate
from the pedagogical \englishgesichtspunkt.\end{enumerate}
\subsubsection{The completely unknown introduction}\label
{section The completely new introduction}%
\notop\halftop\noindent In the completely unknown introduction to  
foundations of mathematics and \hilbertindex\hilbert's \englishbeweistheorie\
(\cfnlb\ \sectref{section A completely new introduction}), \hskip.2em
there are some references to relevant published work that may help to
date our typescript.\begin{enumerate}\item
There is the following reference (\cfnlb\ \noteref{note goedel},
\sectref{section A completely new introduction}):\notop\halftop\begin{quote}
``\germanpointertogoedelsincompletness''\par
``\englishpointertogoedelsincompletness''\notop\end{quote}
We do not think that this could have been written before autumn\,1930, \hskip.1em
when \goedelsincompletenesstheorems\ \shortcite{goedel} became known. \hskip.3em
If this is so, \hskip.1em
then our typescript cannot have been
written before autumn\,1930, \hskip.1em
which excludes the first option of the previous subsection.\item
Moreover,
there is the reference to
\makeaciteoffour{griss-1951a}{griss-1951b}{griss-1951c}{griss-1951d} \hskip.1em
in \noterefs{note griss}{note griss english} \hskip.1em
of \hskip.15em
\sectref{section A completely new introduction}. \hskip.5em
Unless the notes on the extra pages
(\cfnlb\,\sectref{section additional footnotes}) \hskip.1em were added
much later to our typescript
than our typescript itself was originally written
\pagebreak
(which is very unlikely because the position markers of the
footnotes are already part of our typescript), \hskip.3em
this means that our typescript
cannot have been written before 1951.\item
Furthermore,
there is also the following reference (\cfnlb\ \noteref{note hasenjaeger},
\sectref{section A completely new introduction}):\notop\halftop\begin{quote}
``{\germanfont\germanpointertohasenjaegersdiss}''\par
``\englishpointertohasenjaegersdiss{}''%
\notop\end{quote}

\noindent
We have no idea to which publication before \cite{hasenjaeger-1950c} \hskip.1em
the contact with topology could refer.\notop\halftop\end{enumerate}
If our assessment in the last two items is correct, \hskip.1em
then also the second option of the previous subsection
can be excluded and our typescript was probably written by
\hasenjaegerindex\hasenjaeger\ during or after his stay with \bernaysindex\bernays\ in \Zuerich.
\subsection{\litsectref 2 of \hskip.1em Our Typescript
(Hints: \hasenjaegerplain, 1951--\bernaysdeathyear)}\label
{2 of Our Typescript}%
The sentence of our typescript 
found in the second instead of the first edition of ``\grundlagendermathematik''
(which we discussed in \nlbsectref
{section A sentence following the second edition instead of the first one})
\hskip.1em
clearly speaks in favor of a version written during the preparation of the
second edition; \hskip.3em
and the only work in this context we \nolinebreak know about is
the joint work of \hasenjaegerindex\hasenjaeger\ and \bernaysindex\bernays.

The sentence of our typescript improving on both editions of ``\grundlagendermathematik''
(which we discussed in
\nlbsectref{subsection A significant improvement compared to both editions})
\hskip.1em
makes it likely that it was written with the intention\linebreak of
a thorough revision; \hskip.3em
and the only such attempt we know about is
the joint work of
\hasenjaegerindex\hasenjaeger\ and \bernaysindex\bernays. \hskip.4em
Moreover, \hskip.2em
as \bernaysindex\bernays\ used to keep track 
of any correction \mbox{very carefully}
in his author's copies\commanospace\footnote{%
  Such as the one
  of the first edition \makeaciteoftwo
  {grundlagen-first-edition-volume-one}
  {grundlagen-first-edition-volume-two}
  owned by \engelername\ and the one of the second edition \makeaciteoftwo
  {grundlagen-second-edition-volume-one}
  {grundlagen-second-edition-volume-two}  
  owned by \bernaysrenename,
  \cfnlb\ \makeaciteofthree
  {grundlagen-german-english-edition-volume-one-one}
  {grundlagen-german-english-edition-volume-one-two}
  {grundlagen-german-english-edition-volume-one-three}.%
} \hskip.2em
such an improvement over both editions strongly indicates that
our typescript may be the result of a collaboration of \hasenjaeger\
and \bernaysindex\bernays\ after the publication of the second edition of the
first volume (1968), \hskip.23em
\ie\ in the years from 1968 to \bernaysdeathyear\
(when \bernaysindex\bernays\ died).
\halftop\subsection{\litsectref 3 of \hskip.1em Our Typescript (No Hints)}
Although \litsectref 3 \hskip.05em
of our typescript is very interesting and we would be
keen on reading the remainder of it if it were found, \hskip.3em
it breaks off too soon
for giving us a clear hint on the time of its writing.\vfill\pagebreak
\section{\hasenjaegerplain\ and \bernaysplain}\label
{section hasenjaeger and bernays}%
As we have seen, \hskip.1em
for finding out the time of writing of our typescript, \hskip.1em
it may be helpful 
to find out more about \hasenjaeger's scholarship as  \bernaysindex\bernays' assistant
in \Zuerich\ for the preparation of the second edition of 
``\grundlagendermathematik'' in the early 1950s
(\cfnlb\ \sectref
{section The Mentioning of the Given-Up Work on the Second Edition}). \hskip.3em
So the questions is: 
{\em What do we know about the relation of 
  \hasenjaegerindex\hasenjaegerplain\ and \bernaysindex\bernaysplain\
  around the year\,\,1950\,?}

Besides several biographical remarks from \hasenjaeger\
in \cite{menzler-gentzen-german},
\hskip.3em
there seem to be essentially only two non-trivial biographical texts on \hasenjaeger:
\hskip.25em One is a laudation by
\dillerindex\diller\ \shortcite{diller-on-hasenjaeger}; \hskip.25em
the other one is on his time as a cryptologist responsible for the security
of the German Enigma in the second world war and appears in similar forms in 
\makeaciteofthree
    {schmeh-hasenjaeger-0}
    {schmeh-hasenjaeger-1}
    {schmeh-hasenjaeger-2}. \hskip.3em
There we learn that he was born \Jun\,1, \hasenjaegerbirthyear,
in \nolinebreak\Hildesheim\ as a son of a lawyer\commanospace\footnotemark\
\hskip.2em that he was seriously wounded on January\,2, 1942,
as \nolinebreak a German soldier in Russia, \hskip.2em
and that the famous logician, philosopher, and theologian
\scholzname\ \scholzlifetime\
saved him from being ordered to the Russian front again,
by recruiting him for the cryptology department of the High Command
of the German Armed Forces (OKW/Chi) in \Berlin.

As we were neither able to find a reasonable short CV of
\hasenjaeger\ in publications nor in the \WWW\hspace*{-.1em}, \hskip.2em
we have put his own one and a further one compiled by us into our Appendix.
\hskip.3em For a list of his publications,
see our hint at the beginning of our \refname.
\notop\halftop\subsection{Early Relation}
\begin{sloppy}
\hasenjaeger\ and \bernaysindex\bernays\ were exchanging letters on the first edition of
``\grundlagendermathematik''
since\,1943, \hskip.1em
and from\,1949 on also on \hasenjaeger's own work.\footnotemark\ \hskip.3em
They definitely met each other in autumn\,1949 at the
``Kolloquium zur Logistik und der mathematischen
Grundlagen\-forschung\closequotecomma
\Sep\,27\,--\,\Oct\,1,\,1949,
Mathematisches Forschung\esi institut Oberwolfach (MFO),
Oberwolfach (Germany). \hskip.3em
Indeed,
the following photo\footnotemark\
shows \hasenjaeger\ in the back row to the right behind \bernaysindex\bernays\
(with \arnoldschmidtname\ in dark suit
to the right of \bernaysindex\bernays, and \schuettename\ rightmost):\end{sloppy}

\halftop
\noindent\LINEnomath{\includegraphics[bb=0 0 400 277,width=%
.6667%
\linewidth]{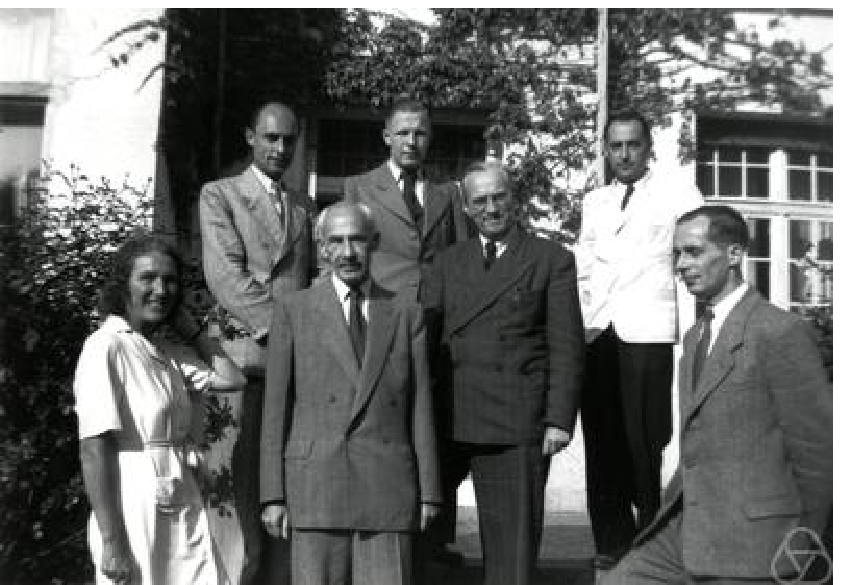}}\pagebreak

\halftop\noindent\bernaysindex\bernays\ \hskip.05em seems to have reviewed
the \PhDthesis\ of \hskip.1em\hasenjaeger\
\shortcite{hasenjaeger-1950c} \hskip.1em
before \May\,24, 1950, \hskip.3em
which is the date of the final report
by \hasenjaeger's supervisor \scholzname.%
\footnotemark

\addtocounter{footnote}{-3}\footnotetext{\sloppy%
  According to \cite[\p 186]{menzler-gentzen-german} \hskip.1em
  and also to \hasenjaeger's daughter \beatebeckername, \hskip.2em
  \hasenjaegername's 
  father was \hasenjaegerfathername\ \hasenjaegerfatherlifetime, \hskip.2em
  who was the mayor of \Muelheim\ from 1936 to 1946,
  \hskip.1em and 
  a very close friend of 
  \hskip.1em \scholzname\
  (according to\linebreak \url
  {https://www.muelheim-ruhr.de/cms/edwin_hasenjaeger_-_portrait_eines_oberbuergermeisters_1936-19461.html}).%
}\addtocounter{footnote}{1}\footnotetext{%
  \majorheadroom
  In particular, \hasenjaeger\ and \bernays\ were exchanging letters
  on the following publications: \makeaciteoffive{hasenjaeger-1950a}
  {hasenjaeger-1950b}{hasenjaeger-1950c}{hasenjaeger-1952a}{hasenjaeger-1953a}.%
}\addtocounter{footnote}{1}\footnotetext{%
  \majorheadroom
  Source: Univ.\,Archiv \Freiburg: Depositalbestand E6.
  Scanned, imprinted, and usage permitted by
  Mathematisches Forschung\esi institut Oberwolfach (MFO),
  Oberwolfach (Germany):
  \url{https://opc.mfo.de/detail?photo_id=11199}
}\addtocounter{footnote}{1}\footnotetext{%
  \majorheadroom
  For more information on \hasenjaeger' promotion,
  see the item labeled \cite{hasenjaeger-1950c}
  in our \refname.%
}
\subsection{Joint Work on ``\grundlagendermathematik''}\label
{section Joint Work}%
\begin{sloppypar}
All in all, \hskip.1em
three scholarships were granted to \hasenjaeger\ by the
\ETHshort\ \hskip.05em Board,\footnote{\majorheadroom%
  In German: ``\ETHshort-Rat\closequotecomma
  at that time called ``Schweizerischer Schulrat''. \hskip.4em
  For the grants see \makeaciteoftwo{hasenjaeger-1951}{hasenjaeger-1952g}.%
} \hskip.05em
with the goal to assist \bernaysindex\bernays\ in the preparation of the
second edition of ``\grundlagendermathematik\closequotefullstopextraspace
\hasenjaeger's first stay
as \bernaysindex\bernays' assistant in \Zuerich\ took place
in winter term 1950/51.\footnote{\majorheadroom
  \Cfnlb\ \cite{hasenjaeger-1950e}.%
} \hskip.4em
As planned from the beginning, \hskip.1em 
a second stay followed in summer term 1951.

Only after this second stay in Z\ue rich, \hskip.1em
however, \hskip.1em
\hasenjaeger\ 
produced the {\em first proper sketch}\/\footnote{%
  \majorheadroom
  This sketch
  attached to \cite{hasenjaeger-1951} \hskip.1em
  may only have been written in response to \scholz' letter
  \shortcite{scholz-1951a} to \bernaysindex\bernays, \hskip.2em
  where \scholz\ expresses his worries raised by an
  oral report from \hasenjaeger\
  (who must have met \scholz\ in \MuensternoWestfalen\
   between the years 1950 and
   1951) \hskip.1em
  that the work on ``\grundlagendermathematik''
  had not even started \nolinebreak yet:\par
  ``Inzwischen ist Herr \hasenjaeger\ hier erschienen,
  und ich habe ihn als aller\-er\esi te\es\ gefragt nach dem Stande der Arbeit an
  dem Grundlagenbuch. Er hat mir gesagt, dass diese Arbeit noch gar nicht
  hat anlaufen k\oe nnen, weil Sie die Mengenlehre
  erst fertig machen m\ue ssen. Die\es\ hat mich nun wirklich so erschreckt,
  dass ich es Ihnen mit der er\esi ten M\oe glich\-keit sagen muss.''\par
  In English: ``In the meantime, \hskip.1em
  \MrUS\ \hasenjaeger\ appeared here, \hskip.15em
  and my very first question was about the state of affairs of the work
  on the foundations book. \hskip.3em
  He \nolinebreak told \nolinebreak
  me that this work could not even have started yet,
  because you first
  have to complete the set theory. \hskip.3em
  Now this has alarmed me to such an extent
  that I \nolinebreak have to tell this to you with the first opportunity.''\par
  The noun phrase ``the set theory'' refers to
  \cite{bernays-axiomatic-set-theory-VII}
  (received \Mar\,30, 1953), \hskip.2em
  where \hasenjaeger\ \hskip.1em is mentioned only in the form of a reference to
  \cite{hasenjaeger-1953a}.\par
  In the eyes of \hskip.1em\scholz, \hskip.2em
  it was even worse for \hasenjaeger's career
  that the work of \hasenjaeger\ and \bernaysindex\bernays\
  on the second edition did not even start
  before summer~1951, \hskip.2em
  \cfnlb\ \makeaciteoftwo{scholz-1951c}{scholz-1951d}.%
}
of the overall layout for the
second edition of ``\grundlagendermathematik\closequotecommasmallextraspace
in the form of a rough table of contents. \hskip.3em
%
This sketch was and still is attached
to \hasenjaeger's  letter \shortcite{hasenjaeger-1951} \hskip.1em
to \bernaysindex\bernays. \hskip.3em
It is an elaborated and augmented protocol of the first detailed\linebreak
discussion with \bernaysindex\bernays\
on the subject of the layout for the
second edition of \hilbertbernays, \hskip.2em
which probably took place in \Zuerich\
at the very end of the
summer term 1951.\footnotemark\pagebreak

\footnotetext{
  For estimating the time of the detailed discussion with \bernays, \hskip.1em
  see also \cite{scholz-1951d} \hskip.1em
  in addition to \cite{hasenjaeger-1951}. \hskip.5em
  In the letter \cite{hasenjaeger-1951}, \hskip.2em
  a very rough previous sketch is mentioned
  to have been brought into accordance with the elaborated sketch attached to
  \cite{hasenjaeger-1951}. \hskip.3em
  This very rough sketch is probably a
  protocol of a very first short discussion with \bernaysindex\bernays\ on the
  \hilbertbernays\ subject, \hskip.1em
  which took place at the end of winter term 1950/51 in \Zuerich.%
}%
During \hasenjaeger's first two semesters in \Zuerich, \hskip.1em
there were neither publications under his name
nor work on the ``\grundlagendermathematik\closequotefullstopextraspace
Therefore, a further stay of \hskip.1em\hasenjaeger\ as \bernaysindex\bernays' assistant
was strongly suggested by \hskip.05em\scholzname\ \hskip.1em
already February\,1, 1951,
\hskip.2em in \hskip.05em\scholz' letter \shortcite{scholz-1951b} \hskip.1em
to \hasenjaeger. \hskip.5em
This third stay indeed took place during winter term 1951/52 \hskip.1em at the
\ETHZshort, \hskip.1em
and \hasenjaeger\ left \Zuerich\
soon after February\,18, 1952.\footnote{\majorheadroom
  \Cfnlb\
  \makeaciteoftwo
    {hasenjaeger-1952d}
    {hasenjaeger-1970b}.%
}\end{sloppypar}

Most surprisingly, \hskip.1em
however, \hskip.1em
this first proper sketch attached to \cite{hasenjaeger-1951}
  {\em does not show the slightest similarity
    with our typescript here.} \hskip.3em
As the non-matching of this sketch (handwritten by \hasenjaeger) \hskip.1em
with our typescript 
is crucial for the timing assessment of our typescript, \hskip.2em
let us quote the ``Introduction'' (``Einleitung'')
right at the beginning of it, following the original line breaking:

\halftop\halftop\noindent\LINEnomath
    {{\bf ``Die logische Struktur der mathematischen Theorien.}\nlnomath
      Gliederung\esi entwurf.}
\par\noindent\mbox{~~~}{\bf Einleitung}~(E)
\begin{enumerate}\notop\item[a)]%
  Endliche Bereiche. Das finite Schliessen in Bezug auf \hskip.1em``offene''
  Gesamt\-heiten\linebreak
  (Spezies). Endliche Zahlen -- Ordnung\esi zahlen.
  Unendliche Bereiche. 
  Beweise\linebreak al\es\ mathematische Objekte.\noitem\item[b)]\sloppy%
  Elementare Verband\es- und/oder Gruppen-theorie.
  Dabei Au\fsi uchen der einfachsten \germanschlussweisen. \hskip.3em
  Au\fsi tellung eine\es\ provisorischen Au\esi sagen-\linebreak Kalk\ue l\es\
  (\evtl\ Dingvariablen \germanzunaechst\ nur als
  \germanmitteilung\esi zeichen; \linebreak
  dann \Ue berlegung zum) elementaren Kalk\ue l
  mit freien Variablen.\noitem\item[c)]%
  Kombinatorischer Reichtum der
  ``\germanbeweistheorie\closequotefullstopextraspace
  Einbettung der Beweis-\linebreak theorie in die \germanzahlentheorie\
  (Arithmetisierung). \hskip.3em
  ?~Elementare Probleme,\linebreak
  die nichtelementare Methoden erfordern; \ 
  Vergleich mit \germanzahlentheorie.\noitem\item[d)]\sloppy%
  \germanwiderspruchsfreiheit. \hskip.3em
  Behandlung durch endliche Modelle und\linebreak
  durch\,\germanentscheidungsverfahren. \hskip.1em
  \germanentscheidungsverfahren\,an\,sich\,%
  \mbox{(Beisp.Verb.\-theorie).{\bf ''}}\end{enumerate}

\halftop\halftop\noindent
English translation:

\halftop\halftop\noindent\LINEnomath
{{\bf ``The logical structure of mathematical theories.}\nlnomath
Layout outline.}
\par\noindent\mbox{~~~}{\bf Introduction}~(E)
\begin{enumerate}\notop\item[a)]%
  Finite \englishbereiche.
  \englishFinit\ \englishschliessen\ \wrt\ ``open''
  \englishgesamtheiten\ (species).
  Finite numbers --- ordinal numbers.
  Infinite \englishbereiche.
  Proofs as mathematical objects.\noitem\item[b)]\sloppy%
  Elementary lattice and/or group theory.
  Thereby exploration of the most simple \englishschlussweisen. \hskip.3em
  Establishing a provisional \englishaussagenkalkul\
  (maybe object variables \englishzunaechstfirst\ only as
  \englishzeichenpluralzurmitteilungnoindex;
  then \englishueberlegungsingularasplural\
  on the) elementary calculus with free variables.\pagebreak\noitem\item[c)]%
  Combinatorial abundance of
  ``\englishbeweistheorie\closequotefullstopextraspace
  Embedding of \englishbeweistheorie\ into \englishzahlentheorie\
  (arithmetization). \hskip.3em
  ?~Elementary problems that require non-elementary methods;
  comparison to \englishzahlentheorie.\noitem\item[d)]\sloppy%
  \englishWiderspruchsfreiheit. \hskip.3em
  Treatment by means of finite models and
  \englishentscheidungsverfahrenplural. \hskip.3em
  \englishEntscheidungsverfahrenplural\ themselves
  (example lattice theory).{\bf ''}\end{enumerate}

\halftop\noindent
Moreover, this sketch attached to \cite{hasenjaeger-1951} \hskip.1em
shows that the goal of the joint work of
\hasenjaegerindex\hasenjaeger\ and \bernaysindex\bernays\
was 
{\em a complete rewriting or even a completely new writing of the book right from the beginning}\/
\ ---~~definitely not just an integration of new parts into the first edition as
found in our typescript and as indicated in the
``Preface to the Second Edition\closequotecommasmallextraspace
\cfnlb\
\sectref{section The Mentioning of the Given-Up Work on the Second Edition}.

From the sequence of letters and postcards
\cite{hasenjaeger-1952e},
\cite{bernays-1952},
\cite{hasenjaeger-1952f}, \hskip.3em
it becomes clear that the project of a complete reorganization of the book
was still continued in 1952; \hskip.4em
\mbox{in particular}
a completely new version of the treatment of the \math\iota-operator
is discussed in this sequence of letters. \hskip.4em
Moreover,
\hasenjaeger's\linebreak 
state-of-the-art introduction of the \math\iota-operator
\shortcite{hasenjaeger-1952h}
---~ready for press and fundamentally different from the obsolete treatment in both \hilbertbernays\ editions~---
was probably already attached to the first letter \cite{hasenjaeger-1952e}. \hskip.4em
Note that, \hskip.1em
according to the sketch attached to 
\cite{hasenjaeger-1951}, \hskip.3em
the \nlbmath\iota\ was to be treated in \litchaprefs{VI}{VII}
of all in all 12~chapters
(\incl\ the unnumbered introduction and the supplement chapters).

All in all, \hskip.1em
{\em \bernaysindex\bernaysplain\ may give us a wrong impression}\/ \hskip.2em
by using the subjunctive \hskip.1em
``er\-fordert h\ae tte'' \hskip.1em(``would have required'') \hskip.25em
in the \hskip.2em``Preface to the Second Edition\closequotecommaextraspace
\cfnlb\ \sectref{section The Mentioning of the Given-Up Work on the Second Edition}. \hskip.5em
Indeed, \hskip.1em
as \nolinebreak a complete rewriting of the book was obviously planned from the
very beginning of the work in summer~1951 and
pursued far
beyond the end of \hasenjaeger's last stay
in \Zuerich\ in February\,1952, \hskip.3em
it cannot be the case that the planned work was {\em reduced}\/ when\linebreak
---~``already back then'' when
``\MrUS\shorth\hasenjaegerindex\hasenjaegershortname\
came to \Zuerich\ for some time''~---\linebreak
``it became obvious that the integration of the many new results in the area 
of \englishbeweistheorie\ would have required
a complete reorganization of the book.'' \hskip.5em
Of course, \hskip.25em
\bernays\ does not explicitly state that the work on a
complete reorganization was {\em reduced}\/ \hskip.15em
``already back then\closequotecommaextraspace
but his directly following sentence already speaks only of the {\em reduced}\/
version of the second edition
---~without mentioning {\em any other point in time}\/ for the decision on this reduction, \hskip.1em
\cfnlb\ \sectref{section The Mentioning of the Given-Up Work on the Second Edition}. \hskip.5em
Considering \bernays\ most contextual style of writing, \hskip.1em
where a sentence often can only be understood correctly
if the reader has the meaning
of many sentences of the context completely in mind,\footnotemark\ \hskip.1em
we never had the slightest doubt on the identity of these two points in time
for
several
decades, \hskip.1em
until we started our investigations on the joint work of \hasenjaeger\
and \bernays.

The truth about the joint work of \hskip.1em
\hasenjaegerindex\hasenjaeger\ and \bernaysindex\bernays\
on the \hskip.1em``\grundlagendermathematik'' \hskip.05em
seems to be that this work was given up 
only after a completely reorganized version
for the first volume was already written, \hskip.2em
including {\em two}\/ completely new chapters on 
the\,\,\nlbmath\iota-operator
(among other subjects).
\hskip.4em
Be reminded that the \math\iota-operator is treated in both published editions
at the very end of the first volume
only in a {\em single}\/ chapter, \hskip.1em
which does not\linebreak even meet the state of the art of
the original treatment of the \nlbmath\varepsilon\
in the first edition of~1939.

In \nolinebreak addition to this huge work on the entire first volume, \hskip.1em
the joint work of \hasenjaeger\ and \bernays\ also includes the writing of
our typescript (not \atall\ following the ``Introduction'' presented above)
\hskip.1em and, moreover, at least the still lost continuation
of its \nlblitsectref 3.\vfill\cleardoublepage
\section{Conclusion}\halftop
\subsection{Concluding Assessment of Our Typescript}%
\halftop\noindent\footnotetext{
  To the best of our knowledge, \hskip.1em
  the only German author who clearly tops \bernaysname\ \hskip.1em
  in this aspect
  is \jeanpaulindexindexsee\jeanpaulname.%
}%
In spite of some minor evidence that our typescript was written
before\,1929 (\cfnlb\ \sectref{section Parallel introductory part II}(1))
\hskip.1em or in the year\,1931
(\cfnlb\ \sectref{section Parallel introductory part II}(2)), \hskip.3em
the overwhelming evidence says
that it was written not before\,1951
(as explicated in \sectrefs{section The completely new introduction}
{2 of Our Typescript}). \hskip.5em
In \sectref{section Joint Work}, \hskip.2em
we then studied the joint work of
\bernaysindex\bernays\ and \hasenjaegerindex\hasenjaeger\ in \Zuerich\
1950--1952, \hskip.2em
where we found some further evidence that our typescript resulted
from that joint work. \hskip.4em
After all, \hskip.2em
the strongest evidence we have is the following: \hskip.3em
Two carbon copies of our typescript
were found 
in \hasenjaeger's legacy\linebreak
---~together with the only source for the footnotes to our typescript! \hskip.6em
As \hasenjaeger\ seems to have added these footnotes to a typescript
corrected by
\bernays' hand with
marks indicating the places of the anchors for these footnotes, \hskip.2em
\hasenjaeger\ must have been more or less
involved in the production of our typescript
---~and, \hskip.1em
for this and several other reasons
(\cfnlb\ \sectrefs{section bernays-no-typewriter}
{section additional footnotes}), \hskip.25em
he is most probably the actual typist of it \aswell. \hskip.2em

What we still do not know, however, is the time of writing, \hskip.1em
and we cannot exclude any year from 1951 to 1968
(when the second edition of the first volume was published). \hskip.3em
Because of the improvements of our typescript over both editions
of \hilbertbernays\
(\cfnlb\
\sectref{subsection A significant improvement compared to both editions},
\sectref{section Parallel introductory part II}(3),
\sectref{2 of Our Typescript}), \hskip.2em
even the years from 1968 to 1977 (when \bernaysindex\bernays\ died) cannot be
excluded. \hskip.3em
Most likely, however, is the time
of the winter term 1951/52, \hskip.2em
the last one \hasenjaeger\ stayed with
\bernaysindex\bernays\ at the \ETHshort\ in \Zuerich\
(as explicated in \sectref{section Joint Work}).

In any case, \hskip.1em
the finding of this incorrectly filed typescript
is essential in the context of \hilbertindex\hilbert\ and
\bernaysindex\bernays' ``\grundlagendermathematik\closequotecommasmallextraspace
because it provides us with some new
\mbox{insights}\footnotemark\footnotetext{\majorheadroom%
  For instance, \hskip.2em
  \litsectref 2 of our typescript has provided us with some new insight into
  the \englishfinit\ \englishstandpunkt,
  which helped us to improve the translation of the first part of \hskip.1em
  \litsectref 2
  of \makeaciteoftwo{grundlagen-first-edition-volume-one}
  {grundlagen-second-edition-volume-one}
  in \cite{grundlagen-german-english-edition-volume-one-one}
  \hskip.05em considerably, \hskip.1em
  \cfnlb\ \sectref
  {subsection A significant improvement compared to both editions}.\par
  Moreover ---~and maybe more important~---
  also \litsectref 3 of our typescript has provided \nolinebreak us
  with some new insight into
  the \englishfinit\ \englishstandpunkt,
  which led to \litnoteref{349.2} in 
  \cite[\p\,349]{grundlagen-german-english-edition-volume-one-three}.%
} 
into the \englishfinit\ \englishstandpunkt\
and its development, \hskip.2em
in particular in connection with the two editions of the {\em first volume}\/ of
``\grundlagendermathematik'' \hskip.05em \makeaciteoftwo
{grundlagen-first-edition-volume-one}
{grundlagen-second-edition-volume-one}.
\halftop\subsection{Three Most Interesting Scripts Still Missing}\halftop\noindent
As we do not know any way to find the following three scripts
with our limited resources,
the treasure quest remains open for future prospectors:\begin{enumerate}\item
It is worthwhile to invest further effort and to search
the legacy of \bernaysindex\bernays\ and the \ETHshort\
archives more broadly: \hskip.3em
There may be a chance to find an incorrectly filed continuation of
our typescript, \hskip.1em
at least the remainder of its \litsectref 3, \hskip.2em
possibly in form of a draft partly 
in \bernaysindex\bernays' \englishgabelsbergerkurzschrift.\par
Such a finding may change our point of view on \bernaysindex\bernays' ideas
on the foundations of mathematics considerably.\yestop\item
Be aware that our typescript must not be mistaken for
the {\em other script}\/ \hskip.05em
that \hasenjaeger\ wrote according to his handwritten sketch
attached to \cite{hasenjaeger-1951} \ 
---~~probably a typescript of a completely rewritten
first volume of \hskip.25em ``\grundlagendermathematik\closequotecommaextraspace
\cfnlb\ \sectref{section Joint Work}.

This other script would be a real bonanza
regarding the views of \bernaysindex\bernays\ on the foundations of mathematics
in the early 1950s. \par

We have no idea on the whereabouts of this other script.\par
  
It would be a real pity if this major work of \hasenjaeger\
and \bernaysindex\bernays,                                      \hskip.2em
which can hardly be overlooked in any library by its mere size, \hskip.2em
really remained lost.

\yestop\item
Finally, \hskip.1em
to find the script for the first edition remains, \hskip.1em
of course, \hskip.1em
one of the
biggest wishes of maybe every historian of modern logic.\end{enumerate}
\vfill\vfill\vfill\vfill\section*{Acknowledgments}%
\addcontentsline{toc}{section}{Acknowledgments}%
We would like to thank
\beatebeckername\ (\nee\ \hasenjaeger),
\bernaysludwigname, \glaschickname,
\mancosuname, \schmehname, and \siegname\ for their most helpful comments,
and \hungerbuehlername\ and \stolzenburgname\ for their kind support.
\pagebreak

\section*{Curriculum Vitae of \hasenjaegernameplain\ by Himself}\bernaysindex
\addcontentsline{toc}{section}{Curriculum Vitae of \hasenjaegernameplain\ by Himself}%
\subsection*{German Original}\Zuerichindex
The following is a reprint of \cite{hasenjaeger-1997},
written by \hasenjaegername\ himself. \hskip.3em
Also the layout is exactly the one of \hasenjaeger's ASCII file.

{\setlength{\baselineskip}{11pt}\footnotesize\tt\mbox{}

\par\noindent
\input{VITA_GH3.tex}
}

\subsection*{English translation}\Zuerichindex\hasenjaegerindex

\par\halftop\noindent Date of Generation: \Aug\,16, 1997
\\\mbox{}
\\\mbox{}
\\\mbox{}~~~~~~~~~~~~~~~~VITA~~~~Gisbert~HASENJAEGER,~born~\May\,1, 1919, in \Hildesheim
\\\mbox{}
\\\mbox{~~~~~~~~}[REM\edcomment{ark?:}
  I agree with corrections in the sense of standardization and abbreviation,
  \\\mbox{~~~~~~~~~~~~~~~~~~~~~~~~~}but ask for communication of the result]
\\\mbox{}
\\\mbox{~~~~~~~~}1937~A-levels~(Abitur)~in~M\ue lheim~a.d.~Ruhr.
\\\mbox{~~~~~~~~}1937--1939~Working and Military Service.
\\\mbox{~~~~~~~~}1942~After~shot in the head~(\Jan\,2)~and~reconvalensence~at~OKW/CHI
\\\mbox{~~~~~~~~~~~~~~~}mission against~A.\,Turing  (as I was told much later).
\\\mbox{~~~~~~~~}1945--1950~Studies~in~M\ue nster/Westf.,~from the beginning
Scientific Assistent at the
\\\mbox{~~~~~~~~~~~~~~~~~~~~~~\,}Seminar/Institute for Mathematical Logic and Foundational Research
\\\mbox{~~~~~~~~}1950~Promotion~to~\Drrernat~by~\scholznameplain~with the~dissertation.
\\\mbox{~~~~~~~~~~~~~~~~}``\hasenjaegernineteenfiftyctitleenglishwithlinebreak\closequotefullstopnospace
\\\mbox{~~~~~~~~}1950~Postdoctoral Research Assistent at the above-mentioned institute.
\\\mbox{~~~~~~~~}1950/51~(the~academic~year)~Scholarship Awardee as a guest of \bernaysnameplain\ at the
\\\mbox{~~~~~~~~~~~~~~~~~~~~~~~~~~~~~~~~~~~~~~~~~~~~~~}\ETHZshort.
\\\mbox{~~~~~~~~}1953~Habilitation~and~venia~legendi~(mathematical~logic~and~foundational research)
\\\mbox{~~~~~~~~~~~~~~~}at the \Univ\,\MuensternoWestfalen
\\\mbox{~~~~~~~~~~~~~~~}(habilitation thesis: ``\hasenjaegerhabiltitleenglish'').
\\\mbox{~~~~~~~~}1955~Position of Supernumerary University Lecturer,
1960~Adjunct~Professor,
\\\mbox{~~~~~~~~~~~~~~~}\Univ\,\MuensternoWestfalen
\\\mbox{~~~~~~~~}1961/62\,\,(summer and winter term)~~Deputy for a newly founded logic professorship
\\\mbox{~~~~~~~~~~~~~~~~~~~~~~~~~~~~~~~~~~~~~~~~~~~~~~~~~~~~~~~~}at the philosophical faculty at \Univ\,Bonn.
\\\mbox{~~~~~~~~}1962~Appointed to this professorschip, for logic and foundational research.
\\\mbox{~~~~~~~~~~~~~~~}[REM~As the adjective ``math.''\ is lacking here and according to today's
\\\mbox{~~~~~~~~~~~~~~~~}\englishsprachgebrauch, this seems to be a
``wide field'' nowadays, which I
\\\mbox{~~~~~~~~~~~~~~~~}would definitely not want to treat anymore.]
\\\mbox{~~~~~~~~}1964~Personal Professor (``Pers\oe nlicher~Ordinarius''),
\\\mbox{~~~~~~~~}1966~Full Professor (``Ordinarius'').
\\\mbox{~~~~~~~~}1964/5~(the academic year)~Guest at the Institut~for~Advanced~Studies
\\\mbox{~~~~~~~~~~~~~~~~~~~~~~~~~~~~~~~~~~~~~~~~~~~~~}Princeton~N.J.
\\\mbox{~~~~~~~~1970/1}~(the academic year) Guest Professorship at the \Univ\ of Illinois in
\\\mbox{~~~~~~~~~~~~~~~~~~~~~~~~~~~~~~~~~~~~~~~~~~~~~}Urbana/Champaign.
\\\mbox{~~~~~~~~}1984~Professor Emeritus, \Univ\ Bonn.
\vfill\cleardoublepage
\section*{Our Curriculum Vitae of\/ \hasenjaegernameplain}\label{section short CV GH}\Zuerichindex\beatebeckerindex
\addcontentsline{toc}{section}{Our Curriculum Vitae of \hasenjaegernameplain}

To provide some more detail and
to overcome errors in other publications,
we carefully compiled the following Curriculum Vitae.\footnotemark

\vfill\noindent\LINEnomath{\begin{tabular}{@{}l@{~\,\,}l@{}}1919
&Born in \Hildesheim, \Jun\,1. 
\\1937
 &A-levels (``Abitur'') in \Muelheim.
\\1937--1939
&Involuntary\footnotemark\
Working and Military Service (Reich\esi arbeit\es- und Wehrdienst).
\\1939&Drafted for Military Service in World War II.
\\1941--1942&Artillerist in the German Attack of Russia (``Russlandfeldzug'').
\\1942
&Shot through his helmet to the head as a German soldier in Russia, \Jan\,2.
\\1942--1945
&Cryptologist at the High Command of the German Armed Forces in \Berlin.
\\1945--1950
&Studies in Mathematics and Physics, \Univ\ of \Muenster.
\\1945--1950
&
Scientific
Assistant (``wissenschaftl.\,Hilfskraft''), \hskip.1em
``Seminar/Institut f\ue r 
\\&\hfill Mathematische Logik und Grundlagenforschung'' \hskip.1em
   at the \Univ\,of \MuensternoWestfalen.
\\1950
&Promoted to \Drrernat\ (\PhD) \cite{hasenjaeger-1950c}
\\&\hfill by \scholzname\ \scholzlifetime.
\\1950
&Postdoctoral Research Assistant (``wissenschaftl.\,Assistent''), \hskip.1em
``Institut f\ue r 
\\&\hfill Mathematische Logik und Grundlagenforschung'' \hskip.1em
   at the \Univ\,of \MuensternoWestfalen.
\\1950--1952
&Three Scholarships in winter term 1950/51, summer 1951,
winter 1951/52,\footnotemark 
\\&\hfill granted by the \ETHshort\ Board,
    for visiting the \ETHZshort\ to assist \bernaysindex\bernays
\\&\hfill in the writing of the second edition of \hskip.15em
``\grundlagendermathematik\closequotefullstopnospace
\\1953
&Habilitation \cite{hasenjaeger-1953b}, \hskip.2em and 
{\it venia legendi}\/ (Mathematische
\\&\hfill  Logik und Grundlagenforschung) \hskip.1em
    at the University of \MuensternoWestfalen.
\\1955--1960
&Lecturer (``Di\ae tendozentur'') \hskip.1em
at the University of \MuensternoWestfalen.
\\1956
&Married {\namefont Irmhild Reinl\ae nder}
(a former student assistant of \scholz,
\\&\hfill mathematics \& physics teacher in Soest (Germany)\,1955--1957,
died 2012). 
\\1957
&Son {\namefont\andreas} is born. 
\\1959
&Daughter {\namefont Beate} is born. 
\\1960--1961
&Associate Professor (``apl.\,\Prof'')\ \hskip.1em
at the University of \MuensternoWestfalen.
\\1961
&Daughter {\namefont Cordula} is born. 
\\1961--1962
&Substitute Professor in summer term 1961 and winter term 1961/62
\\&\LINEnomath{at a newly founded professorship for logic at the philosophical faculty}
\\&\hfill of the  University of Bonn (Germany).
\\1962--1984
&Associate Professor (``a.\,o.\,\Prof'') on this professorship,
\\&\hfill ``Pers\oe nlicher  Ordinarius'' 1964,
Full Professor (``Ordinarius'') 1966.
\\
&Director of the newly founded ``Seminar f\ue r Logik und
\\&\hfill
Grundlagenforschung'' {\em at the  University of  Bonn}.
\\1964--1965
&Visitor at the \IAS,
\\&\hfill winter term 1964/65 and summer term 1965 (with his family).
\\1970--1971
&Visiting Professor at the \UIUC
\\&\hfill winter term 1970/71 and summer term 1971 (with his family).
\\1984--2006
&Professor Emeritus at the University of Bonn.
\\2006
&Died on the estate of his wife in Plettenberg (Germany), \Sep\,2.
\\\end{tabular}}
\vfill\noindent
\pagebreak
\addtocounter{footnote}{-2}\footnotetext{%
  Sources:
  Correspondence \hasenjaeger--\bernaysindex\bernays:
  \makeaciteoffour
  {hasenjaeger-1950d}{hasenjaeger-1950e}{hasenjaeger-1952e}{hasenjaeger-1952f}.
  \hskip.3em
  Three Curricula Vitae by \hasenjaeger\ himself:
  \makeaciteofthree
      {hasenjaeger-1952g}
      {hasenjaeger-1970b}
      {hasenjaeger-1997}. \hskip.3em
  Published \mbox{articles:}
  \makeaciteofthree
      {schmeh-hasenjaeger-0}
      {schmeh-hasenjaeger-1}
      {schmeh-hasenjaeger-2}.
  \hskip.3em
  Unpublished laudation: \cite{diller-on-hasenjaeger}. \hskip.3em
  Communication with \hasenjaegername's daughter \beatebeckername\ by \EMAIL\ and by phone,
  \Jan--\Feb\,2018 and \Jan\,2020.%
}\addtocounter{footnote}{1}\footnotetext{%
  \majorheadroom
  Contrary to a text 
  that wrongly puts \hasenjaeger\ close to being an admirer of the Nazi reign of terror,
  found on \Dec\,28,~2019, in the German,
  English and French Wikipedia entries on \hasenjaegername\ and also 
  in the many places that automatically
  copy this attack to the soundness of public knowledge,
  \hasenjaeger\ did not do these services {\em voluntarily}.
  For instance, in the German version we read:\par
  ``Hasenjaeger war danach freiwillig beim Reichsarbeitsdienst und leistete seinen Wehrdienst.''\par
  Having completed his \nth{18}~year by \Jun\,1,\,1937, \hasenjaeger\
  was forced by law
  to do 6~months service
  according to the Reich\esi arbeit\esi dienst\-gesetz of \Jun\,26, 1935,
  immediately followed by 24~months service
  according to the Wehrgesetz of \May\,21, 1935. \hskip.3em
  So he should have been drafted by the Nazi state by \Jul\,1937,
  unless something relevant for the Nazi regime spoke against it,
  such as activities relevant for warfare,
  which included {\em studying medicine}. \hskip.1em
  According to these laws,
  the length of these services
  are to be set by the ``F\ue hrer and Reich\esi kanzler\closequotecomma
  \ie\ by \hitlername, who chose these times as given here.\par
  A possible source for the claim of a {\em voluntary}\/ participation
  is a new link on 
  \Jan\,10, 2020, in the German Wikipedia, pointing to \cite{schmeh-hasenjaeger-0},
  where we read:
  ``Nach dem Abitur im Jahr 1936 meldete er sich zun\ae chst freiwillig zum Arbeitsdienst,
  um anschlie\sz end studieren zu k\oe nnen.
  Da jedoch der Krieg dazwischen kam,
  wurde er zum Milit\ae rdienst eingezogen.'' \hskip.2em
  In \cite[\p\,343]{schmeh-hasenjaeger-1} we find an English version of these sentences:
  ``After his graduation in 1936,
  he volunteered for labor service, planning to take up his studies afterwards.
  But the war intervened and he was drafted.''\par
  Several aspects should be noted here: \hskip.3em
  (1)~The sentences are not very reliable as
  \hasenjaeger\ writes in his own CV \cite{hasenjaeger-1997} that he graduated in 1937,
  which is perfectly in line with the German standard of entering school with 6
  (\ie\ in late summer\,1925) and graduating 12 years later from 1937 on (before\,1937: 13 years).
  \hskip.3em
  (2)~The war intervened not during \hasenjaeger's labor service,
      but during his later military service.
  (3)~In \cite[\p\,343]{schmeh-hasenjaeger-1} we read:
  ``The author met him in 2005, a year before his death.''
  Therefore, the German sentence in \cite{schmeh-hasenjaeger-0} is closer
  to the communication of \schmehname\ with \hasenjaeger\ in~2005
  and therefore more reliable. \hskip.3em
  (4)~This is relevant because the German version actually says something different:
  Not the labor service was
  voluntary, but his ``Meldung'' (``answer'' or ``notice'').
  Thus,
  after passing his obligatory muster (probably in late\,1936),
  \hasenjaeger\ probably answered 
  that he would prefer an early start
  of his labor service because he wanted
  to study mathematics
  (but not medicine!)\ \hskip.1em
  after his labor and military services
  {\em without any interruption of his studies}.%
}\addtocounter{footnote}{1}\footnotetext{%
  \majorheadroom
  In his own CV \cite{hasenjaeger-1997},
  \hasenjaeger\ omits the last semester of his grants in \Zuerich.
  In \nolinebreak his draft for an English CV \cite{hasenjaeger-1970b}, however,
  he writes:
  ``From autumn\,1950 to spring\,1952 I worked under a scholarship with \bernaysindex\bernaysshortnameplain\
    at the \ETHshort, Switzerland.'' \par
  Up to now,
  we can only be certain that \hasenjaeger\ was in \Zuerich\ for the first 18 days of 
  February\,1952, \ie\ during the very end of the winter term\,1951/52,
  and that he planned to meet \scholz\ in \MuensternoGermany\ not before \Feb\,29, 1952.
  The hard evidence for this is the letter \cite{hasenjaeger-1952d} to \scholz,
  written on \Feb\,18, 2012, in \Zuerich.%
}%
\vfill\cleardoublepage

\nocite{hasenjaeger-1968,hasenjaeger-1976c,bernays-Hs973:41,hasenjaeger-1998,hasenjaeger-1984c,hasenjaeger-1995,hasenjaeger-1976b,hasenjaeger-1976a,hasenjaeger-1990,hasenjaeger-1972b,hasenjaeger-1987,hasenjaeger-1977,hasenjaeger-1968a,hasenjaeger-1967,hasenjaeger-1966b,hasenjaeger-1966a,hasenjaeger-1965,hasenjaeger-1960,hasenjaeger-1959,hasenjaeger-1955,behnke-hasenjaeger-1955,hasenjaeger-1953/54,bernays-questions-methodologiques-actuelles-german-manuscript,bernays-1958,bernays-axiomatic-set-theory-I,bernays-axiomatic-set-theory-II,bernays-axiomatic-set-theory-III,bernays-axiomatic-set-theory-IV,bernays-axiomatic-set-theory-V,bernays-axiomatic-set-theory-VI,scholz-hasenjaeger,hasenjaeger-1962,scholz-Logistik,scholz-geschichte-der-Logik,hasenjaeger-1958a,hasenjaeger-1958b,hasenjaeger-1978,hasenjaeger-1984a,hasenjaeger-1984b}

\def\refinitialtext{\begin{sloppypar}%
  For a hopefully complete list of \hasenjaegername's scientific
  and scholarly publications, \hskip.2em
  check the items with the labels of the forms \hskip.2em
  \cite{behnke-hasenjaeger-1955}, \hskip.25em
  \cite{logic-and-machines}, \hskip.25em
  \startcite\hasenjaeger,~\ldots\stopcite, \hskip.25em
  \cite{hasenjaeger-thyssen-1968}, \hskip.25em
  \cite{scholz-hasenjaeger}.\end{sloppypar}
}

\catcode`\@=11
\renewcommand\@openbib@code{%
      \advance\leftmargin\bibindent
      \itemindent -\bibindent
      \listparindent \itemindent
      \parsep 12.5pt
      \itemsep 0pt
      \baselineskip 13.6pt
      }%
\catcode`\@=12
\addcontentsline{toc}{section}{\refname} 
\bibliography{herbrandbib}\cleardoublepage
\addcontentsline{toc}{section}{Index} 
\printindex
\end{document}

%% file: namedhelper.tex
\def\citep{\cite}
\def\citet#1{\citeauthor{#1} \shortcite{#1}}
\newcommand\startcite{{\raise.2ex\hbox{[}}}
\newcommand\stopcite {\raise.2ex\hbox{]}}
\newcommand\citehelper[1]{\startcite #1\stopcite}
\newcommand\makeaciteoftwo[2]
{\citehelper{\citeauthor{#1}, \citeyear{#1}; \citeyear{#2}}}
\newcommand\makeacitetoftwo[2]
{\citeauthor{#1} \citehelper{\citeyear{#1}; \citeyear{#2}}}
\newcommand\makeaciteofthree[3]
{\citehelper{\citeauthor{#1}, \citeyear{#1}; \citeyear{#2}; \citeyear{#3}}}

\newcommand\makeaciteoffour[4]
{\citehelper{\citeauthor{#1}, \citeyear{#1}; \citeyear{#2}; \citeyear{#3};
\citeyear{#4}}}

\newcommand\makeaciteoffive[5]
{\citehelper{\citeauthor{#1}, \citeyear{#1}; \citeyear{#2}; \citeyear{#3};
\citeyear{#4}; \citeyear{#5}}}

%% file: headernames.tex
\usepackage{amssymb}

\usepackage[T1]{fontenc}

\input specialfonts

\def\mathfrak#1{%
\mathchoice
{{\mfrak{#1}}}
{{\mfrak{#1}}}
{{\mscriptfrak{#1}}}
{{\mscriptscriptfrak{#1}}}
}

\newcommand\mathfrakfootnote[1]{%
\mathchoice
{{\mfootnotefrak{#1}}}
{{\mfootnotefrak{#1}}}
{{\mscriptscriptscriptfrak{#1}}}
{{\mscriptscriptscriptfrak{#1}}}
}

\input headernamesrest

%% file: specialfonts.tex
\newcommand\actualpagebreaknospace{{\rm\textbar}}
\newcommand\actualpagebreaknospacepagenumber
                                     [1]{\actualpagebreaknospace\math{_{#1}}}
\newcommand\actualpagebreak             {\actualpagebreaknospace              \ }
\newcommand\actualpagebreakpagenumber[1]{\actualpagebreaknospacepagenumber{#1}\ }

\newcommand\actualparpagebreakpagenumber[1]
                {\getittotheright{\actualpagebreaknospacepagenumber{#1}\hskip-.5em}}

\newcommand\actualpagebreakindentpagenumber[1]
                   {\noindent\actualpagebreaknospacepagenumber{#1}\hskip-1.1em\indent}

\mathchardef\Gammaoffont="7000
\mathchardef\Gamma="0100
\mathchardef\Deltaoffont="7001
\mathchardef\Delta="0101
\mathchardef\Thetaoffont="7002
\mathchardef\Theta="0102
\mathchardef\Lambdaoffont="7003
\mathchardef\Lambda="0103
\mathchardef\Xioffont="7004
\mathchardef\Xi="0104
\mathchardef\Pioffont="7005
\mathchardef\Pi="0105
\mathchardef\Sigmaoffont="7006
\mathchardef\Sigma="0106
\mathchardef\Upsilonoffont="7007
\mathchardef\Upsilon="0107
\mathchardef\Phioffont="7008
\mathchardef\Phi="0108
\mathchardef\Psioffont="7009
\mathchardef\Psi="0109
\mathchardef\Omegaoffont="700A
\mathchardef\Omega="010A
\mathchardef\itype="017B

\catcode`\@=11

\gdef\allowhyphens{\penalty\@M \hskip\z@skip}

\gdef\set@low@box#1{\setbox\tw@\hbox{,}\setbox\z@\hbox{#1}\dimen\z@\ht\z@
     \advance\dimen\z@ -\ht\tw@
     \setbox\z@\hbox{\lower\dimen\z@ \box\z@}\ht\z@\ht\tw@ \dp\z@\dp\tw@ }
\gdef\set@low@boxsingle#1{\setbox\tw@\hbox{\rm,}\setbox\z@\hbox{#1}\dimen\z@\ht\z@
     \advance\dimen\z@ -\ht\tw@
     \setbox\z@\hbox{\lower\dimen\z@ \box\z@}\ht\z@\ht\tw@ \dp\z@\dp\tw@ }

\gdef\@glqq{%
\ifhmode\edef\@SF{\spacefactor\the\spacefactor}%
\else\let\@SF\empty
\fi
\CheckFamily\font\fraknomath\ifSameFamily ``\relax
\else\CheckFamily\font\swab\ifSameFamily ``\relax
\else\leavevmode\set@low@box{''}\box\z@\kern-.04em\allowhyphens\@SF\relax
\fi\fi}
\gdef\glqq{\protect\@glqq\kern+.07em}
\gdef\@grqq{%
\ifhmode\edef\@SF{\spacefactor\the\spacefactor}%
\else\let\@SF\empty 
\fi 
\CheckFamily\font\fraknomath\ifSameFamily ''\relax
\else\CheckFamily\font\swab\ifSameFamily ''\relax
\else\kern+.07em``\kern.07em\@SF\relax
\fi\fi}
\gdef\grqq{\protect\@grqq}
\gdef\@glq{{\ifhmode \edef\@SF{\spacefactor\the\spacefactor}\else
     \let\@SF\empty \fi \leavevmode
     \set@low@boxsingle{'\/}\box\z@\kern-.04em\allowhyphens\@SF\relax}}
\gdef\glq{\protect\@glq\kern+.07em}
\gdef\@grq{\ifhmode \edef\@SF{\spacefactor\the\spacefactor}\else
     \let\@SF\empty \fi \kern-.0125em`\kern.07em\@SF\relax}
\gdef\grq{\protect\@grq}

\catcode`\@=12

\newcommand\closequotecommanospace{''\nolinebreak\hskip-0.23em,}
\newcommand\closequotecomma      {\closequotecommanospace\         \,}
\newcommand\closequotecommasmallextraspace{\closequotecommanospace\ \,\,}
\newcommand\closequotecommaextraspace{\closequotecommanospace\   \ \,}

\newcommand\closequotefullstopextraspace   
                                 {\closequotefullstopnospace\    \ \,}

\newcommand\closequotefullstopnospace
                                 {''\nolinebreak\hskip-0.20em\@.}

\newcommand\grqqcommanospace       {\grqq\nolinebreak\hskip-0.25em,}

\newcommand\commanospace           {,\nolinebreak\hskip-0.23em}
\newcommand\fullstopnospace        {\@.\nolinebreak\hskip-0.23em}

\makeatletter
\if \@ptsize 0
   \newfont{\scriptscriptscriptgoth}{ygoth scaled 760}
   \newfont{\scriptscriptgoth}{ygoth scaled 833}
   \newfont{\scriptgoth}{ygoth scaled 912}
   \newfont{\gothnomath}{ygoth}
   \newfont{\Goth}{ygoth scaled \magstephalf}
   \newfont{\GOth}{ygoth scaled \magstep1}
   \newfont{\GOTh}{ygoth scaled \magstep2}
   \newfont{\GOTH}{ygoth scaled \magstep3}

   \newfont{\scriptscriptscriptswab}{yswab scaled 760}
   \newfont{\scriptscriptswab}{yswab scaled 833}
   \newfont{\scriptswab}{yswab scaled 912}
   \newfont{\swab}{yswab}
   \newfont{\Swab}{yswab scaled \magstephalf}
   \newfont{\SWab}{yswab scaled \magstep1}
   \newfont{\SWAb}{yswab scaled \magstep2}
   \newfont{\SWAB}{yswab scaled \magstep3}

   \newfont{\scriptscriptscriptfrak}{yfrak scaled 760}
   \newfont{\scriptscriptfrak}{yfrak scaled 833}
   \newfont{\scriptfrak}{yfrak scaled 912}
   \newfont{\fraknomath}{yfrak}
   \newfont{\Frak}{yfrak scaled \magstephalf}
   \newfont{\FRak}{yfrak scaled \magstep1}
   \newfont{\FRAk}{yfrak scaled \magstep2}
   \newfont{\FRAK}{yfrak scaled \magstep3}

   \newfont{\init}{yinit}
   \newfont{\Init}{yinit scaled \magstephalf}
   \newfont{\INit}{yinit scaled \magstep1}
   \newfont{\INIt}{yinit scaled \magstep2}
   \newfont{\INIT}{yinit scaled \magstep3}
\fi
\if \@ptsize 1
   \newfont{\scriptscriptscriptgoth}{ygoth scaled 833}
   \newfont{\scriptscriptgoth}{ygoth scaled 912}
   \newfont{\scriptgoth}{ygoth}
   \newfont{\gothnomath}{ygoth scaled \magstephalf}
   \newfont{\Goth}{ygoth scaled \magstep1}
   \newfont{\GOth}{ygoth scaled \magstep2}
   \newfont{\GOTh}{ygoth scaled \magstep3}
   \newfont{\GOTH}{ygoth scaled \magstep4}

   \newfont{\scriptscriptscriptswab}{yswab scaled 833}
   \newfont{\scriptscriptswab}{yswab scaled 912}
   \newfont{\scriptswab}{yswab}
   \newfont{\swab}{yswab scaled \magstephalf}
   \newfont{\Swab}{yswab scaled \magstep1}
   \newfont{\SWab}{yswab scaled \magstep2}
   \newfont{\SWAb}{yswab scaled \magstep3}
   \newfont{\SWAB}{yswab scaled \magstep4}

   \newfont{\scriptscriptscriptfrak}{yfrak scaled 833}
   \newfont{\scriptscriptfrak}{yfrak scaled 912}
   \newfont{\scriptfrak}{yfrak}
   \newfont{\fraknomath}{yfrak scaled \magstephalf}
   \newfont{\Frak}{yfrak scaled \magstep1}
   \newfont{\FRak}{yfrak scaled \magstep2}
   \newfont{\FRAk}{yfrak scaled \magstep3}
   \newfont{\FRAK}{yfrak scaled \magstep4}

   \newfont{\init}{yinit scaled \magstephalf}
   \newfont{\Init}{yinit scaled \magstep1}
   \newfont{\INit}{yinit scaled \magstep2}
   \newfont{\INIt}{yinit scaled \magstep3}
   \newfont{\INIT}{yinit scaled \magstep4}
\fi
\if \@ptsize 2
   \newfont{\scriptscriptscriptgoth}{ygoth scaled 912}
   \newfont{\scriptscriptgoth}{ygoth}
   \newfont{\scriptgoth}{ygoth scaled \magstephalf}
   \newfont{\gothnomath}{ygoth scaled \magstep1}
   \newfont{\Goth}{ygoth scaled \magstep2}
   \newfont{\GOth}{ygoth scaled \magstep3}
   \newfont{\GOTh}{ygoth scaled \magstep4}
   \newfont{\GOTH}{ygoth scaled \magstep5}

   \newfont{\scriptscriptscriptswab}{yswab scaled 912}
   \newfont{\scriptscriptswab}{yswab}
   \newfont{\scriptswab}{yswab scaled \magstephalf}
   \newfont{\swab}{yswab scaled \magstep1}
   \newfont{\Swab}{yswab scaled \magstep2}
   \newfont{\SWab}{yswab scaled \magstep3}
   \newfont{\SWAb}{yswab scaled \magstep4}
   \newfont{\SWAB}{yswab scaled \magstep5}

   \newfont{\scriptscriptscriptfrak}{yfrak scaled 833}
   \newfont{\scriptscriptfrak}{yfrak}
   \newfont{\scriptfrak}{yfrak scaled \magstephalf}
   \newfont{\fraknomath}{yfrak scaled \magstep1}
   \newfont{\Frak}{yfrak scaled \magstep2}
   \newfont{\FRak}{yfrak scaled \magstep3}
   \newfont{\FRAk}{yfrak scaled \magstep4}
   \newfont{\FRAK}{yfrak scaled \magstep5}

   \newfont{\init}{yinit scaled \magstep1}
   \newfont{\Init}{yinit scaled \magstep2}
   \newfont{\INit}{yinit scaled \magstep3}
   \newfont{\INIt}{yinit scaled \magstep4}
   \newfont{\INIT}{yinit scaled \magstep5}
\fi

\newcommand{\mscriptscriptscriptfrak}[1]{\mbox{\scriptscriptscriptfrak#1}}
\newcommand{\mscriptscriptfrak}      [1]{\mbox{\scriptscriptscriptfrak#1}}
\newcommand{\mscriptfrak}            [1]{\mbox{\scriptscriptscriptfrak#1}}
\newcommand{\mfootnotefrak}          [1]{\mbox{\scriptscriptfrak#1}}
\newcommand{\mfrak}[1]{\mbox{\fraknomath#1}}

\makeatother

\newif\ifSameFamily
\def\CheckFamily#1#2{\GetFamilyName{#1}\ArgOne
        \GetFamilyName{#2}\ArgTwo
        \ifx\ArgOne\ArgTwo\SameFamilytrue\else\SameFamilyfalse\fi}
\def\GetFamilyName#1{\edef\Tempa{#1}\def\Tempb{#1}\ifx\Tempa\Tempb
        \edef\Tempa{\fontname#1}\fi
        \edef\Tempa{\Tempa\space}%
        \expandafter\iGetFamilyName\Tempa\\}
\def\iGetFamilyName#1 #2\\#3{\def#3{#1}}
\def\DefFontName#1#2{{\escapechar-1\expandafter\expandafter\expandafter
        \iDefFontName\expandafter{\csname#2\endcsname}%
        \xdef#1{\expandafter\string\Tempa}}}
\def\iDefFontName{\def\Tempa}

\newcommand\unprotectedae
{\CheckFamily\font\fraknomath\ifSameFamily *a\else
 \CheckFamily\font\swab\ifSameFamily\char'212\else\"a\fi\fi}
\newcommand\unprotectedoe
{\CheckFamily\font\fraknomath\ifSameFamily 
*o\else\CheckFamily\font\swab\ifSameFamily\char'232\else\"o\fi\fi}
\newcommand\unprotectedue
{\CheckFamily\font\fraknomath\ifSameFamily 
*u\else\CheckFamily\font\swab\ifSameFamily\char'237\else\"u\fi\fi}
\newcommand\unprotectedAe
{\CheckFamily\font\fraknomath\ifSameFamily 
Ae\else\CheckFamily\font\swab\ifSameFamily Ae\else\"A\fi\fi}

\newcommand\unprotectedUe
{\CheckFamily\font\fraknomath\ifSameFamily 
Ue\else\CheckFamily\font\swab\ifSameFamily Ue\else\"U\fi\fi}
\DefFontName\eccclarge{eccc1200}
\DefFontName\eccc{eccc1000}
\DefFontName\ecccsmall{eccc0900}
\DefFontName\ecccfootnotesize{eccc0800}
\newcommand\unprotectedsz
{\CheckFamily\font\fraknomath\ifSameFamily\char'032\else
 \CheckFamily\font\swab\ifSameFamily\char'032\else
 \CheckFamily\font\eccclarge\ifSameFamily s
 z\else 
 \CheckFamily\font\eccc\ifSameFamily s
 z\else
 \ss\fi\fi\fi\fi
}
\newcommand\unprotectedes
{\CheckFamily\font\fraknomath\ifSameFamily\char'215\else
\CheckFamily\font\swab\ifSameFamily\char'215\else  
s\fi\fi}

\newcommand\unprotectedesi
{\CheckFamily\font\fraknomath\ifSameFamily\char'215\else
\CheckFamily\font\swab\ifSameFamily\char'215\else  
\mbox{s\hskip.04em}\fi\fi
\-\ignorespaces}

\newcommand\unprotectedsesi
{\CheckFamily\font\fraknomath\ifSameFamily 
s\nolinebreak\hskip-.055em\char'215\else
\CheckFamily\font\swab\ifSameFamily s\char'215\else\mbox{ss\hskip.041em}\fi\fi
\-\ignorespaces}
\newcommand\unprotectedmyparagraphsymbol
{\CheckFamily\font\fraknomath\ifSameFamily 
\char'244\else\CheckFamily\font\swab\ifSameFamily
\char'244\else\S\fi\fi}

\renewcommand\ae{\protect\unprotectedae}
\renewcommand\oe{\protect\unprotectedoe}
\newcommand\ue  {\protect\unprotectedue}
\newcommand\Ae  {\protect\unprotectedAe}

\newcommand\Ue  {\protect\unprotectedUe}
\newcommand\sz  {\protect\unprotectedsz}
\newcommand\es  {\protect\unprotectedes}

\newcommand\esi {\protect\unprotectedesi}  
\newcommand\sesi{\protect\unprotectedsesi} 
 
\newcommand\fbi {\discretionary{f-} {b}{f\mbox {\hskip.07em}b}}

\newcommand\fgi {\discretionary{f-} {g}{f\mbox {\hskip.05em}g}}

\newcommand\fli {\discretionary{f-} {l}{f\mbox {\hskip.07em}l}}

\newcommand\fri {\discretionary{f-} {r}{f\mbox {\hskip.09em}r}} 
\newcommand\fsi {\discretionary{f-} {s}{f\mbox {\hskip.08em}s}}

\newcommand\lli {\discretionary{l-} {l}{l\mbox {\hskip.07em}l}} 
 
\newcommand\tti {\discretionary{t-} {t}{t\mbox {\hskip.03em}t}} 
\newcommand\tzi {\discretionary{t-} {z}{t\mbox {\hskip.03em}z}} 
\newcommand\myparagraphsymbol{\protect\unprotectedmyparagraphsymbol}


%% file: headernamesrest.tex
\newcommand\namefont{}

\newcommand\edcomment[1]{\mbox{}\raise.2ex\hbox{[}#1\raise.2ex\hbox{]}}
\usepackage{url}

\newcommand\majorfootroom{\raisebox{-1.9ex}{\rule{0ex}{.5ex}}}

\newcommand\footroom{\raisebox{-1.5ex}{\rule{0ex}{.5ex}}}

\newcommand\smallfootroom{\raisebox{-1.0ex}{\rule{0ex}{3ex}}}

\newcommand\majorheadroom{\rule{0ex}{3.2ex}}

\newcommand\mediumheadroom{\rule{0ex}{2.4ex}}

\newcommand\achilles {Achille\es} 

\newcommand\adolf    {Adolf}

\newcommand\andreas  {Andrea\es}

\newcommand\charles  {Charle\es}

\newcommand\claus    {Clau\es}

\newcommand\david    {David}

\newcommand\gerhard  {\mbox{Gerhard}}

\newcommand\heinrich {Hein\-rich}

\newcommand\jacques  {Jacque\es}

\newcommand\jean     {Jean}
\newcommand\john     {John}

\newcommand\julius   {Julius}

\newcommand\kurt     {Kurt}

\newcommand\ludwig   {Lud\-wig}

\newcommand\paolo    {Paolo}
\newcommand\paul     {Paul}
\newcommand\peter    {Peter}

\newcommand\rene     {Ren\'e}

\newcommand\wilfried {Wil\-fried}
\newcommand\wilhelm  {\mbox{Wilhelm}}

\newcommand\shorth          {\hskip.2em}

\newcommand\ackermannindex  {\index{Ackermann, Wilhelm (1896--1962)}}

\newcommand\ackermann       {{\namefont Acker\-mann}}

\newcommand\ackermannname     {\ackermannindex{\namefont \wilhelm\ \ackermann}}
\newcommand\ackermannlifetime{(1896--1962)}

\newcommand\archimedesindex {\index{Archimedes (28?--212(?)\,\BC)}}
\newcommand\archimedes      {{\namefont Archimedes}}

\newcommand\beckerplain     {Becker}
\newcommand\becker          {{\namefont\beckerplain}}

\newcommand\beatebeckerindex{\index{Becker, Beate (*1959)}}
\newcommand\beatebeckername {\beatebeckerindex{\namefont Beate \becker}}

\newcommand\bernaysindex    {\index{Bernays, Paul (1888--1977)}}

\newcommand\bernaysplain    {\mbox{Ber\-nay\es}}
\newcommand\bernays         {{\namefont\bernaysplain}}
\newcommand\bernaysshortnameplain
                            {\bernaysindex P\hspace*{-.15em}.~\bernaysplain}%
\newcommand\bernaysnameplain{\bernaysindex\paul\       \bernaysplain}

\newcommand\bernaysname     {{\namefont\bernaysnameplain}}

\newcommand\bernaysdeathyear{1977}
\newcommand\bernayslifetime {(1888--\bernaysdeathyear)}
\newcommand\bernaysludwigindex{\index{Bernays, Ludwig (1924--2020)}}

\newcommand\bernaysludwignameplain{\bernaysludwigindex\ludwig\ \bernaysplain}
\newcommand\bernaysludwigname{{\namefont\bernaysludwignameplain}}
\newcommand\bernaysreneindex                 {\index{Bernays, Ren{\'e} (*1949)}}

\newcommand\bernaysrenenameplain{\bernaysreneindex\rene\ \bernaysplain}
\newcommand\bernaysrenename {{\namefont\bernaysrenenameplain}}

\newcommand\bolzanoindex    {\index{Bolzano, Bernard (1781--1848)}}
\newcommand\bolzano         {{\namefont Bolzano}}

\newcommand\booleindex      {\index{Boole, George (1815--1864)}}
\newcommand\boole           {{\namefont Boole}}

\newcommand\myBoolean       {\boole an}

\newcommand\brouwersecondname{Brouwer} 

\newcommand\brouwerindex    {\index{Brouwer, L. E. J. (1881--1966)}}

\newcommand\brouwer         {{\namefont \brouwersecondname}}

\newcommand\cantorindex     {\index{Cantor, Georg (1845--1918)}}

\newcommand\cantorplain     {Cantor}
\newcommand\cantor          {{\namefont\cantorplain}}

\newcommand\cauchyindex     {\index{Cauchy, Augustin Louis (1789-1857)}}
\newcommand\cauchy          {{\namefont Cauchy}}

\newcommand\dillerindex                 {\index{Diller, Justus (*1936)}}

\newcommand\diller          {{\namefont Diller}}

\newcommand\dedekindindex   {\index{Dedekind, Richard (1831--1916)}}
\newcommand\dedekindplain   {Dedekind}
\newcommand\dedekind        {{\namefont\dedekindplain}}

\newcommand\engelerindex    {\index{Engeler, Erwin (*1930)}}

\newcommand\engeler         {{\namefont Engeler}}
\newcommand\engelername     {\engelerindex{\namefont Erwin \engeler}}

\newcommand\eudoxosindex
                   {\index{Eudoxos of Cnidus (\ca\,408\,BC\,--\,\ca\,355\,BC)}}
\newcommand\eudoxos         {{\namefont Eudoxos}}

\newcommand\jeanpaulindexindexsee{\index
  {Jean  Paul (pseudonym)|see{Richter, Johann Paul Friedrich}}\index
  {Paul, Jean (pseudonym)|see{Richter, Johann Paul Friedrich}}}
\newcommand\jeanpaulindex{\index{Richter, Johann Paul Friedrich (1763--1825)}}
\newcommand\jeanpaulname{\jeanpaulindex{\namefont \jean\ \paul}}

\newcommand\fermatbirthyear 
{160?}

\newcommand\fregeindex      {\index{Frege, Gottlob (1848--1925)}}

\newcommand\fregeplain      {Frege}
\newcommand\frege           {{\namefont\fregeplain}}

\newcommand\gabelsberger    {{\namefont Gabel\esi berger}}
\newcommand\englishgabelsbergerkurzschriftindex
                               {\index{Gabelsberger@Gabel\esi berger shorthand}}
\newcommand\englishgabelsbergerkurzschrift
                  {\englishgabelsbergerkurzschriftindex\gabelsberger\ shorthand}

\newcommand\gentzenindex                 {\index{Gentzen, Gerhard (1909--1945)}}

\newcommand\gentzen         {{\namefont Gentzen}}
\newcommand\gentzenname     {\gentzenindex{\namefont\gerhard\ \gentzen}}
\newcommand\gentzenlifetime {(1909--1945)}

\newcommand\glaschickindex  {\index{Glaschick, Rainer (*1949)}}
\newcommand\glaschickname   {\glaschickindex{\namefont Rainer Glaschick}}

\newcommand\goedelindex     {\index{Goedel, Kurt@G\"odel, Kurt (1906--1978)}}

\newcommand\goedel          {{\namefont G\oe del}}

\newcommand\goedelname      {\goedelindex{\namefont \kurt\   \goedel}}
\newcommand\goedellifetime  {(1906--1978)}
\newcommand\goedelsincompletenesstheorems{\goedel's incompleteness theorems}

\newcommand\secondincompletenesstheorem
                            {second \incompletenesstheorem}
\newcommand\secondIncompletenessTheorem
                            {Second \IncompletenessTheorem}

\newcommand\firstincompletenesstheorem
                            {first \incompletenesstheorem}
\newcommand\firstIncompletenessTheorem
                            {First \IncompletenessTheorem}

\newcommand\incompletenesstheorem{incompleteness theorem}
\newcommand\IncompletenessTheorem{Incompleteness Theorem}

\newcommand\grissindex                 {\index{Griss, George F. C. (1898--1953)}}

\newcommand\grissshortname
                     {\grissindex{\namefont G.\shorth F.\shorth C.\shorth Griss}}

\newcommand\hasenjaegerindex{\index{Hasenjaeger, Gisbert (1919--2006)}}

\newcommand\hasenjaegerplain{Hasen\-jaeger}
\newcommand\hasenjaegernameplain{Gisbert \hasenjaegerplain}
\newcommand\hasenjaeger     {{\namefont\hasenjaegerplain}}
\newcommand\hasenjaegershortnameplain{G.\shorth\hasenjaegerplain}
\newcommand\hasenjaegershortname
                       {\hasenjaegerindex{\namefont\hasenjaegershortnameplain}}
\newcommand\hasenjaegername
                       {\hasenjaegerindex{\namefont\hasenjaegernameplain}}
\newcommand\hasenjaegerbirthyear{1919}

\newcommand\hasenjaegernineteenfiftyctitle{Topologische Untersuchungen zur Semantik und Syntax eine\es\ erweiterten Pr\ae dikatenkalk\ue l\es}
\newcommand\hasenjaegernineteenfiftyctitleenglishparta
       {Topological Investigations on Semantics and Syntax of an Extended}
\newcommand\hasenjaegernineteenfiftyctitleenglishpartb{Predicate Calculus}
\newcommand\hasenjaegernineteenfiftyctitleenglish
                                 {\hasenjaegernineteenfiftyctitleenglishparta\
                                  \hasenjaegernineteenfiftyctitleenglishpartb}
\newcommand\hasenjaegernineteenfiftyctitleenglishwithlinebreak
                                 {\hasenjaegernineteenfiftyctitleenglishparta\\
        \mbox{~~~~~~~~~~~~~~~~~~~}\hasenjaegernineteenfiftyctitleenglishpartb}

\newcommand\hasenjaegerhabiltitleenglish
           {Consistent Axiom Systems Without Standard Model}
\newcommand\hasenjaegerfatherindex{\index{Hasenjaeger, Edwin (1888--1972)}}
\newcommand\hasenjaegerfathername
                        {\hasenjaegerfatherindex{\namefont Edwin \hasenjaeger}}
\newcommand\hasenjaegerfatherlifetime{(1888--1972)}

\newcommand\herbrandindex     {\index{Herbrand, Jacques (1908--1931)}}

\newcommand\herbrand        {{\namefont Herbrand}}

\newcommand\herbrandname     {\herbrandindex{\namefont \jacques\ \herbrand}}
\newcommand\herbranddeathyear{1931}
\newcommand\herbrandlifetime{(1908--\herbranddeathyear)}

\newcommand\herbrandsfundamentaltheoremindex
                                       {\index{Herbrand's Fundamental Theorem}}

\newcommand\herbrandsfundamentaltheoremnoindex
                                             {\herbrand'\es\ \fundamentaltheorem}
\newcommand\herbrandsfundamentaltheorem
           {\herbrandsfundamentaltheoremindex\herbrandsfundamentaltheoremnoindex}

\newcommand\fundamentaltheorem{Fundamental Theorem}

\newcommand\heytingindex                 {\index{Heyting, Arend (1898--1980)}}

\newcommand\heyting         {\mbox{\namefont Heyting}}

\newcommand\hilbertindex                 {\index{Hilbert, David (1862--1943)}}

\newcommand\hilbertplain    {Hilbert}
\newcommand\hilbert         {\mbox{\namefont\hilbertplain}}

\newcommand\hilbertnameplain{\david\ \hilbertplain}
\newcommand\hilbertname     {\hilbertindex{\namefont\hilbertnameplain}}
\newcommand\hilbertlifetime {(1862--1943)}
\newcommand\hilbertsprogram {\hilbert'\es\ \programme}

\newcommand\hilbertsepsilonlongindex{\index{Hilbert!'s epsilon}}

\newcommand\hilbertsepsilon
                     {\hilbertsepsilonlongindex\hilbert'\es\ \nlbmath\varepsilon}

\newcommand\hilbertbernaysplain{\hilbertplain--\bernaysplain}
\newcommand\hilbertbernays  {{\namefont\hilbertbernaysplain}}

\newcommand\hilbertbernaysproject{{\namefont\hilbertbernays\ Pro\-ject}}

\newcommand\titlehilbertnineteenhundredandseventeen{Axiomatisches Denken}
\newcommand\titlehilbertnineteenhundredandseventeenenglish{Axiomatic Tought}

\newcommand\hitlername      {\adolf\ Hitler}

\newcommand\hungerbuehlerindex{\index{Hungerb\"uhler, Norbert (*1964)}}
\newcommand\hungerbuehler   {{\namefont Hunger\-b\ue hler}}
\newcommand\hungerbuehlername
                         {\hungerbuehlerindex{\namefont Norbert \hungerbuehler}}

\newcommand\koenig          {{\namefont K\oe nig}}

\newcommand\koenigfatherindex{\index{Koenig@K{\"o}nig, Julius (1849--1913)}}

\newcommand\koenigfathername{\koenigfatherindex{\namefont\julius\ \koenig}}
\newcommand\koenigfathertitle
    {Neue Grundlagen der Logik, Arithmetik, und Mengen\-lehre}
\newcommand\koenigfathertitleenglish
        {New Foundations of Logic, Arithmetic, and Set Theory}

\newcommand\kroneckerindex  {\index{Kronecker, Leopold (1823--1891)}}
\newcommand\kronecker       {\mbox{\namefont Kronecker}}

\newcommand\loewenheim      {{\namefont L\oe wen\-heim}}

\newcommand\loewenheimskolem{\loewenheim--\skolem}
\newcommand\loewenheimskolemtheorem{\index
                    {Loewenheim-Skolem Theorem@L{\oe}wenheim--Skolem Theorem}%
                                                     \loewenheimskolem\ Theorem}

\newcommand\mancosuindex    {\index{Mancosu, Paolo (*1960)}}

\newcommand\mancosu         {{\namefont Mancosu}}
\newcommand\mancosuname     {\mancosuindex{\namefont\paolo\ \mancosu}}

\newcommand\merayindex      {\index{M{\'e}ray, Charles (1835--1911)}}
\newcommand\meray           {{\namefont M\'eray}}
\newcommand\merayshortname  {{\namefont Ch.\shorth \meray}}
\newcommand\merayname       {\merayindex{\namefont\charles\ \meray}}
\newcommand\meraylifetime   {(1835--1911)}

\newcommand\neumannindex    {\index{Neumann, John von (1903--1957)}}

\newcommand\neumann         {{\namefont Neumann}}

\newcommand\neumannname     {\neumannindex{\namefont\john\ von \neumann}}
\newcommand\neumannlifetime {(1903--1957)}

\newcommand\peanoindex      {\index{Peano, Guiseppe (1858--1932)}}

\newcommand\peanoplain      {Peano}
\newcommand\peano           {{\namefont\peanoplain}}

\newcommand\peirceindex     {\index{Peirce, Charles S. (1839--1914)}}

\newcommand\peirce          {{\namefont Peirce}}

\newcommand\russellindex                 {\index{Russell, Bertrand (1872--1970)}}

\newcommand\russellsparadoxindex{\index{Russell's Paradox}}
\newcommand\russellplain    {Russell}
\newcommand\russell         {{\namefont\russellplain}}

\newcommand\russellsparadox {\russellsparadoxindex\russell'\es\ Paradox}

\newcommand\schmehindex                 {\index{Schmeh, Klaus (*1970)}}

\newcommand\schmeh          {{\namefont Schmeh}}
\newcommand\schmehname      {\schmehindex{\namefont Klaus \schmeh}}
\newcommand\schmidtplain    {Schmidt}
\newcommand\schmidt         {{\namefont\schmidtplain}}
\newcommand\arnoldschmidtindex{\index{Schmidt, Hermann Arnold (1902--1967)}}
\newcommand\arnoldschmidtname
                    {\arnoldschmidtindex{\namefont H.\shorth Ar\-nold \schmidt}}
\newcommand\fkschmidtindex  {\index{Schmidt, Friedrich Karl (1901--1977)}}
\newcommand\fkschmidtshortnameplain
                            {\fkschmidtindex F.\shorth K.\shorth\schmidtplain}
\newcommand\fkschmidtshortname{{\namefont\fkschmidtshortnameplain}}

\newcommand\scholzindex                 {\index{Scholz, Heinrich (1884--1956)}}

\newcommand\scholz          {{\namefont Scholz}}

\newcommand\scholznameplain {\scholzindex\heinrich\ \scholz}
\newcommand\scholzname      {{\namefont\scholznameplain}}
\newcommand\scholzlifetime  {(1884--1956)}

\newcommand\schroederindex  {\index{Schroeder@Schr\"oder, Ernst (1841--1902)}}

\newcommand\schroeder       {{\namefont Schr\oe der}}

\newcommand\schuetteindex   {\index{Schuette@Sch\"utte, Kurt (1909--1998)}}

\newcommand\schuette        {{\namefont Sch\ue tte}}

\newcommand\schuettename     {\schuetteindex{\namefont\kurt\    \schuette}}

\newcommand\siegindex       {\index{Sieg, Wilfried (*1945)}}

\newcommand\sieg            {{\namefont Sieg}}
\newcommand\siegnamenoindex {{\namefont\wilfried\ \sieg}}
\newcommand\siegname        {\siegindex\siegnamenoindex}

\newcommand\skolem          {{\namefont Skolem}}

\newcommand\stolzenburgindex{\index{Stolzenburg, Frieder (*1966)}}
\newcommand\stolzenburg     {{\namefont Stol\-zen\-burg}}
\newcommand\stolzenburgname{\stolzenburgindex{\namefont Frie\-der \stolzenburg}}

\newcommand\weierstrassindex{\index{Weierstrass@Weierstra\sz, Karl (1815--1897)}}
\newcommand\weierstrass     {{\namefont Weierstra\sz}}

\newcommand\weylindex       {\index{Weyl, Hermann (1885--1955)}}

\newcommand\weyl            {{\namefont Weyl}}

\newcommand\whiteheadindex  {\index{Whitehead, Alfred North (1861--1947)}}
\newcommand\whitehead       {{\namefont White\-head}}

\newcommand\PM              {Prin\-ci\-pia Mathematica}

\newcommand\wirthindex                 {\index{Wirth, Claus-Peter (*1963)}}

\newcommand\wirthplain      {Wirth}

\newcommand\wirthnameplainnoindex{\claus-\peter\ \wirthplain}
\newcommand\wirthnamenoindex{{\namefont\wirthnameplainnoindex}}

\newcommand\wirthname       {\wirthindex\wirthnamenoindex}

\newcommand\wittgensteinindex{\index{Wittgenstein, Ludwig (1889--1951)}}

\newcommand\wittgenstein    {{\namefont Wittgen\-stein}}

\newcommand\achillesparadoxindex{\index{Achilles Paradox}}
\newcommand\achillesparadox{\achillesparadoxindex{{\namefont\achilles} Paradox}}

\newcommand\Apr  {April}
\newcommand\atall{at \nolinebreak all}
\newcommand\Aug  {Aug.}

\newcommand\ca   {ca.}

\newcommand\f    {\mbox{}{f.}}   
\newcommand\Feb  {Feb.}

\newcommand\Drrernat{Dr.\,\,rer.\,\,nat.}

\newcommand\EMAIL{E-mail}

\newcommand\Jan  {Jan.}

\newcommand\nee{n\'ee}

\newcommand\onlyif{only \nolinebreak if}

\newcommand\Prof {Prof.}
\newcommand\programme{program}

\newcommand\qedhelp[1]{Q.e.d.~({#1})}
\newcommand\getittotheright[1]  
{\hfill\mbox{}\penalty 100\mbox{\ \,}\nolinebreak
\hfill\hfill\hfill\nolinebreak\mbox{#1}\ignorespaces}

\newcommand\Qedbf    [1]{\mbox{\bf\qedhelp{#1}}}

\newcommand\QEDbf    [1]{\getittotheright{\Qedbf    {#1}}}

\newcommand\role{r\^ole}
\newcommand\Sep  {Sept.}

\newcommand\Univ {Univ.}
\newcommand\Vol  {Vol.}
\newcommand\WWW  {WWW}

\newcommand\aswell{as \nolinebreak well}
\newcommand\aswellas{\aswell\ \nolinebreak as}
\newcommand\BC   {{\sc b.c.}} 

\newcommand\Cf   {Cf.}
\newcommand\cf   {cf.}

\newcommand\Cfnlb{\Cf\nolinebreak}
\newcommand\cfnlb{\cf\nolinebreak}

\newcommand\CS   {Computer \Sci}

\newcommand\Dec  {Dec.}
\newcommand\Dept {Dept.}

\newcommand\ednnodot{edn}
\newcommand\edn  {\ednnodot.}

\newcommand\eg   {e.g.}

\newcommand\etc  {\&c.}

\newcommand\extd {extd.}

\newcommand\grundlagendergeometrieindex{\index{Grundlagen!der Geometrie}}

\newcommand\englishgrundlagendergeometrie
                    {\grundlagendergeometrieindex Foundations of Geometry}
\newcommand\grundlagendermathematikindex{\index{Grundlagen!der Mathematik}}
\newcommand\grundlagendermathematiknoindex{Grund\-lagen der Ma\-the\-ma\-tik}

\newcommand\grundlagendermathematik
             {\grundlagendermathematikindex\grundlagendermathematiknoindex}

\newcommand\ie   {i.e.}

\newcommand\incl {incl.}
\newcommand\infact{in \nolinebreak fact}
\newcommand\Infact{In \nolinebreak fact}

\newcommand\Jun  {June}
\newcommand\Jul  {July}

\newcommand\May  {May}
\newcommand\Mar  {March}

\newcommand\MrUS {Mr.}

\def\note{Note}
\newcommand\notes{Notes}
\newcommand\Oct  {Oct.}
\newcommand\p    {p.}
\newcommand\ppnodot{pp}
\newcommand\pp   {\ppnodot.}
\newcommand\PP[2]{\pp\,\ignorespaces#1--\ignorespaces#2}

\newcommand\PhD  {PhD}
\newcommand\PhDthesis{\PhD\ thesis}

\newcommand\Proc {Proc.}

\newcommand\rev  {rev.}

\newcommand\sect {\myparagraphsymbol} 
\newcommand\sects{\myparagraphsymbol\myparagraphsymbol}

\newcommand\Sci  {Sci.}

\newcommand\wrt  {w.r.t.}

\newcommand\dasheisst{d.\nolinebreak\,h.}

\newcommand\evtl{evtl.}

\newcommand\Vgl{Vgl.}

\newcommand\zB{z.\penalty100\,B.}

\newcommand\thewordand{and}

\newcommand\litspagerefs[2]{Pages #1 \nolinebreak \thewordand\ \nolinebreak #2}

\newcommand\litnoteref[1]{\note\,#1}

\newcommand\itemname{item}

\newcommand\lititemref[1]{\itemname\,#1}

\newcommand\litsectref[1]{\sect\,#1} 
\newcommand\nlblitsectref[1]{\nolinebreak\litsectref{#1}}

\newcommand\litnoterefs[2]
{\notes\         \nolinebreak #1 \thewordand\ \nolinebreak #2}
\newcommand\litsectrefs[2]
{\sects\         \nolinebreak #1 \thewordand\ \nolinebreak #2}

\newcommand\litchaprefs[2]
{Chapters        \nolinebreak #1 \thewordand\ \nolinebreak #2}

\newcommand\litsectrefss
[3]{\sects\ \nolinebreak #1, #2, \thewordand\ \nolinebreak #3}

\newcommand\spagerefs[2]{\litspagerefs{\pageref{#1}}{\pageref{#2}}}

\newcommand\noteref[1]{\litnoteref{\ref{#1}}}
\newcommand\sectref[1]{\litsectref{\ref{#1}}}
\newcommand\nlbsectref[1]{\nolinebreak\sectref{#1}}

\newcommand\noterefs[2]{\litnoterefs{\ref{#1}}{\ref{#2}}}
\newcommand\sectrefs[2]{\litsectrefs{\ref{#1}}{\ref{#2}}}

\newcommand\sectrefss[3]{\litsectrefss{\ref{#1}}{\ref{#2}}{\ref{#3}}}

\newcommand\nthpositioner[2]
{#1\raisebox{0.52ex}{\scriptsize\hspace{0.07em}#2}}
\newcommand\nth[1]{\nthtinypositioner{#1}{\nthstring{#1}}}
\newcommand\nthtinypositioner[2]{#1\raisebox{0.52ex}{\tiny\hspace{0.07em}#2}}

\newcommand\mthpositioner[2]
{\math{#1}\raisebox{0.52ex}{\scriptsize\hspace{0.07em}#2}}
\newcommand\modulointocountzero[2]
{\count1=#1
\count2=#2
\count0=\count1
\divide  \count0 by \count2
\multiply\count0 by-\count2
\advance \count0 by \count1}
\newcommand\absolutevalueintocountzero[1]
{\count0=#1
\ifnum\count0<0\multiply\count0 by -1\fi}
\newcommand\nthstring[1]
{\def\myargone{#1}\ifcat a\myargone th\else\nthstringnochar{#1}\fi}
\newcommand\nthstringnochar[1]
{\absolutevalueintocountzero{#1}%
\modulointocountzero{\count0}{100} 
\ifnum\count0>9\ifnum\count0<20 th\else\stupidnthstring\fi
                                  \else\stupidnthstring\fi}
\newcommand\stupidnthstring
{\modulointocountzero{\count0}{10}
\ifnum\count0=1 \hskip-0.2em st\else
\ifnum\count0=2 nd\else
\ifnum\count0=3 rd\else 
                th\fi\fi\fi}

\newcommand\writeascents
[1]{\count4=#1
\ifnum\count4<0 
-\multiply\count4 by -1\fi
\modulointocountzero{\count4}{10}
\divide\count4 by 10
\count3=\the\count0
\modulointocountzero{\count4}{10}
\divide\count4 by 10
\the\count4
.\the\count0
\the\count3
}

\newcommand\frenchnthstring[1]
{\def\myargone{#1}\ifcat a\myargone th\else\frenchnthstringnochar{#1}\fi}
\newcommand\frenchnthstringnochar[1]
{\absolutevalueintocountzero{#1}%
\modulointocountzero{\count0}{100} 
\ifnum\count0>9\ifnum\count0<20 th\else\frenchstupidnthstring\fi
                                  \else\frenchstupidnthstring\fi}
\newcommand\frenchstupidnthstring
{\modulointocountzero{\count0}{10}
\ifnum\count0=1 \hskip-0.2em re\else
\ifnum\count0=2 me\else
\ifnum\count0=3 rd\else 
                th\fi\fi\fi}

\newcommand\CLAM      {{\rm CL\kern-.36em\raise.39ex\hbox{\sc a}\kern-.15emM}}

\newcommand\TEXMACS   {{\sc T\kern-.1667em\lower.5ex\hbox{E}\kern-.125emX\kern-.1em\lower.5ex\hbox{\textsc{m\kern-.05ema\kern-.125emc\kern-.05ems}}}}

\newcommand\CMU{Carnegie Mellon \Univ}

\newcommand\UIUC{\Univ\ of Illinois at Urbana--Champaign}

\newcommand\Amsterdam      {Amster\-dam}

\newcommand\Berlin         {Ber\-lin\index{Berlin}}

\newcommand\FreiburgnoBreisgau{Freiburg}

\newcommand\Freiburg       {\FreiburgnoBreisgau\ (Brei\esi gau, Germany)}

\newcommand\Goettingen     {\index{Goettingen@{G\oe ttingen}}G\oe ttingen}

\newcommand\Heidelberg     {\index{Heidelberg}Heidel\-berg}

\newcommand\HildesheimnoGermany{Hilde\esi heim}
\newcommand\Hildesheim     {\HildesheimnoGermany\ (Germany)}

\newcommand\Leipzig        {Leip\-zig}

\newcommand\MuelheimnoGermany{M\ue lheim an der Ruhr}
\newcommand\Muelheim       {\MuelheimnoGermany\ (Germany)}
\newcommand\Muenchen       {M\ue n\-chen}
\newcommand\Muensterindex  {\index{Muenster@{M\ue nster (Westfalen, Germany)}}}
\newcommand\MuensternoWestfalen{\Muensterindex M\ue nster}
\newcommand\MuensternoGermany{\Muensterindex\MuensternoWestfalen\ (Westfalen)}
\newcommand\Muenster{\Muensterindex\MuensternoWestfalen\ (Westfalen, Germany)}

\newcommand\Princetonnostate{Princeton}
\newcommand\Princeton      {\Princetonnostate\ (NJ)}

\newcommand\SB             {Saar\-br\ue cken}

\newcommand\Zuerichindex   {\index{Zuerich@{Z\ue rich\,(Switzerland)}}}
\newcommand\Zuerich{\Zuerichindex\mbox{Z\ue rich}}

\newcommand\plzuniSB  {\mbox{66123}}

\newcommand\englishETHZlongwithoutshort
               {Swiss Federal Institute of Technology in \Zuerich}

\newcommand\ETHshort{ETH}
\newcommand\ETHZshort{\ETHshort\ Zurich}
\newcommand\ETHZofficialarchivereference{ETH-Bibliothek, Hochschularchiv}

\newcommand\uniSBenglishshort{Saar\-land \Univ}

\def       \emailcp      {{\tt wirth@logic.at}}

\newcommand\Institutedept
{\Dept\ of \CS}
\newcommand\Instituteinst
{\mbox\uniSBenglishshort}
\newcommand\Instituteplac
{\plzuniSB\,\SB}
\newcommand\Institutecoun
{Germany}
\newcommand\Institutestre
{R\ae mistr.\,101}
\newcommand\Institute
{\Institutedept, \Instituteinst, 
      \Instituteplac, \Institutecoun}

\newcommand\IAS{Institute for Advanced Study, \Princeton}

\newcommand\academicpress{Academic Press (\ELSEVIER), San Diego (CA)}

\newcommand\ELSEVIER{Elsevier}
\newcommand\elsevier{Elsevier, \Amsterdam}

\newcommand\newspaperreference[5]
{\def\nameofjournalpress{#2}#1, #4 #5, #3\if?\nameofjournalpress
\else, #2\fi}

\newcommand\dateinjournal[1]{}

\newcommand\journalreference[6]
{\def\nameofjournalpress{#2}#1\nolinebreak\hskip.2em%
\dateinjournal{(#3) }{\mbox{\bf #4}}, \PP{#5}{#6}\if?\nameofjournalpress
\else, #2\fi}

\newcommand\journalreferenceprintyear[6]
{\def\nameofjournalpress{#2}#1 
#4:#5--#6, #3%
\if?\nameofjournalpress
\else, #2\fi}

\newcommand\journalreferenceprintyearaspartofnumber[6]
{\def\nameofjournalpress{#2}#1 
{#4/#3}, \PP{#5}{#6}\if?\nameofjournalpress
\else, #2\fi}

\newcommand\jscname
{J. Symbolic Computation}

\newcommand\jscprintyear
{\journalreferenceprintyear{\jscname}\academicpress}

\newcommand\tcsname{Theoretical \CS}
\newcommand\tcsjournal
{\journalreference\tcsname\elsevier}
\newcommand\tcsjournalprintyear
{\journalreferenceprintyear\tcsname\elsevier}

\newcommand\archivdeutschesmuseum{Deutsche\es\ Museum, \Muenchen, Archiv}
\newcommand\hasenjaegerlegacy{Legacy of \hasenjaegername, \archivdeutschesmuseum: ``Nachl\ae sse~H, NL\,288''}

\urldef\wirthkuhnurl\url
{http://wirth.bplaced.net/SEKI/welcome.html#SWP-2007-01}

\urldef\urlsrninetythreedashzerofive\url
{http://wirth.bplaced.net/SEKI/welcome.html#SR-93-05}

\urldef\urlsreightyeighttwelve\url
{http://wirth.bplaced.net/SEKI/welcome.html#SR-88-12}

%% file: headerforformulas.tex
\mathcommand\ident[1]{\mathsf{#1}}
\newcommand\plussymbol  {\ident{+}}
\newcommand\minussymbol {\ident{-}}
\newcommand\dividesymbol{\ident{/}}
\newcommand\timessymbol {\ident{*}}

\newcommand\set     {\ident{set}}

\newcommand\naturalssymbol{\ident{naturals}}
\newcommand\gensymsymbol{\ident{gensym}}
\mathcommand\mbpsymbol{\ident{m\hspace{-0.055em}b\hspace{-0.045em}p}}

\newcommand\csymbol     {\ident c}
\newcommand\esymbol     {\ident e}
\newcommand\fsymbol     {\ident f}
\newcommand\gsymbol     {\ident g}
\newcommand\hsymbol     {\ident h}
\newcommand\ksymbol     {\ident k}
\newcommand\psymbol     {\ident p}
\newcommand\ssymbol     {\ident s}
\newcommand\Everysymbol {\ident{Every}}
\newcommand\Permsymbol {\ident{Perm}}
\newcommand\RExistssymbol{\ident{Rexists}}
\newcommand\invertsymbol{\ident{invert}}
\newcommand\invsymbol{\ident{inv}}
\newcommand\abssymbol   {\ident{abs}}
\newcommand\cnssymbol   {\ident{cons}}
\mathcommand\cnsindexsymbol[1]{\ident{cons}_{#1}}
\newcommand\carsymbol   {\ident{car}}

\newcommand\cdrsymbol   {\ident{cdr}}
\newcommand\lengthsymbol{\ident{length}}
\newcommand\sizesymbol{\ident{size}}
\newcommand\dlsymbol    {\ident{dl}}
\newcommand\dloncesymbol{\ident{delfirst}}
\newcommand\rcsymbol    {\ident{rc}}
\newcommand\brsymbol    {\ident{br}}
\newcommand\revtailsymbol{\ident{revtail}}
\newcommand\revsymbol{\ident{rev}}
\newcommand\appendsymbol {\ident{append}}
\newcommand\zeropredicatesymbol{\ident{zerop}}
\newcommand\eqsymbol        {\ident{eq}}
\newcommand\ifthensymbol    {\mbox{\ident{If{}Then}}}
\newcommand\ifthenelsesymbol{\mbox{\ident{If{}ThenElse}}}
\mathcommand\eqindexsymbol        [1]{\eqsymbol        _{#1}}
\mathcommand\ifthenindexsymbol    [1]{\ifthensymbol    _{#1}}
\mathcommand\ifthenelseindexsymbol[1]{\ifthenelsesymbol_{#1}}
\newcommand\orsymbol    {\ident{or}}
\newcommand\andsymbol   {\ident{and}}
\newcommand\leqsymbol   {\ident{leq}}
\newcommand\lessymbol   {\ident{less}}
\newcommand\lexlessymbol{\ident{lexless}}
\newcommand\lexlimlessymbol{\ident{lexlimless}}
\newcommand\lexsymbol   {\ident{lex}}
\newcommand\acksymbol   {\ident{ack}}
\newcommand\switchsymbol{\ident{switch}}
\newcommand\swatchsymbol{\ident{swatch}}
\newcommand\diveinssymbol{\ident{div1}}
\newcommand\divzweisymbol{\ident{div2}}
\newcommand\divrestsymbol{\ident{divrest}}
\newcommand\diveinstailsymbol{\ident{div1tail}}
\newcommand\divzweitailsymbol{\ident{div2tail}}
\newcommand\remsymbol{\ident{rem}}
\newcommand\divsymbol{\ident{div}}

\newcommand\turingmachinesymbol{\ident T}
\newcommand\terminatespsymbol  {\ident{terminatesp}}
\newcommand\statesymbol        {\ident{state}}
\newcommand\cmdsymbol          {\ident{cmd}}
\newcommand\nthsymbol          {\ident{nth}}
\newcommand\doublesymbol       {\ident{double}}

\newcommand\ppsymbol           {\ident{p}}
\newcommand\qpsymbol           {\ident{q}}
\newcommand\Epsymbol           {\ident{E}}
\newcommand\Ppsymbol           {\ident{P}}
\newcommand\Qpsymbol           {\ident{Q}}
\newcommand\Marriessymbol      {\ident{Marries}}
\newcommand\Lovessymbol        {\ident{Loves}}
\newcommand\StolenBysymbol     {\ident{StolenBy}}
\newcommand\Humansymbol        {\ident{Human}}
\newcommand\Evensymbol         {\ident{Even}}
\newcommand\Oddsymbol          {\ident{Odd}}
\newcommand\Primesymbol        {\ident{Prime}}
\newcommand\EveryPairsymbol   {\ident{EveryPair}}
\newcommand\Givesymbol         {\ident{Give}}
\newcommand\Fathersymbol       {\ident{Father}}
\newcommand\Elephantpsymbol    {\ident{Elephant}}
\newcommand\Flowerpsymbol    {\ident{Flower}}
\newcommand\Germanpsymbol      {\ident{German}}
\newcommand\Bicyclepsymbol     {\ident{Bicycle}}
\newcommand\Hugepsymbol        {\ident{Huge}}
\newcommand\Animalpsymbol      {\ident{Animal}}
\newcommand\Malepsymbol        {\ident{Male}}
\newcommand\Boypsymbol         {\ident{Boy}}
\newcommand\Girlpsymbol        {\ident{Girl}}
\newcommand\Femalepsymbol      {\ident{Female}}
\newcommand\Roundpsymbol       {\ident{Round}}
\newcommand\Quadrangularpsymbol{\ident{Quadrangular}}
\newcommand\Metpsymbol         {\ident{Met}}
\newcommand\Kissedpsymbol      {\ident{Kissed}}
\newcommand\Bishopsymbol       {\ident{Bishop}}
\newcommand\mindexsymbol[1]{\existsvari w{#1}}

\newcommand\nonnegpsymbol      {\ident{nonnegp}}
\newcommand\wellsymbol         {\ident{well}}
\newcommand\welltailsymbol     {\ident{welltail}}
\newcommand\varsymbol          {\ident{var}}
\newcommand\aritysymbol        {\ident{arity}}

\newcommand\whilesymbol        {\ident{while}}

\newcommand\nullsymbol         {\ident{null}}
\newcommand\hdsymbol           {\ident{hd}}
\newcommand\tlsymbol           {\ident{tl}}
\newcommand\insymbol           {\ident{in}}
\newcommand\applysymbol        {\ident{app}}
\newcommand\termsymbol         {\ident{term}}
\newcommand\russellsymbol      {\ident{russell}}
\newcommand\sqrtindordsymbol[1]{\ident{sqrtio#1}}
\mathcommand\tightim{\longrightarrow}
\mathcommand\im{\ \tightim\ }
\mathcommand\rs{\:\rulesugar\:\:}
\mathcommand\rulesugar{\longleftarrow}

\mathcommand\doublepp[1]      {\doublesymbol   \beginargs{#1}\allargs}
\mathcommand\aritypp[1]      {\aritysymbol   \beginargs{#1}\allargs}
\mathcommand\lengthpp[1]      {\lengthsymbol   \beginargs{#1}\allargs}
\mathcommand\sizepp[1]      {\sizesymbol   \beginargs{#1}\allargs}
\mathcommand\wellpp[1]      {\wellsymbol   \beginargs{#1}\allargs}
\mathcommand\welltailpp[1]      {\welltailsymbol   \beginargs{#1}\allargs}
\mathcommand\varpp[1]      {\varsymbol   \beginargs{#1}\allargs}
\mathcommand\rempp[2]    {\remsymbol\beginargs{#1}\separgs{#2}\allargs}
\mathcommand\divpp[2]    {\divsymbol\beginargs{#1}\separgs{#2}\allargs}
\mathcommand\divrestpp[2]    {\divrestsymbol\beginargs{#1}\separgs{#2}\allargs}
\mathcommand\diveinspp[2]    {\diveinssymbol\beginargs{#1}\separgs{#2}\allargs}
\mathcommand\divzweipp[3]    {\divzweisymbol\beginargs{#1}\separgs{#2}
\separgs{#3}\allargs}
\mathcommand\diveinstailpp[4]    {\diveinstailsymbol\beginargs{#1}\separgs{#2}
\separgs{#3}\separgs{#4}\allargs}
\mathcommand\divzweitailpp[6]    {\divzweitailsymbol\beginargs{#1}\separgs{#2}
\separgs{#3}\separgs{#4}\separgs{#5}\separgs{#6}\allargs}
\mathcommand\mbppp[2]         {\mbpsymbol   \beginargs{#1}\separgs{#2}\allargs}
\mathcommand\revpp[1]     
{\revsymbol\beginargs{#1}\allargs}
\mathcommand\revppi[2]     
{\revsymbol^{#1}\beginargs{#2}\allargs}
\mathcommand\revtailpp[2]     
{\revtailsymbol\beginargs{#1}\separgs{#2}\allargs}
\mathcommand\revtailppi[3]
{\revtailsymbol^{#1}\beginargs{#2}\separgs{#3}\allargs}
\mathcommand\Permpp[2]     
{\Permsymbol\beginargs{#1}\separgs{#2}\allargs}
\mathcommand\Permppi[3]
{\Permsymbol^{#1}\beginargs{#2}\separgs{#3}\allargs}
\mathcommand\appendpp[2]      
{\appendsymbol \beginargs{#1}\separgs{#2}\allargs}
\mathcommand\appendppi[3]      
{\appendsymbol^{#1}\beginargs{#2}\separgs{#3}\allargs}
\mathcommand\Everypp[2]      
{\Everysymbol \beginargs{#1}\separgs{#2}\allargs}
\mathcommand\RExistspp[1]      
{\RExistssymbol \beginargs{#1}\allargs}
\mathcommand\appendlongpp[2]      
{\appendsymbol\left(\begin{array}{@{}l@{}}{#1}\separgs\\{#2}\end{array}\right)}
\mathcommand\cnspp[2]         {\cnssymbol   \beginargs{#1}\separgs{#2}\allargs}
\mathcommand\cnsppi[3]       {\cnssymbol^{#1}\beginargs{#2}\separgs{#3}\allargs}
\mathcommand\cnsindexpp[3]
{\cnsindexsymbol{#1}\beginargs{#2}\separgs{#3}\allargs}
\mathcommand\dlpp[2]          {\dlsymbol    \beginargs{#1}\separgs{#2}\allargs}
\mathcommand\dloncepp[2]      {\dloncesymbol\beginargs{#1}\separgs{#2}\allargs}
\mathcommand\dlonceppi[3]{\dloncesymbol^{#1}\beginargs{#2}\separgs{#3}\allargs}
\mathcommand\rcpp[2]          {\rcsymbol    \beginargs{#1}\separgs{#2}\allargs}
\mathcommand\brpp[2]          {\brsymbol    \beginargs{#1}\separgs{#2}\allargs}
\mathcommand\orpp[2]          {\orsymbol    \beginargs{#1}\separgs{#2}\allargs}
\mathcommand\andpp[2]         {\andsymbol   \beginargs{#1}\separgs{#2}\allargs}
\mathcommand\shortcnspp[2]    {\csymbol     \beginargs{#1}\separgs{#2}\allargs}
\mathcommand\tightshortcnspp[2]
{\csymbol\beginargs{#1}\tightsepargs{#2}\allargs}
\mathcommand\spp[1]           {\ssymbol     \beginargs{#1}\allargs}
\mathcommand\sppiterated[2]   {\ssymbol^{#1}\beginargs{#2}\allargs}
\mathcommand\sqrtindordpp[3]
                       {\sqrtindordsymbol{#1}\beginargs{#2}\separgs{#3}\allargs}
\mathcommand\ppp[1]           {\psymbol     \beginargs{#1}\allargs}
\mathcommand\pppiterated[2]   {\psymbol^{#1}\beginargs{#2}\allargs}
\mathcommand\zeropp           {\ident 0}
\mathcommand\Julietpp         {\ident{Juliet}}
\mathcommand\Romeopp          {\ident{Romeo}}
\mathcommand\Ipp              {\ident I}
\mathcommand\onepp            {\ident1}
\mathcommand\twopp            {\ident2}
\mathcommand\threepp          {\ident3}
\mathcommand\invertpp[1]      {\invertsymbol\beginargs{#1}\allargs}
\mathcommand\invpp[1]         {\invsymbol\beginargs{#1}\allargs}
\mathcommand\abspp[1]         {\abssymbol\beginargs{#1}\allargs}
\mathcommand\naturalspp[1]    {\naturalssymbol\beginargs{#1}\allargs}
\mathcommand\gensympp[1]      {\gensymsymbol\beginargs{#1}\allargs}
\mathcommand\nilpp            {\ident{nil}}
\mathcommand\falsepp          {\ident{false}}
\mathcommand\truepp           {\ident{true}}
\mathcommand\FALSEpp          {\ident{FALSE}}
\mathcommand\TRUEpp           {\ident{TRUE}}
\mathcommand\UNDEFpp          {\ident{UNDEF}}
\mathcommand\weirdppp         {\ident{weirdp}}
\mathcommand\ambigppp         {\ident{ambigp}}
\mathcommand\zeropredicatepp[1]{\zeropredicatesymbol\beginargs{#1}\allargs}
\mathcommand\cppeins       [1]{\csymbol     \beginargs{#1}\allargs}
\mathcommand\cppzwei       [2]{\csymbol\beginargs{#1}\separgs{#2}\allargs}
\mathcommand\eppeins       [1]{\esymbol     \beginargs{#1}\allargs}
\mathcommand\fppeins       [1]{\fsymbol     \beginargs{#1}\allargs}
\mathcommand\fppeinsindex  [2]{\fsymbol_{#1}\beginargs{#2}\allargs}
\mathcommand\fppeinsiterated[2]{\fsymbol^{#1}\beginargs{#2}\allargs}
\mathcommand\gppeins       [1]{\gsymbol     \beginargs{#1}\allargs}
\mathcommand\gppzwei       [2]{\gsymbol     \beginargs{#1}\separgs{#2}\allargs}
\mathcommand\hppeins       [1]{\hsymbol     \beginargs{#1}\allargs}
\mathcommand\kppeins       [1]{\ksymbol     \beginargs{#1}\allargs}
\mathcommand\appzero          {\ident a}
\mathcommand\bppzero          {\ident b}
\mathcommand\cppzero          {\ident c}
\mathcommand\dppzero          {\ident d}
\mathcommand\eppzero          {\ident e}
\mathcommand\eqindexpp[3]{\eqindexsymbol{#1}\beginargs{#2}\separgs{#3}\allargs}
\mathcommand\ifthenindexpp
[3]{\ifthenindexsymbol{#1}\beginargs{#2}\separgs{#3}\allargs}
\mathcommand\ifthenelseindexpp
[4]{\ifthenelseindexsymbol{#1}\beginargs{#2}\separgs{#3}\separgs{#4}\allargs}
\mathcommand\eqpp[2]{\eqsymbol\beginargs{#1}\separgs{#2}\allargs}
\mathcommand\leqpp[2]{\leqsymbol\beginargs{#1}\separgs{#2}\allargs}
\mathcommand\lespp[2]{\lessymbol\beginargs{#1}\separgs{#2}\allargs}
\mathcommand\lexlespp[2]{\lexlessymbol\beginargs{#1}\separgs{#2}\allargs}
\mathcommand\lexlimlespp[3]
               {\lexlimlessymbol\beginargs{#1}\separgs{#2}\separgs{#3}\allargs}
\mathcommand\lexpp[3]{\lexsymbol\beginargs{#1}\separgs{#2}\separgs{#3}\allargs}
\mathcommand\ackpp[2]{\acksymbol\beginargs{#1}\separgs{#2}\allargs}
\mathcommand\switchpp[1]{\switchsymbol\beginargs{#1}\allargs}
\mathcommand\swatchpp[1]{\swatchsymbol\beginargs{#1}\allargs}
\mathcommand\whilepp[2]{\whilesymbol\beginargs{#1}\separgs{#2}\allargs}
\mathcommand\nullpp[1]{\nullsymbol\beginargs{#1}\allargs}
\mathcommand\nullppiterated[2]{\nullsymbol^{#1}\beginargs{#2}\allargs}
\mathcommand\hdpp[1]{\hdsymbol\beginargs{#1}\allargs}
\mathcommand\hdppiterated[2]{\hdsymbol^{#1}\beginargs{#2}\allargs}
\mathcommand\carpp[1]{\carsymbol\beginargs{#1}\allargs}
\mathcommand\cdrpp[1]{\cdrsymbol\beginargs{#1}\allargs}
\mathcommand\tlpp[1]{\tlsymbol\beginargs{#1}\allargs}
\mathcommand\tlppiterated[2]{\tlsymbol^{#1}\beginargs{#2}\allargs}
\mathcommand\inpp[2]{\insymbol\beginargs{#1}\separgs{#2}\allargs}
\mathcommand\inppiterated[3]{\insymbol^{#1}\beginargs{#2}\separgs{#3}\allargs}
\mathcommand\applypp[2]{\applysymbol\beginargs{#1}\separgs{#2}\allargs}
\mathcommand\applyppiterated
[3]{\applysymbol^{#1}\beginargs{#2}\separgs{#3}\allargs}
\mathcommand\termpp[2]{\termsymbol\beginargs{#1}\separgs{#2}\allargs}
\mathcommand\setpp[1]{\set\beginargs{#1}\allargs}
\mathcommand\russellpp[1]{\russellsymbol\beginargs{#1}\allargs}

\mathcommand\Tpp[6]{\turingmachinesymbol\beginargs{#1}\separgs{#2}\separgs
{#3}\separgs{#4}\separgs{#5}\separgs{#6}\allargs}
\mathcommand\Tppseven[7]{\turingmachinesymbol\beginargs{#1}\separgs{#2}\separgs
{#3}\separgs{#4}\separgs{#5}\separgs{#6}\separgs{#7}\allargs}
\mathcommand\foreverppp[6]{\ident{foreverp}\beginargs{#1}\separgs{#2}\separgs
{#3}\separgs{#4}\separgs{#5}\separgs{#6}\allargs}
\mathcommand\terminatesppp[6]{\terminatespsymbol\beginargs{#1}\separgs
{#2}\separgs{#3}\separgs{#4}\separgs{#5}\separgs{#6}\allargs}
\mathcommand\terminatespppone[1]{\terminatespsymbol \beginargs{#1}\allargs}
\mathcommand\statepp
[3]{\statesymbol\beginargs{#1}\separgs{#2}\separgs{#3}\allargs}
\mathcommand\tightstatepp
[3]{\statesymbol\beginargs{#1}\tightsepargs{#2}\tightsepargs{#3}\allargs}
\mathcommand\cmdpp
[3]{\cmdsymbol  \beginargs{#1}\separgs{#2}\separgs{#3}\allargs}
\mathcommand\tightcmdpp
[3]{\cmdsymbol  \beginargs{#1}\tightsepargs{#2}\tightsepargs{#3}\allargs}
\mathcommand\stoppp           {\ident{stop}}
\mathcommand\leftpp           {\ident{left}}
\mathcommand\rightpp          {\ident{right}}
\mathcommand\nthpp         [2]{\nthsymbol  \beginargs{#1}\separgs{#2}\allargs}
\mathcommand\pppp          [1]{\ppsymbol\beginargs{#1}            \allargs}
\mathcommand\qppp          [2]{\qpsymbol\beginargs{#1}\separgs{#2}\allargs}
\mathcommand\Eppp          [1]{\Epsymbol\beginargs{#1}            \allargs}
\mathcommand\Epppzwei      [2]{\Epsymbol\beginargs{#1}\separgs{#2}\allargs}
\mathcommand\Pppp          [1]{\Ppsymbol\beginargs{#1}            \allargs}
\mathcommand\Ppppeinsindex [2]{\Ppsymbol_{#1}\beginargs{#2}\allargs}
\mathcommand\Ppppdrei      
[3]{\Ppsymbol\beginargs{#1}\separgs{#2}\separgs{#3}\allargs}
\mathcommand\Ppppvier
[4]{\Ppsymbol\beginargs{#1}\separgs{#2}\separgs{#3}\separgs{#4}\allargs}
\mathcommand\Qppp          [2]{\Qpsymbol\beginargs{#1}\separgs{#2}\allargs}
\mathcommand\Qpppeins      [1]{\Qpsymbol\beginargs{#1}\allargs}
\mathcommand\Qpppeinsindex [2]{\Qpsymbol_{#1}\beginargs{#2}\allargs}
\mathcommand\Qpppdrei      
[3]{\Qpsymbol\beginargs{#1}\separgs{#2}\separgs{#3}\allargs}
\mathcommand\Fatherpp      [2]{\Fathersymbol\beginargs{#1}\separgs{#2}\allargs}
\mathcommand\Marriespp     [2]{\Marriessymbol\beginargs{#1}\separgs{#2}\allargs}
\mathcommand\Lovespp       [2]{\Lovessymbol\beginargs{#1}\separgs{#2}\allargs}
\mathcommand\StolenBypp    [2]
{\StolenBysymbol\beginargs{#1}\separgs{#2}\allargs}
\mathcommand\Humanpp       [1]{\Humansymbol\beginargs{#1}\allargs}
\mathcommand\Evenpp        [1]{\Evensymbol\beginargs{#1}\allargs}
\mathcommand\Evenppi       [2]{\Evensymbol^{#1}\beginargs{#2}\allargs}
\mathcommand\Oddpp         [1]{\Oddsymbol\beginargs{#1}\allargs}
\mathcommand\Primepp       [1]{\Primesymbol\beginargs{#1}\allargs}
\mathcommand\EveryPairpp  [2]{\EveryPairsymbol\beginargs{#1}\separgs
{#2}\allargs}
\mathcommand\mindexppeins  [2]{\mindexsymbol{#1}\beginargs{#2}\allargs}
\mathcommand\Givepp        [3]{\Givesymbol
\beginargs{#1}\separgs{#2}\separgs{#3}\allargs}
\mathcommand\mindexppzwei  [3]{\mindexsymbol
{#1}\beginargs{#2}\separgs{#3}\allargs}
\mathcommand\mindexppdrei  [4]{\mindexsymbol
{#1}\beginargs{#2}\separgs{#3}\separgs{#4}\allargs}

\mathcommand\nonnegppp     [1]{\nonnegpsymbol\beginargs{#1}\allargs}

\mathcommand\anonymouscsymbol{c}
\mathcommand\anonymouscindexsymbol[1]{\anonymouscsymbol_{#1}}
\mathcommand\anonymousfsymbol{f}
\mathcommand\anonymouscpp
[2]{\anonymouscsymbol\beginargs{#1}\separgs\ldots\separgs{#2}\allargs}
\mathcommand\anonymouscindexpp
[3]{\anonymouscindexsymbol{#1}\beginargs{#2}\separgs\ldots\separgs{#3}\allargs}
\mathcommand\anonymousfpp
[2]{\anonymousfsymbol\beginargs{#1}\separgs\ldots\separgs{#2}\allargs}
\mathcommand\coerceindexpp[3]{[#3]_{#1}^{#2}}

\mathcommand\Elephantppp    [1]{\Elephantpsymbol\beginargs{#1}\allargs}
\mathcommand\Flowerppp      [1]{\Flowerpsymbol  \beginargs{#1}\allargs}
\mathcommand\Bicycleppp     [1]{\Bicyclepsymbol \beginargs{#1}\allargs}
\mathcommand\Germanppp      [1]{\Germanpsymbol  \beginargs{#1}\allargs}
\mathcommand\Hugeppp        [1]{\Hugepsymbol    \beginargs{#1}\allargs}
\mathcommand\Animalppp      [1]{\Animalpsymbol  \beginargs{#1}\allargs}
\mathcommand\Maleppp        [1]{\Malepsymbol    \beginargs{#1}\allargs}
\mathcommand\Boyppp         [1]{\Boypsymbol     \beginargs{#1}\allargs}
\mathcommand\Girlppp        [1]{\Girlpsymbol    \beginargs{#1}\allargs}
\mathcommand\Femaleppp      [1]{\Femalepsymbol  \beginargs{#1}\allargs}
\mathcommand\Roundppp       [1]{\Roundpsymbol   \beginargs{#1}\allargs}
\mathcommand\Bishoppp       [1]{\Bishopsymbol   \beginargs{#1}\allargs}
\mathcommand\Quadrangularppp[1]{\Quadrangularpsymbol  \beginargs{#1}\allargs}
\mathcommand\Kissedppp[2]{\Kissedpsymbol\beginargs{#1}\separgs{#2}\allargs}
\mathcommand\Metppp[2]   {\Metpsymbol   \beginargs{#1}\separgs{#2}\allargs}

\newcommand\bound     {{\rm bound}}
\newcommand\free      {{\rm free}}

\mathcommand\Vtripleindex[3]{\V\!_{{#1},\,{#2},\,{#3}}}
\mathcommand\Vdoubleindex[2]{\V\!_{{#1},\,{#2}}}
\mathcommand\Vsingleindex[1]{\V\!_{{#1}}}

\mathcommand\Erel[1]{\Gammaoffont\!_{#1}}
\mathcommand\Urel[1]{\Deltaoffont_{#1}}

\mathcommand\theRprimefromstrongtoweak{
  \inparenthesesinlinetight{
     \domres\id{\Vwall\cup\Vsome\setminus\RAN\varsigma}
     \nottight{\nottight\uplus}
     \reverserelation\varsigma
  }
  \nottight{\circ}
  \ranres
    {\transclosureinline R}
    {\Vwall\cup\Vsome\setminus\RAN\varsigma}
  \nottight{\nottight{\nottight{\uplus}}}
  \Vsome\tighttimes\Vsall
}

\mathcommand\deltaminus{\delta^-}
\mathcommand\deltaplus{\delta^+}
\mathcommand\deltaplusplus{\delta^{+^+}}
\mathcommand\deltastar{\delta^*}
\mathcommand\deltastarstar{\delta^{*^*}}

\mathcommand\Vall     {\Vsingleindex\indexdelta         }
\mathcommand\Vwall    {\Vsingleindex\indexdeltaminu     }
\mathcommand\Vsall    {\Vsingleindex\indexdeltaplus     }
\mathcommand\Vgsome   {\Vsingleindex\indexgammaplus     }
\mathcommand\Vsome    {\Vsingleindex\indexgamma         }
\mathcommand\Vfree    {\Vsingleindex\indexfree          }
\mathcommand\Vbound   {\Vsingleindex\indexbound         }
\mathcommand\Vsomesall{\Vsingleindex\indexgammadeltaplus}

\mathapplycommand\VARall      {\VARsingleindex\indexdelta         }
\mathapplycommand\VARwall     {\VARsingleindex\indexdeltaminu     }
\mathapplycommand\VARsall     {\VARsingleindex\indexdeltaplus     }
\mathapplycommand\VARgsome    {\VARsingleindex\indexgammaplus     }
\mathapplycommand\VARsome     {\VARsingleindex\indexgamma         }
\mathapplycommand\VARfree     {\VARsingleindex\indexfree          }
\mathapplycommand\VARbound    {\VARsingleindex\indexbound         }
\mathapplycommand\VARsomesall {\VARsingleindex\indexgammadeltaplus}
\mathcommand\displayVARsall[1]{\VARsingleindex\indexdeltaplus
\!\!\!\:\left(\begin{array}{@{}c@{}}#1\end{array}\right)}

\mathcommand\rigidvari     [2]{#1_{#2}^\indexgammadeltaplus}
\mathcommand\existsvari    [2]{#1_{#2}^\indexgamma    }
\mathcommand\forallvari    [2]{#1_{#2}^\indexdelta    }
\mathcommand\freevari      [2]{#1_{#2}^\indexfree     }
\mathcommand\wforallvari   [2]{#1_{#2}^\indexdeltaminu}
\mathcommand\sforallvari   [2]{#1_{#2}^\indexdeltaplus}
\mathcommand\gexistsvari   [2]{#1_{#2}^\indexgammaplus}
\mathcommand\boundvari     [2]{#1_{#2}}
\mathcommand\vari          [2]{#1_{#2}}
\mathcommand\wforallvarilow[2]{#1_{#2}^
{\raisebox{-.82ex}{\math\indexdeltaminu}}}

\newcommand\indexhelper[1]{{\scriptscriptstyle#1\:\!\!}}
\newcommand\indexdeltaplus
{\indexhelper{\delta^{\raisebox{-.17ex}{\fvesf\hskip-0.14em +}}}}
\newcommand\indexdeltaminu
{\indexhelper{\delta^{\mbox{\fvesf\hskip-0.14em\rule[.2ex]{.7em}{.15ex}}}}}
\newcommand\indexgammaplus
{\indexhelper{\gamma^{\mbox{\fvesf\hskip-0.14em +}}}}
\newcommand\indexgammadeltaplus
{\indexhelper{\gamma\delta^{\raisebox{-.17ex}{\fvesf\hskip-0.14em +}}}}

\newcommand\indexdelta{\indexhelper\delta}
\newcommand\indexgamma{\indexhelper\gamma}
\newcommand\indexfree
{{\scriptscriptstyle\free}}
\newcommand\indexbound
{{\scriptscriptstyle\bound}}

\newcommand\Wellfsymb{\ident{Wellf}}
\mathapplycommand\Wellfpp{\Wellfsymb}

%% file: headersugarterms.tex
\mathcommand\beginargs{(}
\mathcommand\allargs  {)}
\mathcommand\separgs  {,\,}
\mathcommand\tightsepargs{,}

\mathcommand\minusppnoparentheses  [2]{{#1}\,\minussymbol\,{#2}}
\mathcommand\tightminusppnoparentheses  [2]{{#1}\minussymbol{#2}}
\mathcommand\divideppnoparentheses [2]{{#1}\,\dividesymbol\,{#2}}
\mathcommand\plusppnoparentheses   [2]{{#1}\,\plussymbol \,{#2}}
\mathcommand\plusppnoparenthesesi  [3]{{#2}\,\plussymbol^{#1}\,{#3}}
\mathcommand\tightplusppnoparentheses   [2]{{#1}\plussymbol{#2}}
\mathcommand\timesppnoparentheses  [2]{{#1}\,\timessymbol\,{#2}}
\mathcommand\undppnoparentheses    [2]{{#1}\und            {#2}}
\mathcommand\oderppnoparentheses   [2]{{#1}\oder           {#2}}
\mathcommand\impliesppnoparentheses[2]{{#1}\implies        {#2}}
\mathcommand\leqinfixppnoparentheses[2]{{#1}\,\tight\leq\,{#2}}
\mathcommand\geqinfixppnoparentheses[2]{{#1}\,\tight\geq\,{#2}}
\mathcommand\dividepp [2]{(\divideppnoparentheses {#1}{#2})}
\mathcommand\minuspp  [2]{(\minusppnoparentheses  {#1}{#2})}
\mathcommand\pluspp   [2]{(\plusppnoparentheses   {#1}{#2})}
\mathcommand\timespp  [2]{(\timesppnoparentheses  {#1}{#2})}
\mathcommand\undpp    [2]{(\undppnoparentheses    {#1}{#2})}
\mathcommand\oderpp   [2]{(\oderppnoparentheses   {#1}{#2})}
\mathcommand\impliespp[2]{(\impliesppnoparentheses{#1}{#2})}

%% file: header12pt.tex
\mathcommand\notconflu{\mathchoice
             {{\hskip1.5pt\nmid\hskip-4.697545pt\downarrow}}
             {{\hskip1.5pt\nmid\hskip-4.65pt\downarrow}}   
             {{\hskip1pt\nmid\hskip-3.494pt\downarrow\hskip1pt}}  
             {{\hskip1pt\nmid\hskip-3.01pt\downarrow\hskip0.5pt}}   
}
\mathcommand\redpara{\mathchoice
           {{\redsimple\hskip-16pt  \shortparallel}\hskip8.5pt}
           {{\redsimple\hskip-16pt  \shortparallel}\hskip8.5pt}
           {{\redsimple\hskip-8.5pt \shortparallel}\hskip6pt}
           {{\redsimple\hskip-7.5pt \shortparallel}\hskip5pt}
}
\mathcommand\antiredpara{\mathchoice
           {{\antired\hskip-14.6pt  \shortparallel}\hskip7pt}
           {{\antired\hskip-14.6pt  \shortparallel}\hskip7pt}
           {{\antired\hskip-8.pt \shortparallel}\hskip5pt}
           {{\antired\hskip-7.pt \shortparallel}\hskip5pt}
}
\mathcommand\revpara{\mathchoice
           {{\redsimple\hskip-16pt  \infty}\hskip4.8pt}
           {{\redsimple\hskip-16pt  \infty}\hskip4.8pt}
           {{\redsimple\hskip-11.5pt\infty}\hskip4pt}
           {{\redsimple\hskip-9.9pt \infty}\hskip3pt}
}
\mathcommand\antirevpara{\mathchoice
           {{\antired\hskip-15.4pt\infty}\hskip4pt}
           {{\antired\hskip-15.4pt\infty}\hskip4pt}
           {{\antired\hskip-10.8pt\infty}\hskip3pt}
           {{\antired\hskip-9.5pt \infty}\hskip3pt}
}
\mathcommand\simpara{\mathchoice
           {{\redsimple\hskip-13pt  \circ}\hskip7pt}
           {{\redsimple\hskip-13pt  \circ}\hskip7pt}
           {{\redsimple\hskip-11.5pt\circ}\hskip4pt}
           {{\redsimple\hskip-9.9pt \circ}\hskip3pt}
}

%% file: headerproof.tex
\newcommand\Proofof{Proof of}
\renewenvironment{proof}[1]{\yestop\begin{sloppypar}\noindent
{\bf\Proofof\ {#1}}\\}{\end{sloppypar}}
\renewenvironment{proofqed}[1]
{\begin{sloppypar}\def\fooqed{#1}\noindent{\bf\Proofof\ \fooqed}}
{\QEDbf\fooqed\end{sloppypar}}
\renewenvironment{proofparsep}[1]{\parindent=0pt\begin
{sloppypar}\noindent{\bf\Proofof\ {#1}}\nopagebreak\par}
{\end{sloppypar}}
\renewenvironment{proofparsepqed}[1]{\parindent=0pt\begin
{sloppypar}\def\fooqed{#1}\noindent{\bf\Proofof\ \fooqed}\nopagebreak\par}
{\nopagebreak\QEDbf\fooqed\end{sloppypar}}
\renewenvironment{proofshort}[1]{\begin{sloppypar}\noindent
{\bf\Proofof\ {#1}:}}{\end{sloppypar}}

%% file: quotation.tex
\newcommand\englishtextelevenone
{\englishtexteleventwo{Not only the notorious golden mountain is of gold, 
 but also\\the}.}
\newcommand\englishtexteleventwo[1]
{#1 round quadrangle is just as certainly round as it is quadrangular}

\newcommand\englishtextonehundredandthirteenquotation
{our translation}

\newcommand\englishtextninehundred
{Thus, in mathematics, we have no reasons to assume any meaning of
 `existence' that would be fundamentally different from that of 
 `the validity of axiomatic relations'.}

%% file: hbdictionary.tex
\newcommand\englishabhandlung             
{paper}

\newcommand\germanabstrakt                {ab\esi trakt}

\newcommand\englishabstrahierenvon        
                                          {leave out of account}

\newcommand\englishallgemeinuniversal     {universal}

\newcommand\englisheinallgemeinuniversal  {a \englishallgemeinuniversal}

\newcommand\englishallgemeinesurteilindex {\index{universal judgment}}
\newcommand\englishallgemeinesurteil
      {\englishallgemeinesurteilindex\englishallgemeinuniversal\ \englishurteil}

\newcommand\englishallgemeineurteile      {\englishallgemeinesurteil s}

\newcommand\englishdasallgemeineunddasexistentialeurteil
                {\englishallgemeinesurteilindex\englishexistentialesurteilindex
                               the universal and the existential \englishurteil}

\newcommand\englishAllgemeineundexistentialeurteile
                {\englishallgemeinesurteilindex\englishexistentialesurteilindex
                                      Universal and existential \englishurteile}

\newcommand\englishallgemeingueltigkeit 
{\index{validity!universal}universal validity}

\newcommand\englishallgemeingueltig     
{\index{validity!universal}universally valid}

\newcommand\englishalternativedichotomy   
{alternative}
\newcommand\englisheinealternativedichotomy
{an alternative}
\newcommand\englishalternativendichotomies
{alternatives}

\newcommand\germananalysis                {Ana\-lysi\es}
\newcommand\englishanalysis               {Analysis}

\newcommand\germananfangsgruende          {Anfang\esi gr\ue nde}

\newcommand\englishangebbar               
                                          {explicitly specifiable}

\newcommand\englishangeben
                                          {explicitly specify}

\newcommand\germanannahme                 {Annahme}
\newcommand\englishannahme                {assumption}
\newcommand\englishannahmen               {assumptions}
\newcommand\englishAnnahmen               {Assumptions}

\newcommand\germananschaulich             {an\-schau\-lich}
\newcommand\germananschaulichvorgestellt  {\germananschaulich\ vor\-gestellt}
\newcommand\englishanschaulich            {intuitive}

\newcommand\englishanschaulichvorgestellt {intuitively conceived}

\newcommand\englishanwendung             {application}

\newcommand\germananzahl                {Anzahl}

\newcommand\englishanzahldrei             {number}

\newcommand\englishanzahlenlehre          
{theory of cardinal numbers}

\newcommand\englishartsort               {sort}
\newcommand\englishartensorts            {\englishartsort s}

\newcommand\germanaufbau                 {Au\fbi au}

\newcommand\englishauffassenpointofview   {take the point of view}              
\newcommand\englishETWASalsETWASauffassen[2]
{\englishauffassenpointofview\ that #1 is #2}

\newcommand\englishaufloesung             {\index{resolution}resolution}

\newcommand\germanausdruck                {Au\esi druck}

\newcommand\englishausdruckterm           
{expression}
\newcommand\englishausdruck               {expression}

\newcommand\englishausdruckformel         
{expression}

\newcommand\englishausdruckexpression     {expression}

\newcommand\englishausdrucksweise
                                          {mode of \englishausdruckexpression}

\newcommand\englishausgezeichnetcanonical
{full}
\newcommand\englishAusgezeichnetcanonical
{Full}

\newcommand\germanaussage                 {Au\esi sage}

\newcommand\englishaussage                {proposition}
\newcommand\englisheineaussage            {a proposition}

\newcommand\propositionalcalculusindex    {\index{propostional calculus}}

\newcommand\englishaussagenkalkulnoindex  {propositional \englishkalkul}
\newcommand\englishaussagenkalkul         {\propositionalcalculusindex 
                                           \englishaussagenkalkulnoindex}
\newcommand\englishAussagenKalkul         {\propositionalcalculusindex 
                                           Propositional \englishKalkul}

\newcommand\germanaussagenlogik           {Au\esi sagen\-logik}
\newcommand\englishaussagenlogik          {propositional logic}

\newcommand\englishaussagenverbindung     
{propositional \englishverbindung}
\newcommand\englishaussagenverbindungen   
{propositional \englishverbindungen}

\newcommand\englishaussagenverknuepfung   
{propositional \englishverknuepfung}
\newcommand\englishaussagenverknuepfungen 
{propositional \englishverknuepfungen}

\newcommand\germanaxiomatik               {Axiomatik}
\newcommand\englishaxiomatik              {axiomatics}

\newcommand\englishbegriffnotion          {notion}

\newcommand\englishbegriffenotion         {\englishbegriffnotion s}

\newcommand\germanbegriffsbildung         {Begriff\esi bildung}
\newcommand\germanbegriffsbildungen       {Begriff\esi bildungen}
\newcommand\englishbegriffsbildung
                     {\index{concept formation}concept \englishbildungformation}

\newcommand\englishbegriffsbildungen      {\englishbegriffsbildung s}

\newcommand\englishbegruendungtheorie     {grounding}

\newcommand\germanbehauptung              {Be\-hauptung}
\newcommand\englishbehauptung             {assertion}

\newcommand\englishbehaupten              {assert}
\newcommand\englishbehauptenthirdsingular {\englishbehaupten s}

\newcommand\germanbeliebig                {beliebig}
\newcommand\englishbeliebig               {arbitrary}

\newcommand\englishbereich                {domain}

\newcommand\englishbereiche               {domains}
\newcommand\englishBereiche               {Domains}

\newcommand\englishbeschaffenheit
                                          {character}

\newcommand\englishbestandteil            
{part}

\newcommand\englishbestimmt               {definite}
\newcommand\englisheinbestimmterdefinite  {a definite}

\newcommand\englishbestimmtparticular     {particular}

\newcommand\englisheinbestimmterparticular{a \englishbestimmtparticular}

\newcommand\englishnaeherbestimmt      {more specifically \englishbestimmenppp}
\newcommand\englisheinnaeherbestimmter    {a \englishnaeherbestimmt}

\newcommand\englishbestimmenppp           {determined}

\newcommand\germanbetrachtung             {Be\-trach\-tung}
\newcommand\germanbetrachtungen           {Be\-trach\-tun\-gen}

\newcommand\englishbetrachtungeroerterung {discussion}
\newcommand\englishbetrachteneroertern    {discuss}

\newcommand\englishbetrachtungenueberlegungen{considerations}
\newcommand\englishbetrachtungbehandlung  {treatment}

\newcommand\englishbetrachtungbehandlungplusgenitive
                                          {\englishbetrachtungbehandlung\ of}

\newcommand\englishbetrachtungproblem     {problem}
\newcommand\englishbetrachtungenproblem   {\englishbetrachtungproblem s}

\newcommand\germanbetrachtungsweise       {Betrachtung\esi weise}
\newcommand\englishbetrachtungsweise      {\englishweise\ of consideration}

\newcommand\germanbeweistheorie           {Bewei\esi theorie}
\newcommand\englishbeweistheorie          {\index{proof theory}proof theory}
\newcommand\germanbeweistheoretisch       {bewei\esi theore\-tisch}
\newcommand\englishbeweistheoretisch      {\index{proof theory}proof-theoretic}

\newcommand\englishbezeichnen             {designate}

\newcommand\englishbezeichneningform      {designating}

\newcommand\germanbeziehung               {Be\-zie\-hung}

\newcommand\englishbeziehung              {relation}

\newcommand\englisheinebeziehung          {a relation}
\newcommand\englishbeziehungen            {relations}

\newcommand\germanbildung                 {Bildung}

\newcommand\englishbildungformation       {formation}

\newcommand\englishbildungsprozesse       
{formation processes}

\newcommand\englishbildungsregeln         
{formation rules}

\newcommand\englishdarstellendekonjunktion
{\index{Konjunktion!darstellende}%
 \index{conjunction, representing}%
 representing conjunction}

\newcommand\englishderjenigewelcher       {that one, which \ldots}

\newcommand\englishdual                   {\index{duality}dual}

\newcommand\englishsichdualgegenueberstehend
{\englishdual\ to each other}

\newcommand\englishdurchfuehrungimplementation
{implementation}

\newcommand\englisheszurdurchfuehrungbringen{put it into practice}

\newcommand\englishdurchgehen
                                          {go through}

\newcommand\englisheindeutigbestimmt
{\index{uniqueness}unique}

\newcommand\englisheindeutigfestgelegt    {\index{uniqueness}uniquely defined}

\newcommand\englisheingeschachtelt        {nested}

\newcommand\englisheineendformel          {an end formula}

\newcommand\germanendformel               {End\-formel}

\newcommand\englishentscheidungsproblemnoindex{decision \englishproblem}
\newcommand\englishentscheidungsproblem   
                   {\index{decision problem}\englishentscheidungsproblemnoindex}
\newcommand\englishEntscheidungsProblem   
                              {\index{decision problem}Decision \englishProblem}
\newcommand\englishentscheidungsprobleme  
{\englishentscheidungsproblem s}

\newcommand\germanentscheidungsverfahren  {Ent\-schei\-dung\esi ver\-fahren}
\newcommand\englishentscheidungsverfahren 
                                  {\index{decision procedure}decision procedure}
\newcommand\englishEntscheidungsverfahren 
                                  {\index{decision procedure}Decision procedure}

\newcommand\englishentscheidungsverfahrenplural 
                                          {\englishentscheidungsverfahren s}
\newcommand\englishEntscheidungsverfahrenplural 
                                          {\englishEntscheidungsverfahren s}

\newcommand\englisherfahrungskomplexe     
{experience-complexes}

\newcommand\englishunerfuellbar       {\index{satisfiable}un\-satis\-fi\-able}

\newcommand\englisherfuellbarkeit         {\index{satisfiable}satisfiability}

\newcommand\englisherwaegung              {consideration}

\newcommand\germanexistentialeurteil      {exi\esi tentiale Urteil}

\newcommand\germanexistentialenurteils    {exi\esi tentialen Urteil\es}
\newcommand\englishexistentialesurteilindex{\index{existential judgment}}
\newcommand\englishexistentialesurteil
                         {\englishexistentialesurteilindex existential judgment}
\newcommand\englisheinexistentialesurteil {an \englishexistentialesurteil}

\newcommand\englishexistentialformeln     
                              {\index{existential formula}existential \formulae}

\newcommand\englishexistenzsatz            
{existence sentence}

\newcommand\englishfertigegesamtheit      {completed totality}

\newcommand\englishfeststellungstatement  {statement}

\newcommand\englishfeststellungdurchnachweis{verification}

\newcommand\englishfestumgrenzt
{well-determined and fixed}

\newcommand\englishfigur                  {figure}

\newcommand\finitistindex                 {\index{finitist}}
\newcommand\finitismindex                 {\finitistindex\index{finitism}}

\newcommand\germanfinit                   {\finitistindex finit}
\newcommand\englishfinitnoindex                  
                                          {finitist}
\newcommand\englishfinit                  {\finitistindex\englishfinitnoindex}
\newcommand\englishFinit                  {\finitistindex Finitist}
\newcommand\englisheinfinit               {a \englishfinit}
\newcommand\englishEinfinit               {A \nolinebreak\englishfinit}

\newcommand\englishfinitismus             {\finitismindex finitism}

\newcommand\englishfolgerungconclusion    {conclusion}

\newcommand\germanformalismus             {Formali\esi mu\es}

\newcommand\englishformelvariable         {formula variable}

\newcommand\englishformelvariablen        {\englishformelvariable s}

\newcommand\englishformelsprache          
{formula language}

\mathcommand\theFormela{(x)\,A(x)\nottight\implies A(a)}

\mathcommand\theFormelinpiteins{a=b\nottight{\nottight{\implies}}
\inparenthesestight{a=c\nottight{\nottight{\implies}}b=c}}

\mathcommand\theFormelinpitzwei{a\tightequal b\nottight\implies b\tightequal a}

\mathcommand\theFormelinpitzweia
                      {a\tightnotequal b\nottight\implies b\tightnotequal a}
\mathcommand\theFormelinpitdrei{a\tightequal b\nottight\implies 
\inparenthesestight{b\tightequal c\nottight\implies a\tightequal c}}

\newcommand\germanformulierbarkeit          {Formulier\-bar\-keit}
\newcommand\englishformulierbarkeit         {expressibility}
\newcommand\germanformulierung              {Formulie\-rung}
\newcommand\englishformulierung             {expression}

\newcommand\englishfragequestion            {question}

\newcommand\englishfrageproblem           {problem}

\newcommand\englishfragenproblems         {\englishfrageproblem s}

\newcommand\germanfragestellung           {Fragestellung}
\newcommand\englishfragestellungdefinition{problem definition}
\newcommand\englishfragestellungendefinition{\englishfragestellungdefinition s}
\newcommand\englishFrageStellungdefinition{Problem Definition}

\newcommand\englishganzerationalefunktion 
{integer rational function}

\newcommand\englishgedankengang          
                                          {train of thought}

\newcommand\germangegenstand              {\index{Gegenstand}Gegen\-stand}
\newcommand\germangegenstaende            {Gegenst\ae nde}  

\newcommand\englishgegenstand             
{thing}

\newcommand\germangeltung                 {Gel\-tung}
\newcommand\englishgeltung                {\englishgueltigkeit}

\newcommand\englishzurgeltungbringenppp   {brought to bear}

\newcommand\englishzurgeltungkommenthirdsingular{comes to bear}

\newcommand\englishgenetisch              
{genetic}

\newcommand\englishgewuenscht             {desired}

\newcommand\englishgesamtheit             {totality}
\newcommand\englishgesamtheiten           {totalities}

\newcommand\englishGeschaeftsfuehrendeHerausgeber
{Editors-in-Chief}

\newcommand\germangesichtspunkt           {Gesicht\esi punkt}
\newcommand\englishgesichtspunkt          {viewpoint}

\mathcommand\theformelJeins               {a\tightequal a}

\mathcommand\theformelJzwei{a\tightequal b\nottight\implies
\inparenthesestight{A(a)\implies A(b)}}

\newcommand\englishgleichheitsbeziehung
{\index{equality!relation}equality relation}

\newcommand\englishgueltig                {\index{validity}valid}
\newcommand\germangueltigkeit             {G\ue ltig\-keit}
\newcommand\englishgueltigkeit            {\index{validity}validity}

\newcommand\englishhandhabung             {use}

\newcommand\englishmithilfevonexplicitly  {with the help of}

\newcommand\englishhinterglied            
{consequent}

\newcommand\englishhinzuziehendebetrachtung{consideration}

\newcommand\englishIndividuum             {Individual}
\newcommand\englishIndividuen             {\englishIndividuum s}

\newcommand\englishindividuensymbol
                                    {\index{individual symbol}individual symbol}
\newcommand\englisheinindividuensymbol    
                                          {an \englishindividuensymbol} 
 
\newcommand\englishindividuensymbole      
                                          {\englishindividuensymbol s}
\newcommand\englishpraedikatenundindividuensymbole
                                       {predicate and \englishindividuensymbole}
\newcommand\englishpraedikatenoderindividuensymbole
                                       {predicate or \englishindividuensymbole}
\newcommand\englishwederpraedikatennochindividuensymbole
                               {neither predicate nor \englishindividuensymbole}

\newcommand\englishindividuenvariable  
                                {\index{individual variable}individual variable}

\newcommand\englishindividuenvariablen    {\englishindividuenvariable s} 

\newcommand\germaninhaltlich              {inhalt\-lich}
\newcommand\englishinhaltlich             {\index{contentual}contentual}

\newcommand\englishinterpretation         {interpretation}

\newcommand\intuitionismindex            {\index{intuitionism}}
\newcommand\germanintuitionismus       {\intuitionismindex Intuitioni\esi mu\es}
\newcommand\germanintuitionistisch       {\intuitionismindex intuitionistisch}
\newcommand\englishintuitionismus        {\intuitionismindex intuitionism}

\newcommand\englishintuitionistisch      {\intuitionismindex intuitionist}

\newcommand\englishirgendein              {an arbitrary}

\newcommand\englishirgendwelche           {arbitrary}

\newcommand\englishjedebeliebige
                                          {any}

\newcommand\englishjedwede                {any}

\newcommand\englishkalkul                 {cal\-cu\-lus}
\newcommand\englishKalkul                 {Cal\-cu\-lus}

\newcommand\englishkettenschlussindex     {\index{chain inference}}

\newcommand\englishkettenschluesse          
{\englishkettenschlussindex chain inferences}

\newcommand\englishleitgedanke
                                          {main idea}

\newcommand\germanlogistik                {Logi\esi tik}

\newcommand\englishlogistik
                     {\index{Logistik@\germanlogistik}formal mathematical logic}

\newcommand\englishmengentheoretisch      {set-theoretic}

\newcommand\englishmethodischmethodological{methodological}

\newcommand\englishmethodischthesemethods {these methods}

\newcommand\englishmethodischeeinstellung
{\englishmethodischmethodological\ attitude}
\newcommand\englishmethodischeanforderung[1]
{\englishmethodischmethodological#1 requirements}

\newcommand\germanmitteilung              {Mit\-teilung}

\newcommand\englishmitteilungnoindex      {communication}

\newcommand\englishzeichenpluralzurmitteilungnoindex
                          {\englishzeichenplural\ for \englishmitteilungnoindex}

\newcommand\germanmitteilbar              {mit\-teil\-bar}

\newcommand\englishmitteilbar             {communicable}

\newcommand\englishmuendungsstelle
{leaf}
\newcommand\englishmuendungsstellen
                                          {leaves}

\newcommand\englishnachpruefen           
{check}

\newcommand\englishNamenVerzeichnis      {Index of Persons}

\newcommand\englishnormaldisjunktion
{\index{normal disjunction}normal disjunction}

\newcommand\englishnormierung            
{standardization}

\newcommand\englishnumerieren            
                                         {number}

\newcommand\englishobjekt                {object}
\newcommand\englishobjekte               {\englishobjekt s}

\newcommand\englishparagraphsection       {section}

\newcommand\englishPraedikatenKalkul      {Predi\-cate \englishKalkul}
\newcommand\indexenglisheinstelligerpraedikatenkalkul
                                          {\index{predicate calculus!singulary}}

\newcommand\englishEinstelligerpraedikatenkalkulohnepraedikatenkalkul
                          {\indexenglisheinstelligerpraedikatenkalkul Singulary}

\newcommand\englishEinstelligerPraedikatenKalkul
{\englishEinstelligerpraedikatenkalkulohnepraedikatenkalkul\
                                                      \englishPraedikatenKalkul}

\newcommand\germanpraedikatenlogik        {Pr\ae dikaten\-logik}

\newcommand\englishpraedikatenlogik       {predicate logic}

\newcommand\germanpraemisse               {Pr\ae\-misse}
\newcommand\englishpraemisse              {\index{premise}premise}

\newcommand\englishpraemissen             {\englishpraemisse s}

\newcommand\englishproblem                {problem}
\newcommand\englishProblem                {Problem}

\newcommand\germansachverhalt             {Sach\-verhalt}
\newcommand\germansachlage                {Sach\-lage}

\newcommand\englishsachlage               {situation}
\newcommand\englishsachverhalt            {state of affairs}

\newcommand\englishSachVerzeichnis        {Subject Index}

\newcommand\englishsatztheorem            {sentence}
\newcommand\englishsaetzetheoreme         {sentences}

\newcommand\englishsatzproposition        {sentence}
\newcommand\englisheinsatzproposition     {a \englishsatzproposition}
\newcommand\englishsaetzepropositions     {\englishsatzproposition s}

\newcommand\englishsatzverbindung         
{\index{sentential combination}sentential \englishverbindung}
\newcommand\englishsatzverbindungen       
{\index{sentential combination}sentential \englishverbindungen}

\newcommand\germanschliessen              {Schlie\sz en}
\newcommand\englishschliessen             {inference}

\newcommand\germanschlussfolgerungoldspelling{Schlu\sz folgerung}
\newcommand\englishschlussfolgerung       {conclusion}
\newcommand\englishschlussfolgerungen     {conclusions}

\newcommand\englishschlussfigurindex
                                          {\index{proof figure}}
\newcommand\englishschlussfigur
                       {\englishschlussfigurindex\englishschluss\ \englishfigur}

\newcommand\germanschlussoldspelling      {Schlu\sz}
\newcommand\englishschluss                {inference}

\newcommand\englishschluesse              {\englishschluss s}

\newcommand\germanschlussweisen           {Schlu\sesi weisen}
\newcommand\germanschlussweisenoldspelling{Schlu\sz weisen}

\newcommand\englishschlussweisen          {\englishweisen\ of inference}
\newcommand\englishweise                  {mode}
\newcommand\englishweisen                 {\englishweise s}

\newcommand\englishseinszeichen           
{existential quantifier \englishzeichen}

\newcommand\englishsinnlichewahrnehmung   
{sense perception}

\newcommand\germansprachgebrauch           {Sprach\-gebrauch} 
\newcommand\englishsprachgebrauch         
{usage of language}
\newcommand\englishgewoehnlichersprachgebrauch
{common usage of language}

\newcommand\germanstandpunkt              {Stand\-punkt}
\newcommand\englishstandpunkt             {standpoint}
\newcommand\englisheinstandpunkt          {a \englishstandpunkt}

\newcommand\germansystem                  {Sy\esi tem}

\newcommand\englishtatsaechlichkeit       
{actuality}

\newcommand\germanueberblick              {\Ue ber\-blick}

\newcommand\englishueberfuehrbar          
{\index{convertibility}convertible}

\newcommand\englishueberfuehrbarkeit      
{\index{convertibility}convertibility}

\newcommand\englishueberlegung            {consideration}

\newcommand\englishueberlegungsingularasplural
                                          {considerations}

\newcommand\germanueberschreiten          {\ue ber\-schreiten}

\newcommand\germanueberschreitung         {\Ue ber\-schreitung}
\newcommand\englishueberschreitung        {transgression}
\newcommand\englisheineueberschreitung    {a \englishueberschreitung}
\newcommand\englishUeberschreitung        {Transgression}

\newcommand\englishueblich                {cus\-tom\-ary}

\newcommand\englishumfangslogikohnelogik  {extension}

\newcommand\englishumgekehrtorder         {reverse}

\newcommand\englishunmoeglichkeitsbeweis  {\index{impossibility proof}%
impossibility proof}

\newcommand\englisheinunmoeglichkeitsbeweis{an \englishunmoeglichkeitsbeweis}

\newcommand\germanurteil                  {Urteil}
\newcommand\englishurteil                 {judg\-ment}
\newcommand\englishurteile                {\englishurteil s}

\newcommand\englishverbindung             {combination}
\newcommand\englishverbindungen           {com\-bina\-tions}

\newcommand\germanverfahren               {Ver\-fahren}
\newcommand\englishverfahrenvorgehen      {approach}
\newcommand\englishverfahrenvorgehenplural{\englishverfahrenvorgehen es}

\newcommand\englishverfahrenprocedure     {procedure}

\newcommand\englishverfahrenprocedures    {\englishverfahrenprocedure s}

\newcommand\germanverhaeltnis             {Ver\-h\ae lt\-ni\es}
\newcommand\englishverhaeltnis            {relationship}

\newcommand\englishverteilungvonwahrheitswerten{distribution of truth values}

\newcommand\englishprepverteilungvonwahrheitswertenauf{on}

\newcommand\englishsichverifizieren    
{prove true}

\newcommand\englishverknuepfung           
{connective}
\newcommand\englishverknuepfungen         {\englishverknuepfung s}
\newcommand\englishverknuepfungprep       
{of}

\newcommand\englishverschaerfen           
{sharpen}

\newcommand\englishschaerfer  
{sharper}

\newcommand\germanverschaerft             {ver\-sch\ae rft}
\newcommand\englishverschaerft            
{sharp\-ened}

\newcommand\englishverschaerfung          
{sharpen\-ing}

\newcommand\germanvollstaendigeinduktion
  {\index{Induktion!vollstaendige@vollst\"andige}voll\-st\ae ndi\-ge Induktion}
\newcommand\germanvollstaendigeninduktion
  {\index{Induktion!vollstaendige@vollst\"andige}voll\-st\ae ndi\-gen Induktion}
\newcommand\englishvollstaendigeinduktion
                         {\index{induction!mathematical}mathematical induction}

\newcommand\germanvollstaendigkeit        {Voll\-st\ae ndig\-keit}

\newcommand\englishvollstaendigkeit       {\index{completeness}completeness}
\newcommand\englishVollstaendigkeit       {\index{completeness}Completeness}

\newcommand\germanfragedervollstaendigkeit{Frage der \germanvollstaendigkeit} 
\newcommand\englishfrageiSfragedervollstaendigkeit{problem}

\newcommand\englishfragedervollstaendigkeitlabel{\index
                    {completeness!{\englishfrageiSfragedervollstaendigkeit} of}}
\newcommand\englishfragedervollstaendigkeitnolabel
          {\englishfrageiSfragedervollstaendigkeit\ of \englishvollstaendigkeit}

\newcommand\englishfragedervollstaendigkeit
  {\englishfragedervollstaendigkeitlabel\englishfragedervollstaendigkeitnolabel}

\newcommand\englishvoraussetzen           {presuppose}

\newcommand\germanvoraussetzung           {Vorau\esi setzung}
\newcommand\englishvoraussetzung          {presupposition}
\newcommand\englishvoraussetzungen        {presuppositions}

\newcommand\englishvordergliednoindex
{antecedent}

\newcommand\germanvorliegend              {vorliegend}
\newcommand\englishjetztvorliegendpresent {present}
\newcommand\englishvorliegendgiven        {given}

\newcommand\englishvorlegenpppgiven       {given}

\newcommand\germanvorstellung             {Vor\-stellung}
\newcommand\englishvorstellung            {conception}

\newcommand\englishwertevorrat
{range of values}

\newcommand\englishwertsystemdervariablen 
{\index{valuation}valuation}

\newcommand\germanwertverteilungdervariablen
{Wertverteilung der Variablen}
\newcommand\englishwertverteilungdervariablen
{\index{\englishverteilungvonwahrheitswerten}%
 \index{\germanwertverteilungdervariablen}%
 \englishverteilungvonwahrheitswerten}
\newcommand\englishprepwertverteilungauf
                                   {\englishprepverteilungvonwahrheitswertenauf}

\newcommand\englishwertverteilungdervariablenexplizit
{\index{\englishverteilungvonwahrheitswerten}%
 \index{\germanwertverteilungdervariablen}%
 \englishverteilungvonwahrheitswerten\
 \englishprepwertverteilungauf\ the variables}

\newcommand\germanwiderspruch             {Wider\-spruch}
\newcommand\germanwidersprueche           {Wider\-spr\ue che}
\newcommand\englishwiderspruch            {contradiction}

\newcommand\englishwidersprueche          {\englishwiderspruch s}

\newcommand\englishwiderspruchsfrei       {consistent}
\newcommand\englisheinwiderspruchsfrei    {a \englishwiderspruchsfrei}
\newcommand\germanwiderspruchsfreiheit    {Wider\-spruch\esi freiheit}
\newcommand\englishwiderspruchsfreiheit   {consistency}
\newcommand\englishWiderspruchsfreiheit   {Consistency}

\newcommand\englishfrageiSfragederwiderspruchsfreiheiterfuellbarkeit{problem}
\newcommand\englishfrageniSfragenderwiderspruchsfreiheiterfuellbarkeit
                  {\englishfrageiSfragederwiderspruchsfreiheiterfuellbarkeit s}

\newcommand\englishfragederwiderspruchsfreiheitindex
                                                 {\index{consistency!problem of}}

\newcommand\englishfragenderwiderspruchsfreiheit
{\englishfragederwiderspruchsfreiheitindex
 \englishfrageniSfragenderwiderspruchsfreiheiterfuellbarkeit\ 
                                                 of \englishwiderspruchsfreiheit}

\newcommand\englishproblemderwiderspruchsfreiheitindex
                                                 {\index{consistency!problem of}}

\newcommand\englishProblemderWiderspruchsfreiheitnoindex
                                        {Problem of \englishWiderspruchsfreiheit}
\newcommand\englishProblemderWiderspruchsfreiheit
{\englishproblemderwiderspruchsfreiheitindex
                                   \englishProblemderWiderspruchsfreiheitnoindex}

\newcommand\germanzahl                    {Zahl}

\newcommand\germanzahlentheorie           {Zahlen\-theorie}
\newcommand\englishzahlentheorienoindex   {number theory}  
\newcommand\englishZahlenTheorienoindex   {Number Theory}  
\newcommand\englishzahlentheorie    
                              {\index{number theory}\englishzahlentheorienoindex}
\newcommand\englishZahlenTheorie    
                              {\index{number theory}\englishZahlenTheorienoindex}

\newcommand\englishzahligidentisch[1]     
{\math{#1}-identical}
     
\newcommand\englishzahligidentischeins[1] 
{\englishzahligidentisch 1}
\newcommand\englishzahligidentischzwei    
{\englishzahligidentisch 2}
\newcommand\englishzahligidentischdrei    
{\englishzahligidentisch 3}
\newcommand\englishzahligidentischvier    
{\englishzahligidentisch 4}

\newcommand\englishzeichen                {symbol}

\newcommand\englishzeichenplural          {symbols}

\newcommand\englishzenoscheparadoxie      {\achillesparadox}

\newcommand\germanziffer                  {Ziffer}
\newcommand\germanziffern                 {\germanziffer n}
\newcommand\englishziffer                 {numeral}

\newcommand\englishziffern                {numer\-als}

\newcommand\englishzulaessig              {admissible}

\newcommand\germanzunaechst               {zun\ae chst}

\newcommand\englishzunaechstfirst         {first}
\newcommand\englishZunaechstfirst         {First}

\newcommand\englishZunaechstatfirst       {At \nolinebreak first}

\newcommand\englishzurueckfuehreniSlogRed {reduce}
\newcommand\englishzurueckfuehreniSlogRedthirdsingular
                                          {\englishzurueckfuehreniSlogRed s}

\newcommand\englishzurueckfuehreniSlogRedppp{\englishzurueckfuehreniSlogRed d}

\newcommand\formulae                      {formulas}

\newcommand\hbsectionhelper[1]{{\large\par\noindent\LINEnomath{{\bf #1}}\par}}
\newcommand\hbsubsectionhelper[1]{{\par\noindent\LINEnomath{{\bf #1}}\par}}
\newcommand\hbsection[2]{\hbsectionhelper{\S\,\,#1.~~~#2}}
\newcommand\hbsubsection[2]{\hbsubsectionhelper{#1.~~~#2}}
\newcommand\thehbsection[2]{\hbsection{#1}{\csname hbsection#1#2\endcsname}}
\newcommand\thehbsubsection
[3]{\hbsubsection{#2}{\csname hbsubsection#1s#2#3\endcsname}}%
\newcommand\sethbsection[4]{%
\expandafter\newcommand\csname hbsection#1#2\endcsname{#3}%
\expandafter\newcommand\csname tochbsection#1#2\endcsname{#4}%
}
\newcommand\sethbsubsection[5]{%
\expandafter\newcommand\csname hbsubsection#1s#2#3\endcsname{#4}%
\expandafter\newcommand\csname tochbsubsection#1s#2#3\endcsname{#5}%
}
\newcommand\tochbsection[3]{\contentsline
 {section}{\numberline{\S\,\,#1.}{\csname tochbsection#1#2\endcsname}}{#3}{}}
\newcommand\tochbsectionnonumber[3]{\contentsline
 {section}{\csname tochbsection#1#2\endcsname}{#3}{}}
\newcommand\tochbsubsection[3] 
{\contentsline{subsection}{\numberline{(#1)}{#2}}{#3}{}}
\newcommand\tochbIIsubsection[4]
{\contentsline
 {subsection}{\numberline{#2.}{\csname tochbsubsection#1s#2#3\endcsname}}{#4}{}}

\sethbsection{1}{I}
{The \label
{problemderwiderspruchsfreiheit1}\englishProblemderWiderspruchsfreiheit\ 
in Axiomatics \nlnomath
~~~as a Logical \englishEntscheidungsProblem}
{The \englishProblemderWiderspruchsfreiheit\ in Axiomatics
as a \\ Logical \englishEntscheidungsProblem}

\sethbsection{2}{I}
{Elementary \englishZahlenTheorie. 
 \ --- \ 
 \englishFinit\
 Inference
 \nlnomath
 ~~~and its Limits}
{Elementary \englishZahlenTheorie. 
 \ --- \ 
 \englishFinit\
 Inference
 and its Limits \ }

\sethbsection{3}{I}
{Formalization of Logical Inference I: 
 \nlnomath
 ~~~The \englishAussagenKalkul\edfootnotemark{45I1}}
{Formalization of Logical Inference I: \ 
 The \englishAussagenKalkul}

\sethbsection{4}{I}
{Formalization of Logical Inference II: 
 \nlnomath
 ~~~The \englishPraedikatenKalkul\edfootnotemark{86I1}}
{Formalization of Logical Inference II: \
 The \englishPraedikatenKalkul}

\sethbsection{5}{I}
{Adding the\edfootnotemark{163I4}  Identity. \  
 \englishVollstaendigkeit\ of the
 \nlnomath
 ~~~\englishEinstelligerPraedikatenKalkul}
{Adding the Identity.\\
 \englishVollstaendigkeit\  of the \englishEinstelligerPraedikatenKalkul}

\sethbsection{6}{I}
{\englishWiderspruchsfreiheit\ 
of  
Infinite \englishBereiche\ of \englishIndividuen.
 \nlnomath
 ~~~Beginnings of \englishZahlenTheorie.}
{\englishWiderspruchsfreiheit\
of  
Infinite \englishBereiche\ of \englishIndividuen.
\\Beginnings of \englishZahlenTheorie}

\sethbsection{7}{I}
{The Recursive Definitions}
{The Recursive Definitions}

\sethbsection{8}{I}
{The Notion ``\englishderjenigewelcher'' and its Eliminability}
{The Notion ``\englishderjenigewelcher'' and its Eliminability}

\sethbsection{9}{I}
{\englishNamenVerzeichnis}
{\englishNamenVerzeichnis}

\sethbsection{10}{I}
{\englishSachVerzeichnis}
{\englishSachVerzeichnis}

\sethbsection{1}{II}
{The Method of Eliminating the Bound Variable
 \nlnomath 
 \englishmithilfevonexplicitly\
 \hilbert's \math\varepsilon-Symbol.}
{The Method of Eliminating the Bound Variables\\\englishmithilfevonexplicitly\
 \hilbert's \math\varepsilon-Symbol}

\sethbsubsection{1}{1}{II}
{The process of the symbolic 
 \label{resolution page}%
 \englishaufloesung\ of 
 \label{existential formula page}%
 \englishexistentialformeln.}
{The process of the symbolic \englishaufloesung\ of \englishexistentialformeln}

\sethbsubsection{1}{2}{II}
{\hilbertsepsilon-symbol and the \math\varepsilon-formula.}
{\hilbertsepsilon-symbol and the \math\varepsilon-formula}

%% file: VITA_GH3.tex
Erstellt~am:~16.08.1997
\\\mbox{}
\\\mbox{}
\\\mbox{}~~~~~~~~~~~~~~~~VITA~~~~Gisbert~HASENJAEGER,~geb.~01.06.19~Hildesheim
\\\mbox{}
\\\mbox{}~~~~~~~~~~[REM~mit~Korrekturen~im~Sinne~von~Standardisierung~und~von~K\ue rzungen
\\\mbox{}~~~~~~~~~~~~~~~einverstanden;~Mitteilung~des~Ergebnisses~erbeten]
\\\mbox{}
\\\mbox{}~~~~~~~~~~1937~Abitur~in~M\ue lheim~a.d.~Ruhr.1937-1939~Arbeits-~und~Wehr-dienst.
\\\mbox{}
\\\mbox{}~~~~~~~~~~1942~Nach~"Kopfschu\sz"~(2.1.)~und~Rekonvalensenz~bei~OKW/CHI~gegen
\\\mbox{}~~~~~~~~~~~~~~~(wie~ich~viel~sp\ae ter~erfuhr)~A.~Turing~eingesetzt.
\\\mbox{}
\\\mbox{}~~~~~~~~~~1945-1950~Studium~in~M\ue nster/Westf.,~dabei~seit~Beginn~(wissensch.)
\\\mbox{}~~~~~~~~~~~~~~~Hilfskraft~am~Seminar/Institut~f\ue r~Mathematische~Logik~und
\\\mbox{}~~~~~~~~~~~~~~~Grundlagenforschung.
\\\mbox{}
\\\mbox{}~~~~~~~~~~1950~Promotion~zum~Dr.rer.nat.~bei~Heinrich~Scholz~mit~der~Disser-
\\\mbox{}~~~~~~~~~~~~~~~tation~"Topologische~Untersuchungen~zur~Semantik~und~Syntax
\\\mbox{}~~~~~~~~~~~~~~~eines~erweiterten~Pr\ae dikatenkalk\ue ls".
\\\mbox{}~~~~~~~~~~1950~wissenschaflicher~Assistent~am~o.a.~Institut.
\\\mbox{}
\\\mbox{}~~~~~~~~~~1950/51~(das~akademische~Jahr)~Stipendiat~Gast~bei~Paul~Bernays~a.d.
\\\mbox{}~~~~~~~~~~~~~~~ETH~Z\ue rich.
\\\mbox{}
\\\mbox{}~~~~~~~~~~1953~Habilitation~und~venia~legendi~(Mathematische~Logik~und~Grund-
\\\mbox{}~~~~~~~~~~~~~~~lagenforschung)~an~der~Universit\ae t~M\ue nster~(Habilitationsschrift.
\\\mbox{}~~~~~~~~~~~~~~~Widerspruchsfreie~Axiomensysteme~ohne~Standard-Modell).
\\\mbox{}~~~~~~~~~~1955~Di\ae tendozentur,~1960~apl.~Professur~U~M\ue nster.
\\\mbox{}
\\\mbox{}~~~~~~~~~~1961~SS~und~WS~61/62~Vertretung~einer~neu~gegr\ue ndeten~Professur~f\ue r
\\\mbox{}~~~~~~~~~~~~~~~Logik~an~der~Philosophischen~Fakult\ae t~der~Universit\ae t~Bonn.
\\\mbox{}
\\\mbox{}~~~~~~~~~~1962~Berufung~auf~diese~(a.o.)~Professur,~f\ue r~Logik~und~Grundlagen-
\\\mbox{}~~~~~~~~~~~~~~~forschung.~[REM~Durch~das~Fehlen~des~Adjektivs~(math.)~und~den
\\\mbox{}~~~~~~~~~~~~~~~heutigen~Sprachgebrauch~ist~das~nun~wohl~"ein~weites~Feld",
\\\mbox{}~~~~~~~~~~~~~~~welches~zu~beackern~ich~mir~nicht~mehr~zumuten~m\oe chte.]
\\\mbox{}~~~~~~~~~~1964~Pers\oe nlicher~Ordinarius,~1966~o.~Prof.
\\\mbox{}
\\\mbox{}~~~~~~~~~~1964/5~(das~akademische~Jahr)~Gast~am~Institut~for~Advanced~Studies
\\\mbox{}~~~~~~~~~~~~~~~Princeton~N.J.
\\\mbox{}~~~~~~~~~~1970/1~(das~akademische~Jahr)~Gast-Professur~an~der~University~of
\\\mbox{}~~~~~~~~~~~~~~~Illinois~in~Urbana/Champaign.
\\\mbox{}
\\\mbox{}~~~~~~~~~~1984~Emeritierung~U~Bonn.
\\\mbox{}
\\\mbox{}\pagebreak\par